\documentclass[12pt,a4paper]{article}

\usepackage{a4wide,amssymb,amsmath,amsthm,xspace,epsfig}

\usepackage{epigraph, varwidth}

\usepackage{mathrsfs}

\usepackage{tikz-cd}
\usepackage{hyperref}

\newcommand\starbig{\hspace*{-.3ex}
\begin{tikzpicture}
\draw[thick, scale = 1.75](.5ex,-.1ex) -- (.5ex,.9ex);
\draw[thick, scale = 1.75](0,.4ex) -- (1ex,.4ex);
\draw[scale = 1.75](.2ex,.1ex) -- (.8ex,.7ex);
\draw[scale = 1.75](.2ex,.7ex) -- (.8ex,.1ex);
\end{tikzpicture}\hspace*{-.6ex}
}

\newcommand\fanfig{\hspace*{-.3ex}
\begin{tikzpicture}
\draw[thick, scale = 1.8](-.25ex,-.1ex) -- (1ex,-.1ex);
\draw[thick, scale = 1.8](.25ex,-.1ex) -- (.25ex,.65ex);
\draw[scale = 1.8](.25ex,-.1ex) -- (.57ex,.58ex);
\draw[scale = 1.8](.25ex,-.1ex) -- (.78ex,.43ex);
\draw[scale = 1.8](.25ex,-.1ex) -- (.93ex,.23ex);
\draw[scale = 1.8](.25ex,-.1ex) -- (.98ex,.07ex);
\end{tikzpicture}\hspace*{-.6ex}
}

\newcommand\fatS{\hspace*{-.3ex}
\begin{tikzpicture}
\draw[thick] (-.5ex,-.75ex) circle (.7ex);
\foreach \angle in {0,30,...,330}
    {\pgfpathcircle{\pgfpointadd{\pgfpoint{-.5ex}{-.75ex}}{\pgfpointpolar{\angle}{1ex}}}{.5pt}}
  \pgfusepath{fill}
\end{tikzpicture}
\hspace*{-.3ex}
}

\newcommand\Prufer{\hspace*{-.3ex}
\begin{tikzpicture}
\draw[thick, scale = 1.5] (.25ex,-.25ex) -- (-.35ex,.35ex) -- (.75ex,.85ex) -- (1.35ex, .25ex) -- cycle;
\draw[thick, scale = 1.5] (.5ex,.5ex) -- (1.35ex,.87ex) --  (1.35ex,1.47ex) -- (.5ex, 1.1ex) -- cycle;
\draw[thick, scale = 1.5] (.5ex,.3ex) circle[radius=.05ex];
\end{tikzpicture}
\hspace*{-.3ex}
}

\newcommand\Moore{\hspace*{-.3ex}
\begin{tikzpicture}
\draw[thick, scale = 1.5] (.25ex,-.25ex) -- (-.35ex,.35ex) -- (.75ex,.85ex) -- (1.35ex, .25ex) -- cycle;
\draw[thick, scale = 1.5] (.5ex,.55ex) -- (.5ex, 1.2ex) ;
\draw[thick, scale = 1.5] (.5ex,.55ex) circle[radius=.05ex];
\draw[thick, scale = 1.5] (.5ex,.27ex) circle[radius=.05ex];
\end{tikzpicture}
\hspace*{-.3ex}
}

\begin{document}

\newcommand{\N}{\mathbb{N}}
\newcommand{\R}{\mathbb{R}}
\newcommand{\Z}{\mathbb{Z}}
\newcommand{\Q}{\mathbb{Q}}
\newcommand{\C}{\mathbb{C}}
\newcommand{\LL}{\mathbb{L}}
\newcommand{\PP}{\mathbb{P}}

\newcommand{\esp}{\vskip .3cm \noindent}
\mathchardef\flat="115B

\newcommand{\lev}{\text{\rm Lev}}

\def\ut#1{$\underline{\text{#1}}$}
\def\CC#1{${\cal C}^{#1}$}
\def\h#1{\hat #1}
\def\t#1{\tilde #1}
\def\wt#1{\widetilde{#1}}
\def\wh#1{\widehat{#1}}
\def\wb#1{\overline{#1}}

\def\restrict#1{\bigr|_{#1}}

\newtheorem{lemma}{Lemma}[section]

\newtheorem{thm}[lemma]{Theorem}

\newtheorem{defi}[lemma]{Definition}
\newtheorem{nota}[lemma]{Notation}
\newtheorem{conj}[lemma]{Conjecture}
\newtheorem{cor}[lemma]{Corollary}
\newtheorem{prop}[lemma]{Proposition}
\newtheorem{prob}{Problem}
\newtheorem{probs}[prob]{Problems}
\newtheorem{ques}[lemma]{Question}
\newtheorem{rem}[lemma]{Remark}
\newtheorem{rems}[lemma]{Remarks}
\newtheorem{examples}[lemma]{Examples}
\newtheorem{example}[lemma]{Example}

\newtheorem{innercustomthm}{Problem}
\newenvironment{customprob}[1]
  {\renewcommand\theinnercustomthm{#1}\innercustomthm}
  {\endinnercustomthm}

\theoremstyle{remark}
\newtheorem{claim}{Claim}[lemma]
\renewcommand{\theclaim}{\arabic{section}.\arabic{lemma}.\arabic{claim}}

\title{On non-Hausdorff manifolds}
\author{Mathieu Baillif}

\maketitle
\begin{abstract}{\footnotesize 
          In this long note,
          we investigate various 
          purely topological aspects of non-Hausdorff manifolds (NH-manifolds for short).
          Our emphasis is on manifolds which exhibit homogeneity or weakenings thereof, 
          in particular being everywhere non-Hausdorff. Homogeneous NH-manifolds and 
          everywhere non-Hausdorff manifolds are respectively called HNH- and ENH-manifolds.
          We write $NH_X(x)$ for the 
          subset of points of a space $X$ which cannot be separated of $x$ by open sets.
          The topics covered in this note are the following.\\ 
          -- General (basic) properties of manifolds and their quasi-compact or quasi-countably compact subspaces.\\
          -- Covering properties implying the Hausdorffness of (weakly) homogeneous manifolds.\\
          -- (Non-)existence of hereditarily separable ENH-manifolds (under set theoretic hypotheses).\\
          -- Non-existence of a quasi-countably compact ENH-manifold.\\
          -- Properties of NH-manifolds which imply that $NH_M(x)$ is discrete, or at least ``simple''.\\  
          -- Constructions of HNH-manifolds such that $NH(x)$ is non-homogeneous, for instance
             a countable union of closed intervals and $n$-torii.\\
          -- Constructions of NH-manifolds $M$ with a point $x$ such that $NH_M(x)$ is
             homeomorphic to various ``complicated'' spaces, in particular in dimension $1$ and $2$.\\ 
          We use elementary (or at least well known) methods of general or set theoretic topology, 
          with a little bit of conformal theory and dynamical systems (flows) in some constructions.
          Many pictures are given to illustrate the constructions, and the proofs 
          are rather detailed, which is the main reason for the length of this note.
          }
\end{abstract}

{ }
\vskip .3cm
\noindent
{\bf Main changes from v4}: \\
Added Subsection \ref{subsec:hersep2} and Theorem \ref{thm:hersepall}, and amended the rest accordingly.
Also added Corollary \ref{cor:noENHqcc}. Moved (an updated version of) Theorem \ref{thm:Scountcpct}
from Section \ref{sec:cover} to Section \ref{sec:quasi-ctbly-cpct}.
\newpage
\tableofcontents
\newpage


\begin{epigraphs}
\qitem{In Chap. 9 we relax the hypothesis that our manifolds must be Hausdorff. Here I
would agree with Hirsch that it is hard to prove anything about them.}{David Gauld, Non-Metrisable Manifolds.}
\qitem{C'est la porte ouverte \`a tout.}{Strix Grizlyx}
\end{epigraphs}

\section{Introduction}
In this text,
a manifold is a connected space each of whose points has a neighborhood homeomorphic
to $\R^n$ for some fixed $n$. As indicated by the title, we do not assume Hausdorffness.
Our terminology is standard and we refer to Section \ref{sec:def} and \cite{Engelking}
for undefined terms in this introduction.
\\
These notes contain our investigations of purely topological properties of non-Hausdorff manifolds
(NH-manifolds, for short),
that is, we look at them from a general topology perspective.
It can be said that NH-manifolds are not a regular feature of mainstream mathematics, but they do pop up in some
areas. 
To our knowledge, the first paper taking them seriously
is due to A. Haefliger and G. Reeb in 1957 \cite{HaefligerReeb}, in relations to foliations of the plane. 
They have also been investigated since the 1960s by mathematical
physicists for possible models of space-time, and there are recent physics-minded contributions,
in particular by D. O'Connell, see \cite{OConnell:2023deRham, OConnell:2023, OConnell:2023vector}. 
Their homotopic aspects are also sometimes studied
(see e.g. the Bachelor thesis by A.G.M. Hommelberg and F.E.C. Ruijter  
\cite{HommelbergThesis:2014, RuijterThesis:2017}).
These contributions either use NH-manifolds mainly as a tool to study other objects, or 
investigate additional structures on NH-manifolds 
(differential structures, tangent bundles, homotopy groups, etc).
Apart from 
parts of A. Mardani's thesis \cite{MardaniThesis} and Chapter 9 of D. Gauld book \cite{GauldBook},
there does not seem to be many papers that study topological properties of NH-manifolds {\em per se}.
One of the reason is probably related to D. Gauld's quote above: it is difficult to prove anything 
interesting
about (potentially non-Hausdorff) manifolds in all generality.
Hence, our focus will be on manifolds which exhibit 
some kind of (topological) homogeneity, which tend to behave more reasonably.
For instance, it is known that Lindel\"of or even metaLindel\"of homogeneous manifolds
are actually Hausdorff
\cite{BaillifGabard, GartsideGauldGreenwood}. (It is a classical result that 
Hausdorff manifolds are homogeneous -- see e.g. \cite[Cor. 3.7]{GauldBook}.)
Actually, most of this note is a byproduct of our attack of the following problem:
\begin{prob}[{\cite[Problem 4.4]{BaillifGabard}}]
   \label{prob:intro1}
   Is there a homogeneous non-Hausdorff manifold without a (closed) uncountable discrete subset~?
\end{prob}
There is a great deal of parallelism between 
Problem \ref{prob:intro1} and a classical problem 
that asks whether a Hausdorff hereditarily separable manifold must be metrizable.
(It is easy to see that a manifold has no uncountable discrete subspace iff it is hereditarily separable, see
Corollary \ref{cor:hersepdis} below.)
This classical problem was solved in the 1970s:
there are hereditarily separable (perfectly normal, even) Hausdorff non-metrizable manifolds under 
the continuum hypothesis {\bf CH} (Rudin-Zenor 1976 \cite{RudinZenor}) 
and none under {\bf MA$(\omega_1)$} (Martin's axiom restricted to posets
of size $\le\omega_1$, Szentmikl\'ossy 1977 \cite{Szentmiklossy:1980}). \\
An important weakening of homogeneity
is the property of being {\em everywhere non-Hausdorff}, that is,
spaces for which $NH_M(x)$ is non-empty for all $x\in M$. 
In the now fifth version of these notes,
we are finally able to transfer these 70s techniques to this everywhere non-Hausdorff setting,
and our results mirror perfectly the classical aforementioned ones.
But before (partially) succeeding, we often got stuck
when trying to prove that the spaces we were investigating had nice enough properties to be able to carry
some argument or some construction further.
While personal limitations are certainly to blame, we also
believe that our (initial) lack of success is in part due to an inability to form an adequate mental image of
non-Hausdorff manifolds, homogeneous or not.\footnote{A rather unfortunate proof
of this fact is the embarrassing number of errors we found in previous versions of these notes.
We hope that we we able to get rid of the most serious ones.}
We then started to study NH-manifolds more systematically, sometimes drifting quite far away from 
Problem \ref{prob:intro1}, and ended up gathering quite a lot of material
that became the contents of these notes. 
In particular, we pay a close attention to the space $NH_M(x)$ of points of the space $M$
that cannot be separated from $x$ by open sets.
 
Let us now give a more precise description of the topics covered.
\begin{itemize}
\item[(i)] General (basic) properties of manifolds and their quasi-compact or quasi-countably compact subspaces.
           These are results (not all new or due to us) which are used pervasively throughout the paper
           and make it reasonably self-contained.
\item[(ii)] Covering properties implying the Hausdorffness of homogeneous manifolds, or more 
            generally which impede a manifold to be everywhere non-Hausdorff. In particular, 
            we generalize the above mentioned results about the Hausdorffness of homogeneous metaLindel\"of
            manifolds.
\item[(iii)] Properties of NH-manifolds which imply that $NH_M(x)$ is discrete, or at least ``simple'',
           for all $x\in M$. The simplest examples of homogeneous manifolds have discrete $NH_M(x)$, 
           here we look for properties that, together with weakenings of homogeneity,
           ensure the discreteness of $NH_M(x)$.
\item[(iv)] Existence of hereditarily separable everywhere non-Hausdorff manifolds under set theoretic hypotheses.
             That is: partial results about Problem \ref{prob:intro1}. We show in particular
             that they
             exist under {\bf CH}, and that they do not under {\bf MA$(\omega_1)$}. 
\item[(v)]  A closer look at quasi-countably compact subspaces of manifolds. In particular, we show that
            a quasi-countably compact NH-manifold cannot be everywhere non-Hausdorff
            (and thus neither homogeneous).
            This shows that Problem \ref{prob:intro1} has a positive answer
            in {\bf ZFC} if one replaces ``uncountable''
             by ``infinite''.
\item[(vi)] Constructions of homogeneous NH-manifolds such that $NH(x)$ is non-homogeneous, hence in particular
            non-discrete. For instance, we obtain homogeneous manifolds such that $NH_M(x)$ is
            a countable union of closed intervals and $n$-torii.
\item[(vii)] Construction of a homogeneous $1$-manifold $M$ such that for each $x,y\in M$
             there is $z\in M$ such that $NH_M(x)\ni z\in NH_M(y)$. 
\item[(viii)] Closer inspection of ``how complicated'' $NH_M(x)$ can be in not-at-all-homogeneous NH-manifolds.
             We give many examples and show many pictures in dimensions $1$ and $2$.\footnote{We
             apologize in advance to color blind readers: we tried 
             to draw our picture only with shades of grey, but did not obtain good results
             and ended up using colors in most of them.}
             For instance, we show that there are $1$-dimensional NH-manifold manifolds $M$ and some 
             $x\in M$ such that $NH_M(x)$
             is homeomorphic to classical spaces such as the Cantor space, ordinals $<\omega_1$, 
             special Aronszajn trees and Isbell-Mr\'owka $\Psi$-spaces;
             and $2$-dimensional ones such that $NH_M(x)$ is a copy of $\R$, of the long ray or 
             parts of boundaries of open domains in the plane. 

\end{itemize}

For Hausdorff manifolds, complicated behaviour starts in dimension $2$
(one could argue that it only starts in dimension $3$ for metrizable manifolds),
but NH-manifolds exhibit rather strange features already in dimension $1$. 
For this reason, almost all our examples are
of dimensions $1$ and $2$, and are build ``by hand'', using elementary
(or at least well known) techniques,
mainly from general (or set theoretic) topology.
Despite their simplicity, they were somewhat surprising to us: 
we underestimated the degree of weirdness a NH-manifold may exhibit, which might be due
to the fact that we are much more used to studying Hausdorff non-metrizable manifolds.
On the other hand, we were pleasantly surprised
to be able to use results and/or methods (some being new to us) 
that are also elementary, but belong to arguably distinct branches of topology,
for building some of these examples. To name some of these methods: 
set-theoretic inductive constructions of the 1970s in (iii);
Riemann--Caratheodory's classical conformal theory in (viii) and 
the theory of transitive/minimal flows on Hausdorff manifolds in (vi) and (viii).
In short, non-Hausdorff manifolds seem to be a great playground for dabbling in neighboring subjects,
learning something in passing,
and we had fun investigating them.
\vskip .3cm
The structure of this note follows the list above,
except in later sections. 
Here is a kind of roadmap. Section \ref{sec:def} is the only one to which the others constantly refer.
\begin{center}
\begin{tabular}{ccc}
    Topic & Section \\
    \hline
    (i)   & \ref{sec:def} \\
    (ii)  & \ref{sec:cover}\\
    (iii) & \ref{sec:NHdisc}\\
    (iv)  & \ref{sec:conj} \\
    (v)   & \ref{sec:quasi-ctbly-cpct} \\
    (vi)  & \ref{sec:flows} \\
    (vii) & \ref{sec:1dim} \\
    (viii)& \ref{sec:1dim}, \ref{sec:2dim}, \ref{sec:flows}  
\end{tabular}
\end{center}
A very short Section \ref{sec:ques} contains some further questions.
The level of difficulty of our examples is not always increasing:
some more involved spaces are described before simpler ones, when 
they fit in a particular topic. For instance, 
the first examples given in Section \ref{sec:1dim} are arguably simpler than those
in Sections \ref{sec:conj} and \ref{sec:quasi-ctbly-cpct}.
A reader with a background on general or set-theoretic topology should have no trouble 
reading this note.
A reader with no background on these subjects might
have trouble with some parts. For instance, 
ladder systems on $\omega_1$ 
and inductive constructions of lengths $\omega_1$ are used in 
Subsections \ref{subsec:hersep}, \ref{subsec:hersep2},
\ref{subsec:contrarian} and \ref{subsec:noomega_1}.
We hope that our arguments are clear enough to be followed.
\vskip .3cm
\noindent
{\bf Acknowledgements}\\
We thank 
Alexandre Gabard for having 
prompted our interest in the subject by asking us in 2004 which other
property one should add to homogeneity in order to ensure that a manifold is Hausdorff.


\section{Definitions and general properties of manifolds}\label{sec:def}

We follow the set theoretic tradition of denoting the integers by $\omega$,
and $\omega_1$ is the first uncountable ordinal. Ordinals are the sets of their predecessors.
Throughout this text, ``space'' means ``topological space''.
We use brackets $\langle\,\cdot\,,\,\cdot\,\rangle$ for pairs, reserving parenthesis 
$(\,\cdot\,,\,\cdot\,)$ for open intervals.
A space $X$ is {\em homogeneous} 
iff given $x,y\in X$, there is a homeomorphism $h$ of $X$ such that
$h(x) = y$.
Two points in a space may be {\em separated by open sets}, or in short just {\em separated}, iff
there are open disjoint $U,V$ containing one point each. 
A space is {\em Hausdorff} iff any two points in it can be separated.
As usual, the closure of a subspace $A$ in a space clear from the context is denoted by $\wb{A}$.

\subsection{Main definitions}

There is some ambiguity in the phrase ``non-Hausdorff manifold'': one might think of a manifold 
which is not assumed, but could be, Hausdorff, or one which is definitely not, in the sense that 
it contains a pair of points which cannot be separated.
Let us thus fix the terminology.
\begin{defi}\ \\
  $\bullet$
  A manifold, or $n$-manifold if we want to emphasize the dimension,
  is a connected space such that any point has a neighborhood homeomorphic to $\R^n$ for some (fixed)
  integer $n$. \\
  $\bullet$
  An Euclidean subset of an $n$-manifold is a subset homeomorphic to $\R^n$.\\
  $\bullet$
  A non-Hausdorff $n$-manifold, in short NH-$n$-manifold (or NH-manifold), is 
  an $n$-manifold in which at least two points cannot be separated.\\
  $\bullet$
  A HNH-$n$-manifold (or HNH-manifold) is a homogeneous NH-$n$-manifold.
\end{defi}

Notice in particular that Euclidean subsets are connected and open (by invariance of domain). 
The connectedness of a manifold implies that its dimension is fixed.
Any manifold, Hausdorff or not, is $T_1$ (see e.g. the easy proof in \cite[Prop. 9.2]{GauldBook})
and inherits the local properties of the Euclidean space
which do not depend on Hausdorffness: first countability, local second countability,
local separability and local Hausdorffness, for instance. (In each case, ``local'' means that
each point has a neighborhood -- or equivalently a base of neighborhoods -- with the property.) 
Any open subspace of a manifold is Baire
(that is: any countable intersection of open dense subsets is dense). 
The following lemma, which is an immediate consequence of first countability,
will be implicitely used many times.
\begin{lemma}
   If $M$ is a manifold, 
   then $M$ is Fr\'echet-Urysohn, that is: if
   $A\subset M$ and $y\in \wb{A}$, then there is a sequence of points of $A$
   converging to $y$.
\end{lemma}
Notice in passing that limits may be not unique.
The two simplest NH-manifolds are arguably those in the next example.
\begin{example}
   \label{ex:simplestones}
   \ \\
   (a)
   The line with two origins is the union of 
   $\R$ (with the usual topology) and another point $0^*$ whose neighborhoods are 
   the nicked intervals $(-\epsilon,0)\cup\{0^*\}\cup(0,\epsilon)$ ($\epsilon>0$).
   Then $0^*$ cannot be separated from $0$. \\
   (b) The branching line is the union of two copies of
       $\R$ where the strictly negative parts have been identified pointwise. 
       The two copies of $0$ cannot be separated.
\end{example}
Another way to define the line with two origins is by quotienting $\R\times\{0,1\}$ by 
$\langle x,0\rangle \sim \langle x,1\rangle$ whenever $x\not=0$, and we may define
the branching line as $\R\times\{0\}\cup [0,+\infty)\times\{1\}$ with neighborhoods of 
$\langle 0,1\rangle$ given by $(-\epsilon,0)\times\{0\}\cup[0,\epsilon)\times\{1\}$.
When defining NH-manifolds, among the two approaches (quotient space or direct description),
we will choose 
the one that we find more convenient for the example at hand (sometimes combining both).
We will use the following result implicitely 
in many constructions to ensure that when we use the quotient space approach, the resulting space is
indeed a manifold (for a proof, see \cite[Theorem 2.1]{OConnell:2023}).
\begin{thm}
  \label{thm:oconnelladjunctive}
  Let $\{M_i\,:\,i\in I\}$ ($I$ an index set) be Hausdorff manifolds, and 
  $\{\psi_{i,j}\,:\,\langle i,j\rangle\in J\}$ ($J$ an index set $\subset I^2$ which 
  contains the diagonal)
  be homeomorphisms from an open $U_{i,j}\subset M_i$ to an open $V_{i,j}\subset M_j$, such
  that $U_{i,i} = V_{i,i}=M_i$, $\psi_{i,i}=id_{M_i}$
  and $\psi_{i,k}(a) = \psi_{j,k}\circ\psi_{i,j}(a)$ for all $a\in U_{i,j}\cap U_{i,k}$.
  If the space $\sqcup_{i\in I}M_i$ quotiented by $M_i\ni x\sim y\in M_j \Leftrightarrow y=\psi_{i,j}(x)$
  is connected, it
  is a manifold.
\end{thm}
(The collection of $M_i$, $U_{i,j}$ and $\psi_{i,j}$ is called an adjunction system.
The condition on the $\psi_{i,j}$ ensures that our data yields a well defined quotient space.) Notice
that the dimension of each $M_i$ must be the same if the quotient space is connected (and hence a manifold).
\vskip .3cm
\noindent
Recall that every Hausdorff manifold is homogeneous, see e.g. \cite[Cor. 3.7]{GauldBook}.
This property may (and in general does) not hold in NH-manifolds.
The next definition generalizes the idea of the line
with two origins and gives probably the simplest examples of HNH-manifolds. The space
$\mathbf{G}(\R,2)$ is called the everywhere doubled line in \cite[Section 3]{BaillifGabard}.  
\begin{defi} If $M$ is a manifold and $\kappa$ a cardinal (finite or infinite),
             the NH $\kappa$-tuple of $M$, denoted $\mathbf{G}(M,\kappa)$,
             is the NH-manifold $M\times\kappa$ 
             where open neighborhoods of $\langle x,\alpha\rangle$ are of the 
             form $(U-F)\times\{0\}\cup (U\cap F)\times\{\alpha\}$
             for $U$ an open neighborhood of $x$ in $M$ and $F\ni x$ a finite subset.
\end{defi}             

\begin{lemma} \label{lemma:THNH}
              If $M$ is a homogeneous manifold and $\kappa$ a cardinal then 
              $\mathbf{G}(M,\kappa)$ 
              is a HNH-manifold.
\end{lemma}
\begin{proof}
   Interchanging $\langle x,\alpha\rangle$ with $\langle x,\beta\rangle$ (and leaving the other points fixed)
   is a homeomorphism. If $x,y\in M$, choose a homeomorphism $h$ such that $h(x) = y$, 
   then sending $\langle z,\alpha\rangle$ to $\langle h(z),\alpha\rangle$ 
   for each $\alpha$ defines an homeomorphism as well.
\end{proof}
Every countable subset of a Hausdorff manifold is contained in an Euclidean subset (see \cite[Corollary 3.4]{GauldBook}).
This is not true for NH-manifolds as shown by the next example.
\begin{example}
   \label{ex:HausNotReg}
   A countable Hausdorff subspace of a NH-manifold may be non regular. 
\end{example}
\begin{proof}[Details]
   Let $D_0,D_1$ be disjoint dense subsets of the rationals $\Q$, with $0\in D_0$.
   Consider the Hausdorff subset $X = D_0\times\{0\}\cup D_1\times\{1\}$
   of $\mathbf{G}(\R,2)$.
   Then $D_0\times\{0\}$ is open in $X$, contains $\langle 0,0\rangle$
   but no closed neighborhood of it. Indeed, any neighborhood of $\langle 0,0\rangle$ 
   in $X$
   has infinitely many points upstairs in its closure.  
\end{proof}

Let us borrow some notation from \cite{GartsideGauldGreenwood, MardaniThesis} (and add some tweaks to it).
\begin{defi}
   Let $X$ be a space and $x,y\in X$ with $x\not= y$.
   Write $x\bumpeq y$ iff $x$ cannot be separated from $y$ (by open sets), 
   and let $\Bumpeq$ be the transitive closure of $\bumpeq$.
   We set $x\not\bumpeq x$.
   If $A\subset Y \subset X$, let $NH_Y(A) = NH_Y^1(A) = \{y\in Y\,:\,\exists x\in A\text{ with } y\bumpeq x\}$,
   $NH^n_Y(A) = NH_Y(NH^{n-1}_Y(A))$, and 
   $H_Y^n(A) = X-NH^n_Y(A)$. 
   Finally, set $NH^\infty_Y(A) = \cup_{n\in\omega} NH_Y^n(A)$.
   $NH_Y^n(\{x\})$ is abbreviated as $NH_Y^n(x)$, and $NH^n_X(A),H^n_X(A)$ as $NH^n(A),H^n(A)$.
\end{defi}
\begin{rem} \ \\
   $NH(A)\cap A\not=\varnothing$ if $A$ is not Hausdorff.
\end{rem}
By definition, $NH^\infty_Y(A)=\{y\in Y\,:\,\exists x\in A\, x\Bumpeq y\}$.

\begin{lemma}
   \label{lemma:sep}
   Let $X,Y$ be spaces and $x,y\in X$ with $x\bumpeq y$. \\
   (a) For each open $U\ni x$ we have $y\in\wb{U}$. Hence, $NH(U)\subset\wb{U}$. \\
   (b) If $f:X\to Y$ is continuous, then either $f(x) = f(y)$ or $f(x)\bumpeq f(y)$ (in $Y$).\\
   (c) If $h:X\to Y$ is a homeomorphism, then $h(NH_X(x)) = NH_Y(h(x))$.\\
   (d) If $A\subset X$ satisfies $NH(A)\cap A=\varnothing$ and $NH(A)$ is closed, then $NH(A)$ is nowhere dense. \\
   (e) If $U\subset X$ is open Hausdorff, then $NH(U)\cap U=\varnothing$.\\
   (f) If $A$ is a subset of an open Hausdorff $U\subset X$, then $NH(A)$ is nowhere dense.
\end{lemma}
Notice that $NH(A)\cap A=\varnothing$ is a stronger assumption than $A$
being Hausdorff in the induced topology: for instance, in a $T_1$ space,
if $A=\{x,y\}$ with $x\bumpeq y$, then $A$ is discrete in the induced topology but $NH(A)\cap A = A$.
It happens that (d) does not hold if $A$ is only Hausdorff in the induced topology, even if $X$ is a manifold
(see Remark \ref{rem:FatS1} (3)). 

\begin{proof}
   (a), (b), (c) and (e) are immediate. For (f), by (a) and (e) 
   we have that $NH(A)\subset\wb{U}-U$, which is nowhere dense.
   Let us prove (d).
   Let $U$ be the the interior of $NH(A)$. If $U$ is non-empty,
   let $x\in U$ and $y\in A$ with $x\bumpeq y$.
   By (a) and closedness of $NH(A)$, we have
   $NH(U)\subset \wb{U}\subset NH(A)$. 
   But $y\in NH(U)$, thus $y\in NH(A)\cap A$ which is absurd since $NH(A)\cap A=\varnothing$.\\
\end{proof}

Quick consequences of the previous lemma are the following.

\begin{lemma}
   \label{lemma:fxfy}
   Let $X$ be a space, $x,y\in X$ with $x\Bumpeq y$, and $H$ be a Hausdorff space.
   Then $f(x) = f(y)$ for any continuous $f:X\to H$. 
\end{lemma}
\begin{proof}
   By Lemma \ref{lemma:sep} (b), $f(x) = f(y)$ if $x\bumpeq y$. The result follows by induction.
\end{proof}

\begin{lemma}
   \label{lemma:denseconstant}
   Let $X$ be a space, $D\subset X$ be such that $\forall x,y\in D$, $x\Bumpeq y$.
   Then for any Hausdorff space $H$ and continuous $f:X\to H$, $f$ is constant on $\wb{D}$.
\end{lemma}
\begin{proof}
   Immediate by the previous lemma.
\end{proof}

\begin{defi}
  A space $X$ is everywhere non-Hausdorff if $NH(x)\not=\varnothing$ for each $x\in X$.
  An ENH-manifold (or ENH-$n$-manifold) is an everywhere non-Hausdorff $n$-manifold.
\end{defi}
Notice that $A\subset NH_Y^2(A)$ for each $A\subset Y\subset X$ if $Y$ is everywhere non-Hausdorff.
It is immediate that an HNH-manifold is an ENH-manifold. A simple example of an ENH-manifold which is not homogeneous
is provided by $\mathbf{G}(\R,3)$ minus one point.

\begin{defi}
   If $M$ is a space and $A\subset M$, let $\mathcal{E}(A,M)$ [resp. $\mathcal{OH}(A,M)$]
   be the smallest cardinal 
   of a cover of $A$ by open Euclidean subsets of $M$ [resp. by open Hausdorff subsets of $M$] 
   Denote $\mathcal{E}(M,M)$ simply by $\mathcal{E}(M)$
   and the same for $\mathcal{OH}(M,M)$. 
\end{defi}

Of course, $\mathcal{E}(A,M)\ge\mathcal{OH}(A,M)$. 
The following lemma is immediate, since Euclidean subsets are Lindel\"of.
\begin{lemma}
   If
   the Lindel\"of number of a manifold $M$ is infinite, it is equal to $\mathcal{E}(M)$.
\end{lemma}
There are simple examples of manifolds with 
e.g. $\mathcal{E}(M)=\aleph_1>\mathcal{OH}(M)$, 
actually any Hausdorff non-metrizable manifold 
has the property that $\mathcal{E}(M)\ge\aleph_1>1 =\mathcal{OH}(M)$,
since Lindel\"ofness is equivalent to metrizability for Hausdorff manifolds,
see e.g. \cite[Chapter 2]{GauldBook}. Notice also the following 
consequence of this latter fact.
\begin{rem}
   If $M$ is a manifold and
   $\mathcal{E}(M)\le\aleph_0$,
   then any connected open Hausdorff subset of $M$ is a metrizable submanifold of $M$.
\end{rem}

\subsection{General properties of manifolds (and related spaces)}

Recall that a space is {\em locally Hausdorff} iff each point has a Hausdorff neighborhood
(or equivalently a base of Hausdorff neighborhoods).
A subset $U$ of a space $X$ is {\em H-maximal} 
iff it is open Hausdorff and there is no open Hausdorff $V$ properly containing $U$, and
$U$ is {\em CH-maximal} 
iff it is connected open Hausdorff and there is no connected open Hausdorff 
$V$ properly containing $U$.
The following lemma, quite useful in our context, 
was already noted in 1971 (for manifolds) by Hajicek. We repeat the proof for completeness.

\begin{lemma}[{\cite[Lemma 4.2]{BaillifGabard}, \cite[Theorem 1]{Hajicek:1971}}]
   Let $X$ be a space.
   \\
   (a) For each open Hausdorff $U$, 
   there is a H-maximal $V\supset U$.\\
   (b) For each connected open Hausdorff $U$, 
   there is a CH-maximal $V\supset U$.\\
   (c) If $X$ is a locally Hausdorff space, then
   each H-maximal open Hausdorff subset is dense,
   and each point of $X$ is contained in a H-maximal subset.\\
   (d) If each point of $X$ has an open connected Hausdorff neighborhood, then
   each point of $X$ is contained in a CH-maximal subset.
   \label{lemma:locHaus}
\end{lemma}
\begin{proof}\ \\
   (a) Let $\mathfrak{U}=\{V\supset U\,:\,V\text{ open Hausdorff}\}$ ordered by the inclusion.
   Any totally ordered subset of $\mathfrak{U}$ has an upper bound in $\mathfrak{U}$ given by its union,
   hence by Zorn's lemma
   there is a maximal element.\\
   (b) Any increasing union of connected open Hausdorff subsets is connected open Hausdorff, hence the proof is 
       the same as in (a).\\
   (c) If an open Hausdorff $V$ is not dense, take an open $W$ with $W\cap V=\varnothing$.
   Choose some $x\in W$ and an open Hausdorff $A\ni x$, then $(A\cap W)\cup V$
   is a Hausdorff open set properly containing $V$, which is thus not maximal. 
   The other claim follows immediately by (a).\\
   (d) Immediate by (b).
\end{proof}

A CH-maximal subset of a manifold is not H-maximal in general: in the branching line,
the bottom copy of $\R$ is CH-maximal but not H-maximal (and not dense).

\begin{lemma}
   \label{lemma:Yclosed} 
   Let $X$ be a space and $x\in X$.
   Then, the subspace $NH(x)\cup\{x\}$ is closed.
   If $X$ moreover locally Hausdorff, then $NH(x)$ is closed nowhere dense.
\end{lemma}
\begin{proof}
   If $z\not\in NH(x)\cup\{x\}$, there are open disjoint $U\ni z$, $V\ni x$. Hence $U\cap NH(x)=\varnothing$. 
   If $X$ is locally Hausdorff, then any open Hausdorff $U\ni x$
   does not intersect $NH(x)$ which is thus closed. 
   The rest follows either by Lemma \ref{lemma:sep} (d) or Lemma \ref{lemma:locHaus} (c).
\end{proof}

The next two lemmas concern limits of sequences. 
\begin{lemma}
   \label{lemma:twolimits}
   Let $X$ be first countable and 
   $x,y\in X$ with $x\bumpeq y$.
   Then there is a sequence $z_i\in X$ ($i\in\omega$)
   which converges to both $x$ and $y$. 
\end{lemma}
\begin{proof}
   Choose countable bases of neighborhoods $U_i,V_i$ of $x,y$ and take 
   $z_i\in U_i\cap V_i$.
\end{proof}

\begin{lemma}
   \label{lemma:NHseq}
   Let $X$ be a space. Suppose that $x_i,y_i$ ($i\in\omega$) are sequences of points of $X$,
   that
   $y\in X$ is a limit of the $y_i$, and that $x$ is a cluster point of the $x_i$.
   Then the following hold.
   \\
   (a) Assume that for each $i$, either $x_i\bumpeq y_i$ or $x_i=y_i$.
   Then either $x\bumpeq y$ or $x=y$. \\
   (b) Assume that $X$ is moreover locally Hausdorff
   and $x_i\bumpeq y_i$ for each $i\in\omega$. Then $y\not= x\bumpeq y$.
\end{lemma}
\begin{proof}\ \\
   (a)
   Suppose that $x\not= y$.
   Let $U\ni y$ and $V\ni x$ be open.
   Then $U$ contains all but finitely many $y_i$, and $V$ contains infinitely many $x_i$.
   Hence for some $i$, $x_i\in V$ and $y_i\in U$. Since either $x_i\bumpeq y_i$ or $x_i=y_i$,
   $U\cap V\not=\varnothing$. Hence, $x\bumpeq y$. \\
   (b)
   If $X$ is locally Hausdorff and $x=y$, take an open Hausdorff neighborhood $W\ni x$. 
   As in (a), both $x_i$ and $y_i$ are in $W$ for some $i$. But this is impossible since $x_i\bumpeq y_i$.
\end{proof}

Hausdorffness is assumed in the usual definition of compactness. To avoid ambiguities, we
say that as space is {\em quasi-[countably-]compact} iff 
every [countable] open cover of it has a finite subcover, and {\em [countably]
compact} if it is moreover Hausdorff.
Recall that a $T_1$-space is quasi-countably compact iff every infinite subset has a cluster point
iff there is no infinite closed discrete subspace.
A space is {\em locally compact} iff any neighborhood of a point contains a compact neighborhood.
Hence, locally compact spaces are locally Hausdorff.
Euclidean subsets are locally compact, hence so are manifolds; {\em however}
the reader must keep in mind that
compact subspaces of a manifold may be non-closed, 
and the intersection of $2$ compact subsets may be non-compact.

\begin{lemma}
   \label{lemma:NHK}
   If $X$ is a Fr\'echet-Urysohn locally Hausdorff space, and $K\subset X$ is quasi-countably compact,
   then $NH(K)$ is closed. If moreover $NH(K)\cap K=\varnothing$,
   then $NH(K)$ is  nowhere dense in $X$.
\end{lemma}
\begin{proof}   
   Let $K$ be quasi-countably compact.
   Let $y_i\in NH(K)$ ($i\in\omega$) be a sequence converging to $y\in\wb{NH(K)}$.
   Let $x_i\in K$ such that $y_i\bumpeq x_i$.
   Since $K$ is quasi-countably compact, there is a cluster point $x\in K$
   of the $x_i$ (locally Hausdorff spaces are $T_1$).   
   By Lemma \ref{lemma:NHseq}, 
   $y\bumpeq x$, and hence $y\in NH(K)$, which is thus closed.
   Now, if $NH(K)\cap K=\varnothing$, then $NH(K)$ is nowhere dense by Lemma \ref{lemma:sep} (d).
\end{proof}

\begin{cor}
   \label{cor:quasi-compact}
   If $X$ is a Fr\'echet-Urysohn locally Hausdorff space, and 
   $K\subset M$ is quasi-compact. Then $NH(K)$ is closed nowhere dense.
\end{cor}
\begin{proof}
   $NH(K)$ is closed by the previous lemma.
   By quasi-compactness and local Hausdorffness we may cover $K$ by finitely many open Hausdorff $V_0,\dots,V_n$.
   By Lemma \ref{lemma:sep} (f), $NH(K\cap V_i)$ is nowhere dense, and thus so is $NH(K) = \cup_{i=0,\dots,n}NH(K\cap V_i)$.
\end{proof}

\begin{cor}
   \label{cor:no-quasi-cpct}
   Let $M$ be a NH-manifold with a subset $A$ which is dense in an Euclidean subset $U$,
   and such that $NH(x)\not=\varnothing$ when $x\in A$.  
   Then $M$ is not quasi-compact. In particular, there are no quasi-compact ENH-manifolds.
\end{cor}
\begin{proof}
   $NH(M)$ contains $A$, which is impossible if $M$ is quasi-compact by the previous corollary.
\end{proof}

Lemma \ref{lemma:NHK} and its corollaries \ref{cor:quasi-compact}--\ref{cor:no-quasi-cpct}
also hold for quasi-countably compact subspaces, the proof is given in 
Section \ref{sec:quasi-ctbly-cpct}, see Theorem \ref{thm:noENHqcc} and Corollary
\ref{cor:noENHqcc}.
Even if $K$ is compact,
$NH(K)$ is not quasi-countably compact in general. For instance, 
in $\mathbf{G}(\R,2)$, $NH([0,1]\times\{0\})$ is $[0,1]\times\{1\}$ (with the discrete topology).
Also, the conclusion does not hold if $K$ is only assumed to be closed: any
subset of $\R\times\{1\}$ is closed in $\mathbf{G}(\R,2)$, but 
$NH(\Q\times\{1\}) = \Q\times\{0\}$ or $NH(\R\times\{1\}) = \R\times\{0\}$ are
not, for instance. 
Compact subsets have the following property (compare with Example \ref{ex:HausNotReg})
which does not hold for countably compact subspaces (see Example \ref{ex:omega1notinchart} below).

\begin{lemma}
   \label{lemma:compactinchart}
   Let $M$ be a manifold and 
   $K\subset M$ be compact with $NH(K)\cap K=\varnothing$. Then there is an open Hausdorff $U\supset K$.
\end{lemma}

\begin{proof}
   For each $x\in K$ choose an Euclidean $U_x\ni x$. By Lemma \ref{lemma:NHK}, 
   $V_x = U_x-NH(K)$ is open and $U_x\cap K = V_x\cap K$ since $NH(K)\cap K=\varnothing$. 
   Notice that $V_x\cup K$ is Hausdorff for each $x$:
   if $v\in V_x$ and $w\in K$, they can be separated, since otherwise $v\in NH(K)$.
   Let $W_x\subset V_x$ be open such that $A_x = \wb{W_x}\cap V_x$
   is compact and $x\in W_x$. 
   Choose a finite subcover $\{W_{x_i}\,:\,i=0,\dots,n\}$ of $K$. 
   Notice that $NH(\cup_i A_{x_i})\cap K = \varnothing$ and that 
   $NH(\cup_i A_{x_i})\supset NH(\cup_i W_i)$.
   By Lemma \ref{lemma:NHK} again,
   for each $i$, $O_i = W_{x_i} - NH(\cup_{i} A_{x_i})$ is open,
   hence $\cup_i O_i$ is open Hausdorff and
   contains $K$.
\end{proof}

\begin{lemma}
   \label{lemma:cpctVU}
   Let $V\subset U\subset M$ be such that $U,V$ are open Hausdorff in the manifold $M$.
   If $\wb{V}\cap U$ is compact, then $\wb{V}-U \subset NH(\wb{V}\cap U)$.
\end{lemma}
\begin{proof}
   Let $x\in\wb{V}-U$. 
   Since $V\subset U$, any neighborhood of $x$ must intersect $V = V\cap U$ and thus also $\wb{V}\cap U$.
   If $x\not\in NH(\wb{V}\cap U)$, for each $y\in\wb{V}\cap U$ 
   there are open disjoint $A_y\ni x$ and $B_y\ni y$.
   By compactness there are disjoint open $A\ni x$ and $B\supset \wb{V}\cap U$,
   a contradiction.
\end{proof}
Notice that the fact that $\wb{V}\cap U$ is compact 
does not imply $\wb{V}-U \supset NH(\wb{V}\cap U)$.
Take for instance the branching line defined in Example \ref{ex:simplestones} (b), with
$U=\R\times\{0\}$, $V=(0,1)\times\{0\}$. Then $\wb{V} = [0,1]\times\{0\}$, hence
$\langle 0,1\rangle$ is in $NH(\wb{V}\cap U)$
while $\wb{V}-U=\varnothing$.
Also, taking $V=U=(0,+\infty)\times\{0\}$ in the same space yields
$\wb{V}-U=\{\langle 0,0\rangle\}\not\subset NH(\wb{V}\cap U) = \varnothing$; hence, some form of compactness of
$\wb{V}\cap U$ is needed in general.
Lemma \ref{lemma:cpctVU} has a nice variant due to 
Hajicek
which gives a characterisation of CH-maximal subsets. 
(Be aware that Hajicek does use the terms manifolds and submanifolds,
  but does not give their formal definition. He however  
  always assumes them to be connected -- otherwise strange things may happen, see Example \ref{ex:notHmax}.) 

\begin{thm}[{{\cite[Theorem 2]{Hajicek:1971}}}]
   \label{thm:Hajicek}
   Let $M$ be a manifold.
   \\
   (a) 
   A connected open Hausdorff $U\subset M$ is CH-maximal iff $\wb{NH(U)} = \wb{U}-U$. \\
   (b) If $U\subset M$ is H-maximal, then $\wb{NH(U)} = \wb{U}-U = M-U$.
\end{thm}
\begin{proof}
  We repeat Hajicek's proof.\\
  (a) 
  Let $U$ be open, connected and Hausdorff. Assume that $\wb{NH(U)} = \wb{U}-U$,
  we show that $U$ is CH-maximal. Suppose it is not the case, then there is an open connected Hausdorff
  $V$ properly containing $U$. It follows that there is some $x\in (\wb{U}-U)\cap V$, otherwise
  $V=(V-\wb{U})\cup U$ is not connected. But $x\not\in\wb{NH(U)}$, as $V\cap NH(U)=\varnothing$.
  It follows that $U$ is CH-maximal.\\
  We now prove (b) and the converse implication of (a) together.
  Assume thus that $U$ is either CH-maximal or H-maximal.
  By Lemma \ref{lemma:sep} (a), $NH(U)\subset\wb{U}-U$,
  so by closedness $\wb{NH(U)}\subset \wb{U}-U$.
  Take now $x\in\wb{U}-U$.
  Suppose that there is a Hausdorff neighborhood $V$ of $x$ such that $V\cap NH(U) =\varnothing$.
  We may assume $V$ to be connected.
  Then $V\cup U$ is Hausdorff: given two points $p,q\in U\cup V$, 
  if both are in $V$ or both are in $U$, they can be separated;
  but if $p\in U$, $q\in V$, then $q\not\in NH(p)\subset NH(U)$, hence they can also be separated.
  Notice that $V\cup U$ is connected if $U$ and $V$ are connected.
  But $V\cup U$ strictly contains $U$, a contradiction with its CH- or H-maximality.
  Hence any neighborhood of $x$ intersects $NH(U)$ and thus $x\in \wb{NH(U)}$.
  It follows that $\wb{U}-U\subset \wb{NH(U)}$.
  \\
  We finish by noting that, in (b), $\wb{U}=M$ by
  Lemma \ref{lemma:locHaus} (c).
\end{proof}

The converse of (b) does not hold.

\begin{example}\label{ex:notHmax}
  A manifold with an open Hausdorff $U$ satisfying $\wb{NH(U)} = \wb{U}-U$, but $U$ 
  is neither connected nor H-maximal.
\end{example}
\begin{proof}[Details]
   We take $M$ to be the doubly-branching line, that is,
   $M = \R\times\{0,1,2\}/\sim$, where $\langle x,i\rangle\sim \langle y,j\rangle$
   iff $x=y <0$ (i.e. we identify pointwise the strictly negative parts).
   Let $U = \R\times\{0\}\cup(0,+\infty)\times\{1\}$. Then $U$ is clearly Hausdorff but neither
   connected nor H-maximal, and satisfies $\wb{NH(U)} = \wb{U}-U = \{\langle 0,1\rangle,\langle 0,2\rangle\}$.
\end{proof}

A last remark on local compactness, which is a counterpoint to Example \ref{ex:HausNotReg}.
\begin{lemma}
   \label{lemma:openclosedloccpct}
   Let $M$ be a manifold, $B\subset M$ be closed and $U\subset M$ be open.
   Then $U\cap M$ is locally compact (in the induced topology).
\end{lemma}
\begin{proof}
   Indeed, let $V\subset U\cap B$ be open in the induced topology, and $x\in V$.
   Let $W$ be open in $M$ such that $W\cap U\cap B = V$. We may assume that $W\subset U$, hence $W\cap B = V$.
   Take $W_1\subset W$ containing $x$ such that $\wb{W_1}\cap W$ is compact.
   Then $\wb{W_1}\cap U\cap B = \wb{W_1}\cap B$ is compact and contains $x$.
\end{proof}

Let us end this section with a simple observation about 
products.

\begin{lemma}
  \label{lemma:NHproduct}
  Let $\kappa$ be a cardinal, and $X_\alpha$ be a space for each $\alpha\in\kappa$.
  Let $X = \Pi_{\alpha<\kappa}X_\alpha$ with either the box topology or the usual product topology.
  Then for each $x\in X$, 
  $NH_X(x) = \Pi_{\alpha<\kappa} (NH_{X_\alpha}(x_\alpha)\cup\{x_\alpha\}) - \{x\}$, where $x_\alpha$ is the projection
  of $x$ on $X_\alpha$.
\end{lemma}
\begin{proof}
   Let $Y(x) = 
   \Pi_{\alpha<\kappa} (NH_{X_\alpha}(x_\alpha)\cup\{x_\alpha\}) - \{x\}$. We show that
   $Y(x) = NH_X(x)$ 
   for the box topology, hence a basic open set in $X$ is the product of open sets of $X_\alpha$.
   Let $u\in Y(x)$, and let $U = \Pi_{\alpha\in\kappa}U_\alpha$,
   $V = \Pi_{\alpha\in\kappa}V_\alpha$ be basic open sets respectively containing $u,x$.
   Since either $u_\alpha\bumpeq x_\alpha$ or $u_\alpha=x_\alpha$, 
   $U_\alpha\cap V_\alpha\not=\varnothing$ for each $\alpha$,
   hence $U\cap V\not=\varnothing$ and $u\in NH_X(x)$.\\
   Assume now that $u\not\in Y(x)$.
   Let $\alpha$ be such that $u_\alpha\not\bumpeq x_\alpha$ and $u_\alpha\not= x_\alpha$.
   There are thus open disjoint $U_\alpha,V_\alpha\subset X_\alpha$ with $u_\alpha\in U_\alpha$ and $x_\alpha\in V_\alpha$.
   Set $U'_\alpha = U_\alpha$, $V'_\alpha = V_\alpha$ and $V'_\beta = U'_\beta = X_\beta$ when $\beta\not=\alpha$.
   Then $\Pi_{\beta\in\kappa} U'_\beta$ and $\Pi_{\beta\in\kappa} V'_\beta$ separate $u$ and $x$, hence $u\not\in NH_X(x)$.
\end{proof}


\section{Covering properties implying Hausdorffness}\label{sec:cover}

This short section is concerned with the following question.
\begin{ques}
   Which covering properties do imply the Hausdorfness of homogeneous manifolds,
   or more generally impede the existence of everywhere non-Hausdorff ones~?
\end{ques}

Recall that $\mathcal{E}(M)\le\omega$ iff $M$ is Lindel\"of.
The main result of \cite{BaillifGabard} is that if $M$ is a Lindel\"of homogeneous manifold, 
then $M$ is Hausdorff.
This result was improved by P. Gartside, D. Gauld and S. Greenwood \cite{GartsideGauldGreenwood}:
A metaLinde\"of homogeneous manifold is Hausdorff. 
Both assumptions also imply that the manifold is actually metrizable, 
since Hausdorff metaLindel\"of manifolds are such \cite[Chapter 2]{GauldBook}.
The proofs of these theorems are actually rather easy, and
we give slight generalizations below, which have the advantage of 
not implying metrizability.

\begin{defi}
   Let $X$ be a space and $Y\subset X$ be a subset such that $Y\cap NH(Y)=\varnothing$.
   A scan of $Y$ in $X$ is a subset $S\subset NH(Y)$ such that 
   for each $x\in Y$, 
   $$ NH(x)\not=\varnothing\quad\Longleftrightarrow\quad NH(x)\cap S\not=\varnothing.$$
\end{defi}
In particular, $NH(Y)$ itself is a scan of $Y$ in $X$, but taking one point in $NH(x)$
(when available) for each $x\in Y$ is also one.
Notice that $S\cap Y=\varnothing$ since we assume $Y\cap NH(Y)=\varnothing$.
Say that a subspace of a space $X$ is {\em $\kappa$-Baire in $X$} ($\kappa$ a cardinal)
iff it cannot be covered by $<\kappa$ nowhere dense sets of $X$.
\begin{thm}
   \label{thm:ClH}
   Let $X$ be a locally Hausdorff, everywhere non-Hausdorff space.
   Let $A\subset U\subset X$ be such that $U$ is open Hausdorff and $A$ is $\kappa$-Baire in $X$
   for some $\kappa\ge\omega_1$. Finally, let
   $S$ be a scan of $A$ in $X$. 
   Then $S$ cannot be covered by $<\kappa$-many open Hausdorff 
   sets. 
\end{thm}
In other words: $\mathcal{OH}(S,X)\ge\kappa$. Notice that
$\wb{U}-U\supset S$ by Lemma \ref{lemma:sep}.

\begin{proof} 
   Notice that $NH(A)\cap A=\varnothing$ since $A\subset U$.
   By way of contradiction,
   let $\mathcal{V}=\{V_\alpha\,:\,\alpha\in\lambda\}$, $\lambda<\kappa$, 
   be a cover of $S$ by open Hausforff sets. 
   For each $\alpha\in\lambda$ let $W_\alpha\supset V_\alpha$ be 
   open Hausdorff and dense in $X$ (given by Lemma \ref{lemma:locHaus}).
   Then $\mathcal{W} = \{W_\alpha\,:\,\alpha\in\lambda\}$ is a cover of $S$ by 
   $\lambda$-many open dense Hausdorff sets.
   Since $\lambda < \kappa$, there is some $x$ in the
   non-empty intersection $A\cap\left(\cap_{\alpha\in\lambda}W_\alpha\right)$.   
   Since $x$ can be separated from each member of $S$, $S\cap NH(x)=\varnothing$.
   This implies that $NH(x) = \varnothing$ as well, hence $X$ is not everywhere non-Hausdorff.
\end{proof}

Recall that a space is {\em $\omega_1$-Lindel\"of} iff any cover by at most $\omega_1$ open sets has
a countable subcover iff any subset of cardinality $\omega_1$ has a point of complete accumulation.
The proof of the next lemma is directly taken from \cite[Theorem 4.1]{ArhanBuzy:1998b}, 
we include it for completeness.
\begin{lemma}
   \label{lemma:omega_1L}
   An $\omega_1$-Lindel\"of subset of a manifold is hereditarily Lindel\"of.
\end{lemma}
\begin{proof}
   Let $N\subset M$ be $\omega_1$-Lindel\"of, with $M$ a manifold.
   If $N$ is not hereditarily Lindel\"of, it contains a subspace $Y=\{y_\alpha\,:\,\alpha<\omega_1\}$
   which is right separated, i.e. $\{y_\alpha\,:\,\alpha\le\beta\}$ is open in $Y$ for each $\beta<\omega_1$.
   (Indeed, there is a strictly increasing $\omega_1$-sequence
   of open sets $U_\alpha$. Take $y_\alpha\in U_{\alpha+1}-U_\alpha$.)
   Then there is a point of complete accumulation $y$ of $Y$ in $N$. 
   The intersection of any Euclidean
   neighborhood of $y$ with $N$
   contains uncountably many $y_\alpha$ and hence a non-Lindel\"of subspace, a contradiction.
\end{proof}

As usual, the smallest cardinal $\kappa$ such that $\R^n$ ($1\le n<\omega$)
may be covered by $\le\kappa$ nowhere dense subsets
is denoted by $\text{cov}(\mathcal{B})$. 
(It is well known that it does not depend on $n$.)
Of course, $\omega_1\le\text{cov}(\mathcal{B})$ by Baire's theorem.

\begin{cor}
   \label{cor:ClH}
   Let $M$ be an ENH-manifold. 
   Let $A\subset U\subset M$ with $U$ open Hausdorff and $A$ 
   $\kappa$-Baire in $M$ (e.g. $A=U$ and $\kappa = \text{cov}(\mathcal{B})$).
   Let $S$ be a scan of $A$.
   Then $S$ neither (i) is $\omega_1$-Lindel\"of, nor 
   (ii) can be covered by $<\kappa$-many Hausdorff 
   open sets.
\end{cor} 

\begin{proof}
   It is enough to prove (ii) by
   Lemma \ref{lemma:omega_1L}, and it follows immediately from Theorem \ref{thm:ClH}.
\end{proof}

A collection $\mathcal{U}$
of open subsets of some space $X$ is said to be {\em point-$\kappa$ on $Y\subset X$} iff
every member of $Y$ is contained in $\le\kappa$ members of $\mathcal{U}$.
If $\mathcal{U}$ is a collection of subsets of some space $X$ and $A\subset X$, we set
$\mathcal{U}_A = \{U\in\mathcal{U}\,:\,U\cap A\not=\varnothing\}$.
The following easy lemma is well known and almost immediate, see e.g. \cite[Lemma 3.4]{Gauld:1998}. 
\begin{lemma} 
   \label{lemma:point-countable}
   Let $\kappa\ge\omega$ be a cardinal,
   $Y$ be a separable subspace of a space $X$, and $E$ be countable and dense in $Y$.
   Let $\mathcal{U}$ be an open cover of $X$ which is point-$\kappa$ on $E$. 
   Then $|\mathcal{U}_Y|\le\kappa$.
\end{lemma}

The most general statement about these covering properties we could prove is the following.
\begin{thm}
   \label{thm:pointkappa}
   Let $\kappa\le\text{cov}(\mathcal{B})$ be a cardinal,
   $M$ be an ENH-manifold and $S$ be a scan of
   some subset $A$, $\kappa$-Baire in $M$, contained in a
   Euclidean subset $U\subset M$ (e.g. $A=U$, $\kappa=\text{cov}(\mathcal{B})$). 
   Finally, let $E$ be dense in $A$ and $\lambda<\kappa$. Then
   no cover of $S$ by open Hausdorff sets is point-$\lambda$
   on $E$. 
\end{thm}
\begin{proof}
   We may assume that $E$ is countable. 
   Since $E$ is dense in $A$, it is dense in 
   $A\cup S$.
   Let $\mathcal{V}$ be a cover of $A\cup S$ which is point-$\lambda$ on $E$.
   We may assume that $U\in\mathcal{V}$, otherwise we just add it.
   By Lemma \ref{lemma:point-countable}, $|\mathcal{V}|\le\lambda$, and $U\cup S$  
   is covered by $<\kappa$ many Hausdorff open sets.
   We conclude with Corollary \ref{cor:ClH} (ii).
\end{proof}

In Theorems \ref{thm:ClH}, \ref{thm:pointkappa} and Corollary \ref{cor:ClH}, we may 
replace the scan $S$ by $\wb{U}$ in the statements and have the same conclusions.
The last result of this section is another variation in the same spirit.

\begin{lemma}
   \label{lemma:opennotFsigma}
   Let $U$ be Euclidean in the manifold $M$. 
   If $U$ is the union of $<\text{cov}(\mathcal{B})$
   closed sets of $M$, then there is an open $V$, dense in $U$,
   such that $NH(V) = \varnothing$. In particular,
   $M$ is not everywhere non-Hausdorff.
\end{lemma}
\begin{proof}
   Let $\{F_\alpha\,:\,\alpha<\kappa\}$, $\kappa<\text{cov}(\mathcal{B})$, be closed
   in $M$ such that $\cup_{\alpha<\kappa}F_\alpha = U$.
   Let $V = \cup_{\alpha<\kappa}\text{int}(F_\alpha)$, where $\text{int}$ denotes the interior. 
   Then $V$ is dense in $U$, otherwise there is some 
   open $W\subset U$ such that $F_\alpha\cap W$ is nowhere dense for each $\alpha$, 
   which is a contradiction with $\kappa<\text{cov}(\mathcal{B})$.
   Let $x\in V$, then $x\in \text{int}(F_\alpha)$ for some $\alpha$. 
   By closedness and Lemma \ref{lemma:sep} (a), $U\supset F_\alpha\supset \wb{\text{int}(F_\alpha)}\supset NH(x)$. 
   Since $U\ni x$ is Hausdorff, $NH(x) = \varnothing$.
\end{proof}

We close this section with a related question.
\begin{ques}
   Is there an ENH- or even HNH-manifold with $\mathcal{E}(M)>\mathcal{OH}(M)$~?
\end{ques}


\section{Discreteness of $NH(x)$ and weakenings of homogeneity in NH-manifolds}\label{sec:NHdisc}

This section is dedicated to investigating the problem below.

\begin{prob}
   In a HNH-manifold $M$, which additional properties ensure that $NH(x)$ is discrete,
   or at least homogeneous, for every $x\in M$~?
   \label{prob:NHd}
\end{prob}

Recall that a subset of a manifold is discrete iff there is no converging sequence in it.
We ask for additional properties because there are examples 
of homogeneous $2$-manifolds with non-homogeneous $NH(x)$ 
(Example \ref{ex:sines} below and those after it)
although we were not able to settle the $1$-dimensional case. 
We came up with three possible such properties: the existence of a sorted neighborhood 
in any dimension
(Definition \ref{defi:sort}), and, in dimension $1$,
f-rim-simplicity (a notion inspired 
by Hajicek, see Definition \ref{defi:simple}) and the existence of homeomorphisms fixing $x$ and shuffling $NH(x)$.
We dedicate a subsection to each one. Before describing them,
notice that the naive approach to proving that $NH(x)$ must be discrete for HNH-$1$-manifolds is
to show that if $NH(x)$ contains a convergent sequence, then no homeomorphim can send $x$ to the limit 
point of the sequence. This does not work so directly, as seen in 
Example \ref{ex:NHnotdiscrete} below (which is not homogeneous, but has such homeomorphism).
Notice also the following easy lemma.

\begin{lemma}
   \label{lemma:NHopendisjoint}
   Let $\mathcal{D}=\{D_\alpha\,:\,\alpha\in\lambda\}$ be a disjoint family of open sets
   in a NH-manifold $M$, and let $x\in M$.
   Then $NH(x)$ intersects at most countably many $D_\alpha$.
\end{lemma}
\begin{proof}
   Let $U$ be an Euclidean neighborhood of $x$, if $D_\alpha$ intersects $NH(x)$, then 
   $D_\alpha\cap U\not=\varnothing$. The lemma follows by disjointness of $\mathcal{D}$
   and secound countability of $U$.
\end{proof}

\subsection{Sorted neighborhoods}\label{subsec:sort}

\begin{defi}
   \label{defi:sort}
   Let $U$ be open Hausdorff in a manifold $M$.
   We say that $U$ is sorted [resp. weakly sorted] 
   in $M$ iff there is a continuous $f:\wb{U}\to H$ called the [weakly] sorting map, 
   where $H$ is a Hausdorff space,
   such that $f\upharpoonright U$ is $1$-to-$1$ [resp. finite-to-one] on its image.
   We say that $M$ is  [weakly] sortable iff each $x\in M$ has a neighborhood 
   which is a [weakly] sorted. 
   An open Hausdorff [weakly] sorted subset is called a [weakly] sorted neighborhood. 
\end{defi}

As we shall see, a Hausdorff open subsets $U$ is sorted whenever closure is in a sense ``well organized''.
For instance, 
$N\times\{0\}$ (if $N$ is Hausdorff) 
is a sorted neighborhood in $\mathbf{G}(N,\kappa)$, the sorting map being the projection
on the first coordinate.
The existence of a (weakly) sorted neighborhood is a quite strong assumption.

\begin{lemma}
   \label{lemma:sortBumpeq}
   Let $U\subset X$ be open Hausdorff and $x\in U$. If $U$ is weakly sorted, then 
   $NH^\infty(x)\cap U$ is finite. 
   If $U$ is a sorted, then $NH^\infty(x)\cap U\subset \{x\}$.
\end{lemma}
\begin{proof}
   Immediate by Lemma \ref{lemma:fxfy}. Notice that $NH^\infty(x)$ is empty if $NH(x)$ is.
\end{proof}

The following proposition shows the impact of sorted neighborhoods on Problem \ref{prob:NHd}.

\begin{prop}
   \label{prop:homogensort}
   Let $M$ be a manifold and $U\subset M$ be a 
   weakly sorted neighborhood 
   containing $x$. 
   Suppose that for each $y\in NH(x)$ there is a continuous $h:M\to M$
   such that $h(y) = x$ and $h$ restricted to $\{x\}\cup NH(x)$ is $1$-to-$1$. 
   Then $NH(x)$ is discrete. 
\end{prop}
\begin{proof}
   Suppose that there is a converging sequence $y_n\to y$ in $NH(x)$. 
   Let $h$ be as in the statement of the theorem. 
   Then, $h(y_n)\bumpeq h(x) \bumpeq h(y) = x$
   and $h(x)\in\wb{U}$
   (because $NH(x)\subset\wb{U}$).
   Since $h(y_n)$ converges to $h(y) = x$, there are infinitely many $n$ such that $h(y_n)\in U$,
   hence any map $f:\wb{U}\to H$ with Hausdorff $H$ takes the same value on infinitely many distinct points of $U$.
   Hence no weakly sorting map can exist and $U$ is not weakly sorted. 
\end{proof}

\begin{cor}
   If $M$ is a weakly sortable homogeneous manifold, then $NH(x)$ is discrete for each $x\in M$ .
\end{cor}

Note that the converse does not hold: Example \ref{ex:FatS1} is homogeneous with discrete $NH(x)$, 
but has no weakly sorted neighborhood.
Also, the assumption on the restriction of $h$ being $1$-to-$1$ 
cannot be entirely lifted, see e.g. Example \ref{ex:BadPrufer} below.
We finish this subsection by proving some features of sorted neigborhoods.

\begin{lemma}
   \label{lemma:sorthomeo}
   Let $M$ be a manifold, $f:\wb{U}\to H$ be a sorting map, and $x\in U$. Then 
   for any open $V\subset U$ with $\wb{V}\cap U$ compact, the following hold.\\
   (a) $f(\wb{V}\cap U) = f(\wb{V})$.\\
   (b) There is a sorting map 
       $g:\wb{V}\to \wb{V}\cap U$ such that $g\upharpoonright \wb{V}\cap U$ is the identity.
\end{lemma}
In other words for (b): $\wb{V}\cap U$ is a retract of $\wb{V}$. 
\begin{proof}
   Notice that as a subspace of $U$, $\wb{V}\cap U$ is Hausdorff.\\
   (a) By Lemma \ref{lemma:cpctVU}, $\wb{V}-U \subset NH(\wb{V}\cap U)$.
   Lemma \ref{lemma:fxfy} then implies that $f(\wb{V}\cap U) = f(\wb{V})$.\\
   (b) A $1$-to-$1$ map from a compact to a Hausdorff space is a homeomorphism onto its image,
   hence by postcomposing with the inverse of 
   $f\upharpoonright \wb{V}\cap U$, we obtain the desired $g$.      
\end{proof}

\begin{prop}
    \label{prop:NHdisjoint}
    Let $V\subset U$ be open Hausdorff in the manifold $M$ 
    such that $\wb{V}\cap U$ is compact. Then the following is are equivalent.\\
    (a) $V$ is a sorted neighborhood;\\
    (b) $NH(x)\cap NH(y) = \varnothing$ for all $x,y\in \wb{V}\cap U$, $x\not=y$;\\
    (c) There is a sorting map $f:\wb{V}\to\wb{V}\cap U$ such that $f\upharpoonright \wb{V}\cap U$
        is the identity and for each compact Hausdorff $K\subset \wb{V}$,
        $f\upharpoonright K$ is a homeomorphism onto its image.
\end{prop}
This proposition explains our choice of vocabulary: 
the points in the closure of $\wb{V}\cap U$ are sorted into the disjoint families $\{x\}\cup NH(x)$ for 
$x\in \wb{V}\cap U$.

\begin{proof}  \ \\
   (a) $\Rightarrow$ (b) 
   Assume that $f:\wb{V}\to H$ is a sorting map.
   Then $f(NH(x)) = f(\{x\}) \not= f(\{y\}) = f(NH(y))$ whenever $x\not=y$ are points of $\wb{V}\cap U$.
   It follows that $NH(x) \cap NH(y) = \varnothing$.\\
   (b) $\Rightarrow$ (c)
   Assume that $NH(x) \cap NH(y) = \varnothing$ for each distinct $x,y\in \wb{V}\cap U$, and 
   set $f(u) = x$ for $x\in\wb{V}\cap U$ and $u\in NH(x)\cup\{x\}$.
   We show that $f$ is continuous.
   Let $u_i$ ($i\in\omega$) be a sequence in $\wb{V}$ with limit $u\in \wb{V}$.
   Recall that $\wb{V}-U\subset NH(\wb{V}\cap U)$ by Lemma \ref{lemma:cpctVU}.
   Let $x_i\in\wb{V}\cap U$ be such that either $x_i = u_i$ or $x_i\bumpeq u_i$.
   Let $x\in\wb{V}\cap U$ be a cluster point of the $x_i$. Then either  
   $x = u$ (if $u\in \wb{V}\cap U$) or $x\bumpeq u$ (if $u\in \wb{V}- U$)
   by Lemma \ref{lemma:NHseq}. It follows that $x$ is the unique cluster point
   of the $x_i$ and is thus their limit, and that 
   $f(u) = f(x) = \displaystyle\lim_{i\to+\infty} f(x_i) = \lim_{i\to+\infty} f(u_i)$.
   Hence, $f$ is continuous.
   \\
   Note that $f$ is $1$-to-$1$ on each Hausforff subset of $\wb{V}$
   since $f(x) = f(y)$ iff $x=y$ or $x\bumpeq y$. Hence, $f\upharpoonright K$ is a homeomorphism 
   onto its image for any
   compact Hausdorff $K\subset\wb{V}$.\\
   (c) $\Rightarrow$ (a) is immediate.
\end{proof}


\subsection{Rim-simplicity in dimension $1$ (mainly)}

Rim simplicity is one way of imposing that neighborhoods of two points $x\bumpeq y$ cannot 
intersect too erratically. The idea dates back to Hajicek \cite[Def. 3]{Hajicek:1971},
and is especially powerful (albeit somewhat tautologically) in dimension $1$.

\begin{defi}
  \label{defi:simple}
  Let $U$ be an open Hausdorff subset of a manifold $M$, $p\in\wb{U}-U$ and $n\in\omega$.
  Then $U$ is $n$-simple [resp. f-simple] 
  at $p$ iff given any open $V\ni p$, there is an open $W\subset V$ containing $p$
  such that $W\cap U$ is has $\le n$ many [resp. finitely many] connected components. 
  We say that $U$ is  
  $n$-rim-simple [resp. f-rim-simple]
  iff it is $n$-simple [resp. f-simple] at each point of its boundary $\wb{U}-U$.  
\end{defi}

We abbreviate ``$1$-simple'' as ``simple''.
Of course, in a Hausdorff manifold, not all open sets are f-rim-simple
(take $M$ to be the plane and $U$ to be $\{\langle x,y\rangle\,:\,y>\sin(1/x),\, x\not=0\}$),
but there is a basis of rim-simple open subsets (the regular open ones).
It is not the case for NH-manifolds, as seen immediately in the next example.
\begin{example}
   \label{ex:notrimsimple}
   No Hausdorff open subset of $\mathbf{G}(\R,\kappa)$ is rim-simple (for any $\kappa\ge 2$), but each point has
   a $2$-rim-simple neighborhood.
\end{example}

Actually, in dimension $1$ things are not really complicated\footnote{The next lemma 
even suggests that things are too simple. We apologize for the terrible pun.}.
\begin{lemma}
   \label{lemma:2simple}
   In a $1$-manifold, if an open subset $U$ is f-simple at $p\in\wb{U}-U$, then it is $2$-simple.
\end{lemma}
In this proof and the following ones of this subsection, we use the term ``interval'' as a shorthand
for ``subspace homeomorphic to an interval'' (open or closed). 
\begin{proof}
   Let $V\ni p$ be an open interval such that $U\cap V$ has finitely many components.
   Hence $U\cap V$ is a finite union of open intervals. 
   Since $p\in\wb{U}-U$, $p$ must be the endpoint of at least $1$ and at most $2$
   of these intervals, we may take $W\subset V$ to be the union of these to conclude.
\end{proof}

We note that rim-simplicity of an open Hausdorff set does not 
imply that things go well for its open subsets: you may have an Euclidean $U$ which is rim-simple, with
another Euclidean $V\subset U$ and $x\in V$ such that no open subset of $V$ containing $x$
is rim-simple (see Example \ref{ex:BadMoore}).
However, this cannot happen in dimension $1$.

\begin{lemma}
  \label{lemma:1rimsimple}
  Let $M$ be a $1$-manifold and $U\subset M$ be Hausdorff open such that $U$ is rim-simple
  [resp. f-rim-simple].
  Then each $x\in U$ has a local base $\mathcal{B}$ (of open subsets of $U$) such that 
  each member of $\mathcal{B}$ is rim-simple [resp. $2$-rim-simple]. 
\end{lemma}
\begin{proof}
  Let $x\in U$ and let $B\subset U$ be open and contain $x$.
  We may assume $B$ to be an interval. 
  Assume first that $U$ is rim-simple.
  In what follows, when $N\subset B$, then $\wb{N}^B$ denotes its closure in $B$ and $\wb{N}$ its closure in
  the whole space $M$.
  We want to find an open rim-simple subset of $B$ containing $x$.
  For that, fix an open interval $A\subset\wb{A}^B\subset B\subset U$, containing $x$.
  If $p\in(\wb{A}-A)\cap B$, then $p\in \wb{A}^B-A$ and $A$ is simple at $p$.
  Take $p\in \wb{A}-B$, then $p\bumpeq q$ for some $q\in A$. Since $U$ is Hausdorff and contains $q$,
  then $p\in \wb{U}-U$ as well. Fix an open interval $V\ni p$. Since $U$ is rim-simple there is another interval
  $W\subset V$ with $x\in W$ and $U\cap W$ connected.
  Hence $U\cap W$ is also an interval,
  and so is $A\cap W=A\cap U\cap W$.
  This shows that $A$ is rim-simple.\\
  The proof when $U$ is
  f-rim-simple is the almost the same, changing ``interval'' to ``finite union of intervals''
  (and ``connected'' to ``finite union of connected subspaces'')
  where relevant, and applying Lemma \ref{lemma:2simple} at the end.
\end{proof}

\begin{lemma}
   Let $M$ be a NH-$1$-manifold, and $p\in M$ be such that $NH(p)$ 
   contains three distinct points $a,b,c$. 
   Suppose that there are pairwise disjoint open Hausdorff subsets $U_a,U_b,U_c$ respectively containing $a,b,c$.
   Then $p\in\wb{U_q}-U_q$ for each $q\in\{a,b,c\}$ and $U_q$ is not f-simple at $p$
   for at least one $q\in\{a,b,c\}$.   
\end{lemma}
\begin{proof}
   The fact that $p\in\wb{U_q}-U_q$ for each $q\in\{a,b,c\}$ is obvious by Lemma \ref{lemma:sep}.
   Let $V\ni p$ be an open interval
   and $\varphi:V\to\R$ be an homeomorphism which sends $p$ to $0$.
   Since $a,b,c\in NH(p)$, any neighborhood of each point intersects any neighborhood of $p$ included in $V$.
   Since an interval has only two sides and we have three points,
   we may assume wlog that $U_a$ and $U_b$ both intersect $\varphi^{-1}\left((0,\frac{1}{n})\right)$ for each $n$.
   This implies that $U_a\cap V$ and $U_b\cap V$ have infinitely many components.
\end{proof}

\begin{cor}
   \label{cor:rimsimdis}
   Let $M$ be a $1$-manifold with a basis $\mathcal{B}$ of open sets such that
   each $U\in\mathcal{B}$ is f-rim-simple.
   Then $NH(p)$ is discrete for each $p\in M$. 
\end{cor}
\begin{proof}
   If $NH(p)$ is not discrete, by first countability it contains a converging sequence $y_n\to y$.
   Let $U\ni y$ be an open interval, then $U$ contains all but finitely many $y_n$.
   Choose three of them, take open pairwise disjoint subsets of $U$ containing them respectively and
   apply the previous lemma.
\end{proof}
Corollary \ref{cor:rimsimdis} does not hold in dimension $2$, as show by Example \ref{ex:BadPrufer}.

\begin{cor}
   Let $M$ be an homogeneous $1$-manifold. If there is an open Hausdorff $U\subset M$ 
   which is f-rim-simple,
   then $NH(p)$ is discrete for each $p\in M$. 
\end{cor}
\begin{proof}
   If $U$ is f-rim-simple at $p\in\wb{U}-U$, then by Lemmas \ref{lemma:2simple}--\ref{lemma:1rimsimple}
   each member of $U$ has a $2$-rim-simple local base.
   Hence by homogeneity $M$ has a basis of open $2$-rim-simple subsets.
\end{proof}

Say that a manifold is locally (f-)[$n$-]rim-simple iff each point has a local base of 
(f-)[$n$-]rim-simple open sets.
Example \ref{ex:FatS1} shows that a locally $2$-rim-simple
HNH-manifolds does not always have a weakly sorted neighborhood.
If one drops homogeneity, there are examples of NH-manifolds with a sorted neighborhood 
and a point without any local base of f-rim-simple open sets (Example \ref{ex:fan}), or 
with a local base of f-rim-simple but no local base of $n$-rim-simple open sets for any $n$
(Example \ref{ex:*}).
But we were unable to answer the following question.

\begin{ques}
   Let $M$ be a HNH-manifold and $f:\wb{U}\to U$ be a (weakly) sorting map.
   Do points in $U$ have a local base of f-rim-simple open sets~?
\end{ques}


\subsection{Homeomorphisms fixing $x$ while shuffling $NH(x)$} 

\begin{prop}
   \label{prop:transhom}
   Let $M$ be a NH-$1$-manifold and $x\in M$ such that for each $y_0,y_1\in NH(x)$ there is
   a homeomorphism $h$ of $M$ with $h(x) = x$, $h(y_0)=y_1$.
   Then $NH(x)$ is discrete (if non-empty).
\end{prop}
This proposition does not hold in dimension $2$, see Example \ref{ex:BadPrufer}.
First, we prove a lemma about subsets of $\R$.
\begin{lemma}
    \label{lemma:ndRdiscrete}
    Let $A\subset\R$ be closed such that for each 
    $x,y\in A$, there is an interval $(a,b)\ni x$ and $h:(a,b)\to\R$
    which is a homeomorphism onto its image, such that $h(A\cap(a,b))\subset A$ and 
    $h(x) = y$.
    Then either $A=\R$ or $A$ is discrete.
\end{lemma}
\begin{proof}
    Assume $A\not=\R$. If $A$ contains a maximal interval of the form $[u,v]$, $(-\infty,v]$ or $[u,+\infty)$,
    there cannot be any local homeomorphism as in the assumptions that sends 
    an interior point 
    to a boundary point; hence $A$ is nowhere dense. Now,
    since $A$ is closed, its complement is a countable union of maximal intervals $(u_n,v_n)$, $n\in\omega$.
    These intervals are disjoint. By maximality,
    $u_n,v_n\in A$ for each $n$. 
    Given $x\in A$, take $h$ as in the statement sending $x$ to $u_0$.
    Then there is $\epsilon>0$ such that
    either $[x,x+\epsilon)\cap A = \{x\}$ or $(x-\epsilon,x]\cap A = \{x\}$.
    It follows that $A = \cup_{n\in\omega}\{u_n,v_n\}$.
    Any closed countable subset of $\R$ is homeomorphic to a countable ordinal and hence has an isolated point
    $y$.
    Taking local homeomorphisms sending any other point to $y$ while preserving $A$ yields that each point  
    of $A$ is isolated.
\end{proof}

\begin{proof}[Proof of Proposition \ref{prop:transhom}]
   If a homeomorphism $h:M\to M$ fixes $x$, then $h(NH(x)) =  NH(x)$
   by Lemma \ref{lemma:sep} (c).
   Fix an interval $U\subset M$,
   then $NH(x)\cap U$ is closed nowhere dense in $U$.
   For each $y_0,y_1\in U\cap NH(x)$, take $h$ as in the statement. 
   By continuity, there is an interval around $y_0$ contained in $h^{-1}(U)$.    
   Hence, $A = NH(x)\cap U$ and $U\simeq\R$ satisfy the assumptions of Lemma \ref{lemma:ndRdiscrete},
   so $NH(x)\cap U$ is discrete.
\end{proof}


\section{ENH-manifolds and hereditary separability}\label{sec:conj}

Let us be a bit more precise than in the introduction and separate 
Problem \ref{prob:intro1} in two distinct questions, to which we add two natural variations.

\begin{ques}\ \\
   (a) Does an HNH-manifold contain a closed discrete uncountable subset~?\\
   (b) Does an HNH-manifold contain a discrete uncountable subset~?\\
   (c) Does an ENH-manifold contain a closed discrete uncountable subset~?\\
   (d) Does an ENH-manifold contain a discrete uncountable subset~?
   \label{ques:cdu}
\end{ques}

Let us first make a simple remark in the form of a lemma, whose proof is almost immediate.

\begin{lemma}
   Let $X$ be a locally hereditarily separable space. 
   Then $X$ is hereditarily separable iff $X$ contains no uncountable
   discrete subspace.
\end{lemma}
\begin{proof} 
   Suppose that $X$ is not hereditarily separable. Then
   $X$ contains a left separated space, i.e. a subspace $Y = \{y_\alpha\,:\,\alpha<\omega_1\}$
   such that $\wb{ \{y_\beta\,:\,\beta<\alpha\}}\cap\{y_\beta\,:\,\beta\ge\alpha\}=\varnothing$.
   For each $\alpha\in\omega_1$ choose an open hereditarily separable $U_\alpha\ni y_\alpha$, 
   then 
   $U_\alpha\cap Y$ is countable. 
   Let $f:\omega_1\to\omega_1$ be such that 
   $\cup_{\gamma<\alpha}U_\gamma \cap\{y_\beta\,:\,\beta\ge f(\alpha)\}=\varnothing$.
   Then $\{y_{f(\alpha)}\,:\,\alpha\in\omega_1\}$ is discrete in $X$.
\end{proof}

\begin{cor}
   \label{cor:hersepdis}
   A manifold is hereditarily separable iff it contains no uncountable discrete subspace.
\end{cor}

So, Question \ref{ques:cdu} (b) and (d) actually ask for hereditarily separable HNH- or ENH-manifolds.
We have consistent answers in both the positive and negative side for (c) and (d).
\begin{thm}
   \label{thm:hersepall}
   \ \\
   (a) Under {\bf CH}, there is a hereditarily separable ENH-manifold.\\
   (b) Under $\diamondsuit$, there is an ENH-manifold which is not hereditarily separable but does not
       contain any uncountable closed discrete subspace.\\
   (c) Under {\bf MA$(\omega_1)$}, there are no hereditarily separable ENH-manifolds.
\end{thm}
Here, {\bf CH} is the continuum hypothesis, $\diamondsuit$ is Jensen's diamond axiom and 
{\bf MA$(\omega_1)$} is Martin's axiom for posets of cardinality $\le\omega_1$. 
This situation of having positive results under {\bf CH} or $\diamondsuit$ and negative ones
under {\bf MA$(\omega_1)$} is incredibly common and was a constant theme in the 70s and 80s.
Our results are actually based on constructions and theorems of these times.
The proof of (a) and (b) are done in 
Subsection \ref{subsec:hersep}
(see Example \ref{ex:hersep}), but we stress that it is a quick consequence of 
the existence of a (so-called) {\em Kunen line}.
In Subsection \ref{subsec:ovenrack}, we show that our construction cannot yield 
an homogeneous manifold (at least in dimension $1$). 
The proof of (c) is done in 
Subsection \ref{subsec:hersep2}.
But
first, let us 
show some simple lemmas which take care of the most obvious cases.

\begin{lemma}\label{lemma:disc}
   If $Z$ is a subset of a locally Hausdorff space such that
   $z_0\bumpeq z_1$ for each $z_0,z_1\in Z$, then $Z$ is closed discrete.
\end{lemma}
\proof
   Any Hausdorff open subset can contain at most one point of $z\in Z$.
\endproof

\begin{lemma}\label{lemma:cupHx}
   If $M$ is a NH-manifold and $F\subset M$ is closed such that 
   $F\not\subset \cup_{n\in\omega} H(x_n)$ for any choice of $x_n\in F$,
   then
   $F$ contains a closed discrete uncountable subspace.
\end{lemma}
\begin{proof}
   Choose $x_\alpha$ ($\alpha\in\omega_1$) by induction such that $F\ni x_\alpha\not\in \cup_{\beta<\alpha} H(x_\beta)$.
   Then $x_\alpha\bumpeq x_\beta$ for each $\alpha,\beta$.
\end{proof}

A simple example of a HNH-manifold satisfying the assumptions of Lemma \ref{lemma:cupHx}
is $\mathbf{G}(\R,\omega_1)$, but this space
trivially contains a lot of closed uncountable discrete subsets.

\subsection{Towel-rack manifolds}\label{subsec:ovenrack}

A particular case of ENH-manifolds with a H-maximal sorted neighborhood 
are the towel-rack manifolds.
They trivially have discrete $NH(x)$. 
In what follows, $\kappa>1$ is a fixed cardinal.
(Our main examples will have $\kappa = 2$, but we decided to 
present the more general construction.)
A simple way of defining a NH-manifold is to start with $M\times\kappa$ (where $M$ is a Hausdorff manifold).
The $0$-th floor $M\times\{0\}$ is just a copy of $M$;
to define a neighborhood of $\langle x,\alpha\rangle$,
start with an open set $U\subset M$ containing $x$, put it on the $0$-th floor and
``lift'' a closed and nowhere dense subset $K\subset U$ containing $x$ 
to the $\alpha$-th floor,
akin to a towel rack where there are parallel lines and water drops fall on the ground (the $0$-th floor). 
If done correctly, we indeed obtain a manifold, and
if $K$ is a singleton (or discrete) and we do it at each $x\in M$, we obtain exactly $\mathbf{G}(M,\kappa)$.
In this brief subsection, we show that in dimension $1$, $\mathbf{G}(M,\kappa)$.
is actually the only homogeneous such manifold. 
\\
Let us give a more precise definition.
A family $\mathcal{K}$ of closed subsets in a space 
$X$ is {\em locally fixed} iff for each $x\in X$ and $K_0,K_1\in\mathcal{K}$
with $x\in K_0\cap K_1$, there is an open $U\ni x$ such that $K_0\cap U = K_1\cap U$.
Notice that
if $\mathcal{K}$ is locally fixed, then the family containing each finite intersection of members
of $\mathcal{K}$ is also locally fixed. We thus assume that each locally fixed family is closed
under finite intersections.
Let $M$ be a Hausdorff $n$-manifold and
$\mathscr{K} = \{\mathcal{K}_\alpha\,:\,1\le\alpha<\kappa\}$ be
such that

(TR1) each $\mathcal{K}_\alpha$ is a locally fixed familly
of non-empty closed nowhere dense subsets of $M$;

(TR2) if $K_0\in\mathcal{K}_\alpha$, $K_1\in\mathcal{K}_\beta$
with $\alpha\not=\beta$, $\alpha,\beta\ge 1$, then $K_0\cap K_1$ is finite;

(TR3) each $K\in\mathcal{K}_\alpha$ is $0$-dimensional, that is:
has a basis of clopen sets (in the subspace topology).
\vskip .3cm
\noindent
In dimension $1$ (TR3) follows from (TR1): any closed nowhere dense subset
of a $1$-manifold is automatically $0$-dimensional, as well known.
The {\em towel-rack $n$-manifold of height $\kappa$ 
associated to $M,\mathscr{K}$} (or just{ \em towel-rack manifold}
if $M,\mathscr{K}$ are clear from the context) is the space
$N(M,\mathscr{K}) = M\times\{0\}\cup_{1\le\alpha<\kappa}\cup\mathcal{K}_\alpha\times\{\alpha\}$
such that 
a neighborhood base $\mathcal{B}_{x,\alpha}$ of $\langle x,\alpha\rangle$ is given by 
$$ \{U_{K}^\alpha\,:\,U\cap K\ni x,\, U\text{ open in }M,\,K\in\mathcal{K}_\alpha\},$$
where
\begin{equation}
   \label{eq:ovenrack}
   \tag{TR4}
   U_{K}^\alpha = (U-K)\times\{0\}\cup (K\cap U)\times\{\alpha\}.
\end{equation}

\begin{lemma}
   \label{lemma:TRmanifolds}
   Let $N=N(M,\mathscr{K})$ 
   be a towel-rack manifold of height $\kappa$ 
   defined as above. Then $N$ is indeed a NH-manifold, and the following hold.
   \\
   (a) The projection on the first coordinate $\pi:N\to M$ is a sorting map.\\
   (b) If $z=\langle x,\alpha\rangle\in N$, then
       $NH(z) = \pi^{-1}(\{x\})- \{z\}$ is discrete.\\
   (c) For each $K\in\mathcal{K}_\alpha$,
   the map $h_{K,\alpha}$ exchanging $\langle x,\alpha\rangle$ with $\langle x,0\rangle$
   for each $x\in K$ and the identity elswewhere is a homeomorphism,\\
   (d) $N$ is an ENH-manifold iff $\cup_{\alpha<\kappa}\cup\mathcal{K}_\alpha = M$.
\end{lemma}
\begin{proof}
   Notice that since $\mathcal{K}_\alpha$ is closed under finite intersections,
   then so is $\mathcal{B}_{x,\alpha}$. \\
   (a) This follows from (\ref{eq:ovenrack}).\\
   (b) Equation (\ref{eq:ovenrack}) and (TR1) imply that $\langle x,\alpha\rangle\bumpeq\langle x,\beta\rangle$
   whenever $\alpha\not=\beta$,
   while $\langle x,\alpha\rangle$ and $\langle y,\beta\rangle$ 
   can be separated whenever $x\not= y$ since $M$ is Hausdorff.
   Then, (\ref{eq:ovenrack}) imply that $NH(z)$ is discrete.
   \\
   (c)
   Let $x\in M,\alpha<\kappa$, $U$ be open in $M$ and $K\in\mathcal{K}_\alpha$ 
   with $x\in K\cap U$.
   There is some open $V\subset U$ 
   with compact closure in $U$ such that
   $V\cap K$ is compact.
   Let now $z = \langle y,\beta\rangle\in N$ and $L\in\mathcal{K}_\beta$ containing $y$.
   If $\beta\not=\alpha$, by (TR2) there is some open $W\subset M$ containing $y$
   such that $L\cap K\cap W\subset\{y\}$. Then, 
   $h_{K,\alpha}\upharpoonright W_L^\beta$ is the identity,
   where $W_L^\beta$ is defined by (\ref{eq:ovenrack}).
   If $\alpha=\beta$, then either $y\not\in K$, in which case 
   $h_{K,\alpha}\upharpoonright (M-K)_L^\beta$ is also the identity, or
   $y\in K$, in which case there is an open $W\subset M$ with $W\cap L=W\cap K$ by (TR1)
   and $h_{K,\alpha}$ interchanges $W_L^\alpha$ with $W_L^0$. 
   In each of these cases, $h_{K,\alpha}$ is continuous at $z$ by (\ref{eq:ovenrack}).
   Since $h_{K,\alpha}$ is its own inverse, it is a homeomorphism.\\
   Now, 
   Since the $0$-th floor is a copy of $M$, $N$ is a NH-manifold by (c). 
   Then (d) is a direct consequence of (b).
\end{proof}

It is immediate that if $N=N(M,\mathscr{K})$ is a
towel-rack manifold of height $\kappa$, then 
for any $A\subset \kappa$ containing $0$, the subspace $M\times A\cap N$ is open in $N$ and
is itself a towel-rack manifold (of height $\le\kappa$).
Notice also that $N_\alpha = M\times\{\alpha\}\cap N$ is closed and $NH(N_\alpha)\cap N_\alpha =\varnothing$.
Example \ref{ex:hersep} below is an hereditarily separable towel-rack $1$-manifold 
of height $2$ (built with {\bf CH}).
Our next proposition shows in particular that we may not obtain an homogeneous one this way. 

\begin{prop}
   \label{prop:ovenrack} 
   Let $N=N(M,\mathscr{K})$ be a
   $1$-dimensional homogeneous towel-rack manifold of height $\kappa$. Then $N$ 
   is homeomorphic to a submanifold of $\mathbf{G}(M,\kappa)$.
\end{prop}
\begin{proof}
   So, $N = M\times\{0\}\cup_{1\le\alpha<\kappa}\cup\mathcal{K}_\alpha\times\{\alpha\}$,
   where $M$ is $1$-dimensional.
   Fix $K\in\mathcal{K}_\alpha$. We show that $K$ is discrete. This implies
   that $N$ is a submanifold of $\mathbf{G}(M,\kappa)$: given any open $U\subset M$
   with compact closure, $U\cap K$ is finite, and we thus have the same topology as 
   $\mathbf{G}(M,\kappa)$, the only difference being that the $\alpha$-th floor of $N$ is not assumed to 
   contain $\langle x,\alpha\rangle$ for {\em each} $x\in M$.\\
   Let $y_0,y_1\in K$.
   By homogeneity there is a homeomorphism $h$ of $N$ sending 
   $\langle y_0,\alpha\rangle$ to $\langle y_1,\alpha\rangle$.
   Then there are open intervals $I\ni y_0,J\ni y_1$ such that $h( I_K^\alpha )\subset J_K^\alpha$.
   Since $\langle y_0,0\rangle\bumpeq \langle y_0,\alpha\rangle$,
   $h(\langle y_0,0\rangle) = \langle y_1,\beta\rangle$ with $\beta\not=\alpha$.
   Let $L\in\mathcal{K}_\beta$ be such that $y_1\in L$.
   Up to postcomposing with $h_{L,\beta}$ we may assume that $\beta = 0$.
   Shrinking 
   $I,J$ if necessary, we have that $h(I\times\{0\})\subset J\times\{0\}$.
   Take any $z\in I_K^\alpha$, then $h(\langle z,0\rangle)\bumpeq h(\langle z,\alpha\rangle)$,
   hence we must have $h(\langle z,\alpha\rangle)\in K\times\{\alpha\}$
   and thus $h(\langle z,0\rangle)\in K\times\{0\}$. 
   Taking a subspace of $M\times\{0\}$ homeomorphic to $\R$ and containing $I\cup J\cup K$,
   this shows that $K$ satisfies the assumptions of Lemma \ref{lemma:ndRdiscrete} and is thus
   discrete.   
\end{proof}

Let us end this subsection with a brutally honest remark: (TR2) is so restrictive that we do not know if
one really gains something by considering towel-rack manifolds 
of height $\kappa>2$.
Of course, (TR2) holds vacuously when $\kappa=2$, and our examples in the next subsection are such.
If $\kappa>2$, we really need (TR2) to show that $h_{K,\alpha}$ is continuous
(and hence a homeomorphism),
and were unable to prove Proposition \ref{prop:ovenrack} without this fact, which is the reason
why we chosed to keep (TR2).
\\
Another strong restriction is the fact that in (\ref{eq:ovenrack}), only the $0$-th and $\alpha$-th floor 
are involved. It would probably be more interesting if we allow neighborhoods of points in the $\alpha$-th
floor to contain points in other floors than just the $0$-th. We did not find any convincing set
of rules that enables this without sacrificing too much.


\subsection{A hereditarily separable ENH-manifold under {\bf CH}}\label{subsec:hersep}

Our goal in this subsection is given in its title.
It is a straightforward use of the famous ``Kunen line'' construction, which dates back to
the 70s, and appeared in print in \cite[p. 1000]{JuhaszKunenRudin:1976}\footnote{Actually,
the description starts on p. 999, but the author never had the
opportunity to cite the page 1000 of anything.}.
We recall the construction, because we are going to use a slight
modification in a further example. 
We recall that Hausdorff manifolds have cardinality $\mathfrak{c}$, 
(see e.g. \cite[Proposition 1.8]{GauldBook}),
hence $\omega_1$
under {\bf CH}.
We start the construction with a Hausdorff manifold $X$
with original topology $\rho$.
Let $\{x_\alpha\,:\,\alpha\in\omega_1\}$ be a $1$-to-$1$ enumeration of $X$
and write $X_\alpha$ for $\{x_\beta\,:\,\beta<\alpha\}$.
We shall define a refinement $\tau$ of $\rho$ on $X=X_{\omega_1}$.
For this, we let $\{H_\mu\,:\,\mu<\omega_1\}$ be an enumeration of all
countable subsets of $X$, so that $H_\mu\subset X_\mu$.
We are going to define the topology $\tau_\alpha$ on $X_\alpha$ by induction, such that the following holds
for each $\xi<\eta\le\omega_1$. (We write $\text{Cl}_\rho(A)$ for the $\rho$-closure of $A$,
and correspondingly for $\tau,\tau_\alpha$.)
\vskip .3cm
(K1) $\tau_\xi = \tau_\eta\cap\mathcal{P}(X_\xi)$, that is: 
           $\tau_\xi$ is the topology induced by $\tau_\eta$
           on $X_\xi$.

(K2) $\tau_\eta$ is first countable, locally compact and Hausdorff.

(K3) $\tau_\eta\supset\rho_\eta$, where $\rho_\eta$ is the topology induced by $\rho$ on $X_\eta$.

(K4) For each $\nu\le\xi$, if $x_\xi\in\text{Cl}_\rho(H_\nu)$, then 
           $x_\xi\in\text{Cl}_{\tau_\eta}(H_\nu)$.

(K5) There is a countable $\tau_\eta$-compact subset $K_\xi$ of $X_\eta$ such that 
     a base of $\tau_\eta$-neighborhood of $x_\xi$ is given by $\{U_n\cap K_\xi\,:\,n\in\omega\}$, 
           where $U_n$ is a $\rho$-base of $\rho$-open neighborhoods in $X$.
\vskip .3cm\noindent
Notice that sets that are compact in a topology are compact in any weaker topology in which they remain Hausdorff,
hence by (K1), $K_\xi$ is $\tau$- and $\rho$-compact for each $\xi$.
In (K4), notice that $H_\nu\not\ni x_\xi$ when $\nu\le\xi$, since $H_\nu\subset X_\nu$.
The induction starts by letting $\tau_\beta$ be discrete for $\beta\le\omega$.
We now assume that the construction is done below $\beta$, with (1)--(5) holding for $\xi<\eta\le\beta$.
If $\beta$ is limit, take $\tau_\beta$ to be the union of those $U\subset X_\beta$
whose intersection with $X_\eta$ is in $\tau_\eta$ for each $\eta<\beta$.
Then (K1)--(K5) follow immediately.
\\
Let now $\beta = \alpha + 1$. Observe that the conditions (K1)--(K4) imply that $\tau_\alpha$ is actually
regular and $0$-dimensional (i.e. has a base of clopen sets).
We now define $\tau_\beta$ on $X_\beta = X_\alpha\cup\{x_\alpha\}$.
The only difficulty is to make (K4) and (K5) hold for $\xi = \alpha$.
Let $E(\alpha) = \{\nu\le\alpha\,:\,x_\alpha \in\text{Cl}_\rho(H_\nu)\}$.
If $E(\alpha)=\varnothing$, we let $\tau_\beta$ be the topology with base $\tau_\alpha\cup\{x_\alpha\}$.
Otherwise, we let $\langle\nu_n\,:\,n\in\omega\rangle$ enumerate $E(\alpha)$ with each member listed infinitely
many times.
In the $\rho$-topology, $X$ is a Hausdorff manifold and we may choose
a nested $\rho$-open base $\{U_n\,:\,n\in\omega\}$ at $x_\alpha$.
Let $\gamma(n)<\alpha$ be such that $x_{\gamma(n)}\in H_{\nu_n}\cap U_n$. Then 
$\{x_{\gamma(n)}\,:\,n\in\omega\}$ is $\rho$-discrete in $X_\alpha$,
and is thus also $\tau_\alpha$-discrete. Choose $\rho$-open $V_n\ni x_{\gamma(n)}$ 
which are pairwise disjoint in $X_\alpha$. By shrinking further each $V_n$ if necessary,
by (K5) we may assume that $K_{\gamma(n)}\cap V_n$ is compact open in $X_\alpha$.
Define then $\tau_\beta$ to have as a base the sets of the form $\{x_\alpha\}\cup \cup_{m>n} V_n\cap K_{\gamma(m)}$
for $n\in\omega$, together with $\tau_\alpha$.
Then (K1)--(K4) follow easily, and for (K5) 
we may define $K_\alpha$ as $\{x_\alpha\}\cup \cup_{n\in\omega} V_n\cap K_{\gamma(n)}$.
We let $\tau =\tau_{\omega_1}$ on $X=X_{\omega_1}$. Then 
$\tau$ is first countable, non-Lindel\"of and locally compact as $K_\alpha$ is
a compact countable neighborhood of $x_\alpha$ for each $\alpha<\omega_1$.
It is then easy to prove the following property of $\tau$. 
\begin{equation}
    \tag{*}
    \label{eq:closurecount}
    \begin{array}{l} \text{If $A\subset X$ is $\rho$-separable, then Cl$_\rho(A) - $
          Cl$_\tau(A)$}\\
          \text{is at most countable, and $A$ is $\tau$-separable.} 
    \end{array}
\end{equation}
Indeed, let $A_0\subset A$ be countable and $\rho$-dense in $A$.
Then $A_0 = H_\nu$ for some $\nu$.
By condition (K4), if $x_\xi\in\text{Cl}_\rho(A_0) = \text{Cl}_\rho(A)$, 
and $\xi>\nu$, then $x_\xi\in\text{Cl}_\tau(A_0)$.
Hence $\text{Cl}_\rho(A)-\text{Cl}_\tau(A)\subset X_\nu$.
It follows that $A_0\cup \left(A\cap X_\nu\right)$ is $\tau$-dense in $A$.
\\
Hence, $\langle X,\tau\rangle$ is hereditarily separable whenever $\langle X,\rho\rangle$ is.
Taking $X=\R$, this yields the following.

\begin{example}[{\bf CH}]
   \label{ex:hersep}
   There is a hereditarily separable ENH-$1$-manifold with a H-maximal sorted neighborhood.
\end{example}
\begin{proof}[Details]
We let $\R = X = \{x_\alpha\,:\,\alpha<\omega_1\}$ and $K_\alpha$ be defined as above.
Let $\mathcal{K}=\{K_\alpha\,:\,\alpha<\omega_1\}$,
then $\mathcal{K}$ is locally fixed
by (K1) and (K5). We let 
$\mathcal{K}_1$ contain the finite intersections of members of $\mathcal{K}$ and
$M$ be the towel-rack manifold associated to $\R$, $\{\mathcal{K}_1\}$. 
Notice that a neighborhood base for $\langle x_\alpha,1\rangle$ is given by
$$ 
   N(x,\epsilon) = \left((x_\alpha-\epsilon,x_\alpha+\epsilon)\cap K_\alpha\right)\times\{1\}\cup
   \left((x_\alpha-\epsilon,x_\alpha+\epsilon)- K_\alpha\right)\times\{0\},
$$
the subspace topology on $\R\times\{1\}$ is $\tau$, thus $M$ is hereditarily separable.
Moreover, $\R\times\{0\}$ is CH- and H-maximal in $M$.
\end{proof}

We could have started with any hereditarily separable Hausdorff manifold $X$ instead of $\R$,
because countable compact subspaces of manifolds are $0$-dimensional.
Hence, there are examples in any dimension.
We now present a simple modification of the above construction 
to show that there are manifolds for which, in Question \ref{ques:cdu}, 
the answer of (c) is negative but that of (d) is positive: 
ENH-manifolds with uncountable discrete subsets but no closed one.
We will perform the Kunen line construction with $X=\LL_+$ (the long ray),
so let us recall some facts about this space (see e.g. \cite[Section 1.2]{GauldBook} for proofs).
Let $\LL_{\ge 0}$ be
the set $\omega_1\times[0,1)$ with
the lexicographic order topology. 
When convenient, we identify $\alpha\in\omega_1$ with $\langle\alpha,0\rangle\in\LL_{\ge 0}$,
which explains what we mean when we
write for instance $\alpha\in\LL_{\ge 0}$ or an interval $[0,\alpha]\subset\LL_{\ge 0}$ for an ordinal $\alpha$.
Then, $[0,\alpha] \subset\LL_{\ge 0}$ is homeomorphic to $[0,1]$ for each $\alpha$,
$\LL_{\ge 0}$ is countably compact and
any closed unbounded set in $\LL_{\ge 0}$ intersected with $\omega_1\subset\LL_{\ge 0}$
is closed unbounded (in both $\omega_1$ and $\LL_{\ge 0}$).
We then let $\LL_+$ be $\LL_{\ge 0}-\{0\}$. 
Notice that $\LL_+$ contains uncountable discrete subsets (e.g. the successor ordinals), but no
closed one.
We shall use the additional axiom $\clubsuit$.
\begin{defi}[{\bf Axiom $\clubsuit$}] 
   \label{def:clubsuit}
   For each limit $\alpha<\omega_1$ there is $L_\alpha\subset\alpha$,
   called a $\clubsuit$-sequence, such that:
   \\
   (1) $L_\alpha$ has order type $\omega$ and is cofinal in $\alpha$;\\
   (2) for any uncountable $A\subset\omega_1$, there is a stationary subset $S\subset\omega_1$
   such that $L_\alpha\subset A$ for $\alpha\in S$.
\end{defi}
The original definition of $\clubsuit$ is formally weaker, but equivalent, to this one;
recall also that {\bf CH} + $\clubsuit$ is equivalent to $\diamondsuit$, see e.g. 
\cite[Definition 21.8, Lemma 22.9 and Theorem 22.10]{JustWeeseII}.
We recall that a sequence $L_\alpha$ which satisifies (1) is called a {\em ladder system} (on $\omega_1$).
We now assume {\bf CH} + $\clubsuit$ (which, again, is equivalent to $\diamondsuit$), 
and let $L_\alpha$ be a $\clubsuit$-sequence.
We let again
$\{x_\alpha\,:\,\alpha\in\omega_1\}$ be an enumeration
of $\LL_+$ used in the Kunen line construction and write $X_\alpha = \{x_\beta\,:\,\beta<\alpha\}$.
We may assume that $x_\alpha\in(0,\alpha)\subset\LL_+$ for each $\alpha>0$.
Set $G_\alpha = \{x_\alpha\,:\,\alpha\in L_\alpha\}$ for limit $\alpha$.
When taking our sequence of countable subsets $H_\alpha$, we now assume that
\begin{equation}
   \label{eq:doublestar}\tag{**}
   H_\alpha = G_\alpha\text{ when $\alpha$ is limit.}
\end{equation}
(Notice that in the construction of the $\tau$-topology, there is no need for the $H_\alpha$ to be distinct, we only need
each countable subset of $\LL_+$ be included.)
Let thus $\tau$ be the topology on $\LL_+$ obtained as above.

\begin{lemma}
   If (\ref{eq:doublestar}) holds, then 
   $\langle \LL_+,\tau\rangle$ 
   has no closed uncountable discrete subset.   
\end{lemma}
\begin{proof}
   Let $A\subset\LL_+$ be uncountable,
   we show that it is either not $\tau$-closed or not $\tau$-discrete.
   We denote the order on $\LL_+$ by $\ge^\ell$.
   If there is some $\alpha$ such that
   $A$ is included in $(0,\alpha)$ (the interval is taken in $\LL_+$) which is $\rho$-hereditarily separable,
   then by (\ref{eq:closurecount})
   $A$ is $\tau$-separable and thus not $\tau$-discrete.
   We may thus assume that $A$ is $\ge^\ell$-unbounded in $\LL_+$.
   \\
   Let $\xi_0\in\omega_1$ be minimal such that $x_{\xi_0}\in A$.
   For each $\alpha\in\omega_1$, choose $\xi_\alpha$ such that $x_{\xi_\alpha}\in A$, and
   
   (a) $\xi_\alpha > \sup_{\beta<\alpha}\xi_\beta$,
   
   (b)
   $x_{\xi_\alpha} \ge^\ell \sup^\ell_{\beta<\alpha}x_{\xi_\beta}$,
   
   (c) $x_{\xi_\alpha} \ge^\ell\alpha$.
   \\
   Here, $\sup^\ell$ is the supremum in $\LL_+$ according to $\ge^\ell$.
   Thus, when $\alpha$ grows, $\xi_\alpha$ grows in $\omega_1$ and $x_{\xi_\alpha}$
   grows in $\LL_{\ge 0}$ for $\ge^\ell$.
   A classical leapfrog argument (using e.g. \cite[Corollary 21.5]{JustWeeseII}) shows that
   there is a closed unbounded $D\subset\omega_1$ such that for each $\alpha\in D$,
   $ \alpha =\sup_{\beta<\alpha}\xi_\beta = \sup^\ell_{\beta<\alpha} x_{\xi_\beta} = \sup^\ell_{\beta<\alpha}x_\beta$.
   Hence, 
   seen as a member of $\LL_+$, $\alpha$ is equal to $x_\xi$ for some $\xi\ge\alpha$.
   Take $\alpha\in D$ such that $L_\alpha\subset \{\xi_\beta\,:\,\beta<\omega_1\}$. 
   Thus $G_\alpha=\{x_\beta\,:\,\beta\in L_\alpha\}$ 
   is a subset of $\{x_{\xi_\beta}\,:\,\beta<\omega_1\}\subset A$ and has $\rho$-limit point $\alpha = x_\xi$.
   But by (K4) and (\ref{eq:doublestar}), $x_\xi\in\text{Cl}_\tau(G_\alpha)$. Hence
   $x_\xi$ is a $\tau$-cluster point of $A$ in $\LL_{\ge 0}$. 
   It follows that $A$ is either non-closed or non-discrete for $\tau$.   
\end{proof}

\begin{example}[$\diamondsuit$]
   \label{ex:nocloseddis}
   There is an ENH-$1$-manifold $M$
   which contains uncountable discrete subsets
   but none that is closed. Moreover, $M$ contains a H-maximal sorted neighborhood. 
\end{example}
\begin{proof}[Details]
Define $M$ as $\LL_+\times\{0,1\}$, with the same topology as in Example \ref{ex:hersep}.
As said above, there are uncountable discrete subsets in $\LL_+$, and hence in $M$ as well. 
\end{proof}

\begin{ques}
   \label{q:homkun}
   \
   \\
   (a)
   Is there an homogeneous Kunen line and, if yes, can we obtain an homogeneous manifold~?\\
   (b) 
   Is {\bf CH} necessary to obtain an hereditarily separable ENH-manifold as Example \ref{ex:hersep}~?\\
   (c) Is 
   $\diamondsuit$ necessary to obtain an ENH-manifold with the properties of Example \ref{ex:nocloseddis}~?
\end{ques}
Notice that an homogeneous Kunen line cannot be locally countable (because there would be isolated points),
so what we ask in (a) is whether there is a homogeneous hereditarily separable strenghtening of the 
topology on some $\R^n$ 
such that neigborhoods of points are compact nowhere dense sets in the original topology.
(Although such a space would be very far from being a line when $n>1$,
we believe that the original name would still be adequate.)
Notice also that homogeneity of such a Kunen line would not immediately transfer to that of its associated manifold. 
In fact, we saw with Proposition \ref{prop:ovenrack}
that it is not possible to obtain
a homogeneous $1$-manifold with the method of 
Example \ref{ex:hersep}.
This does not a priori precludes possibilities in higher dimensions. 

\subsection{No hereditarily separable ENH-manifolds under {\bf MA($\omega_1$)}}\label{subsec:hersep2}

In this subsection, we show that {\bf MA($\omega_1$)} implies that there are no 
hereditarily separable ENH-manifolds.
Since {\bf MA($\omega_1$)} negates {\bf CH}, it gives a partial answer to Question \ref{q:homkun}.
Let $\mathscr{U} = \langle U_\alpha\,:\,\alpha\in\omega_1\rangle$ be an 
$\omega_1$-sequence of open sets of a given manifold $M$.
We say that $\mathscr{U}$ is {\em step-by-step dense} iff the following three conditions hold:
\vskip .3cm
(SbSD0) $U_0$ is H-maximal in $M$;

(SbSD1) $U_\alpha\subset U_\beta$ whenever $\alpha<\beta$;

(SbSD2) $U_{\alpha+1} - U_\alpha$ is dense in $M-U_\alpha$ (in the induced topology) for each $\alpha<\omega_1$.

(SbSD3) $U_{\alpha} - U_0$ is Lindel\"of for each $\alpha<\omega_1$.
\vskip .3cm
\begin{lemma}
   \label{lemma:SbSDa}
   Let $M$ be an ENH-manifold. Then either (a) or (b) below holds.\\
   (a) There is subspace $V$ of $M$ which is Hausdorff, non-Lindel\"of and locally compact.\\
   (b) There is a step-by-step dense $\mathscr{U} = \langle U_\alpha\,:\,\alpha\in\omega_1\rangle$.
\end{lemma}
\begin{proof}
   Set $U_0$ to be H-maximal in $M$ and let $B_0 = M-U_0\not=\varnothing$.
   We define $U_\alpha$ by induction on $\alpha$ such that $U_\alpha-U_0$ is Lindel\"of.
   If $\alpha$ is limit, set $U_\alpha = \cup_{\beta<\alpha}U_\beta$, then 
   $U_\alpha-U_0$ is Lindel\"of by induction.
   Given $U_\alpha$, set $B_\alpha = M-U_\alpha$.
   Then $M-U_\alpha$ is non-empty by Theorem \ref{thm:ClH} (because $U_\alpha-U_0$ can be covered by countably
   many Euclidean subsets).
   By Lemma \ref{lemma:locHaus} (c), there is 
   $V_{\alpha+1}\subset B_\alpha$ which is Hausdorff, open and dense in $B_\alpha$ in the induced topology.
   Then $U_{\alpha+1} = V_{\alpha+1}\cup U_\alpha$ is open in $M$.
   If $V_{\alpha+1}$ is non-Lindel\"of, then by Lemma \ref{lemma:openclosedloccpct}
   $V_{\alpha+1} = U_{\alpha+1}\cap B_\alpha$ is locally compact, and hence (a) holds
   and the proof is complete.
   If $V_{\alpha+1}$ is Lindel\"of then so is $U_{\alpha+1}-U_0$, hence $B_{\alpha+1}=M-U_{\alpha+1}\not=\varnothing$
   by Theorem \ref{thm:ClH}, and we may proceed.
   If $U_\alpha-U_0$ is Lindel\"of for each $\alpha$, then $\mathscr{U} = \langle U_\alpha\,:\,\alpha\in\omega_1\rangle$
   satisfies (SbSD0) to (SbSD3), and (b) holds.
\end{proof}

\begin{lemma}
   \label{lemma:SbSDb}
   If there is an ENH-manifold $M$ which contains a step-by-step dense 
   $\mathscr{U} = \langle U_\alpha\,:\,\alpha\in\omega_1\rangle$,
   then $\text{cov}(\mathcal{B}) = \omega_1$.
\end{lemma}
\begin{proof}
   By Lindel\"ofness and local compactness, $U_\alpha-U_0$ may be covered by countably many compact sets, hence
   by Lemma \ref{cor:quasi-compact} $NH(U_\alpha-U_0)\cap U_0$ is a countable union of nowhere dense sets.
   Notice that $NH(U_\alpha-U_0)\cap U_0 = NH(U_\alpha)\cap U_0$ since $U_0$ is open Hausdorff.
   If $\text{cov}(\mathcal{B}) > \omega_1$, then there is some $x\in U_0 - \cup_{\alpha<\omega_1}NH(U_\alpha)$.
   Let $y\bumpeq x$
   (which exists since $M$ is everywhere non-Hausdorff), 
   then $y\not\in\cup_{\alpha<\omega_1}U_\alpha$.
   Fix an Euclidean $V\ni y$. Since $y\in M-U_{\alpha}$ for each $\alpha$,
   by (SbSD2) $V\cap U_{\alpha+1}-U_\alpha\not=\varnothing$. 
   It follows that $\langle V\cap U_\alpha\,:\,\alpha\in\omega_1\rangle$
   is an $\omega_1$-sequence of open subsets of $V$ which is strictly increasing for the inclusion.
   This is impossible since $V$ is Euclidean and thus hereditarily Lindel\"of.
\end{proof}

Recall that an {\em S-space} is a regular hereditarily separable non-Lindel\"of space.
None can be constructed in {\bf ZFC}, for instance
the proper forcing axiom {\bf PFA} impedes their existence 
and {\bf MA$(\omega_1)$} that of locally compact ones \cite{Szentmiklossy:1980},
while a Kunen line is an example under {\bf CH}.
See e.g. \cite{Roitman:1984} (in particular just after 6.4) for an account of these classical results,
and \cite{Tall:PFAforthemasses} for similar results on more recent models of set theory.
Recall also that {\bf MA$(\omega_1)$} implies that $\text{cov}(\mathcal{B}) > \omega_1$.

\begin{thm}
   \label{thm:nohersep}
   If there is a hereditarily separable ENH-manifold, then either $\text{cov}(\mathcal{B}) = \omega_1$
   or there is a first countable locally compact S-space.
\end{thm}
\begin{proof}
   Let $M$ be a hereditarily separable ENH-manifold.
   Then either (a) or (b) of Lemma \ref{lemma:SbSDa} holds.
   In the first case,
   $V$ given by (a) is a first countable locally compact S-space.
   If (b) holds, then $\text{cov}(\mathcal{B}) = \omega_1$ by Lemma \ref{lemma:SbSDa}.
\end{proof}

\begin{cor}[{\bf MA$(\omega_1)$}]
   There are no hereditarily separable ENH-manifolds.
\end{cor}

Notice that (consistently) there are manifolds with a step-by-step dense 
$\mathscr{U} = \langle U_\alpha\,:\,\alpha\in\omega_1\rangle$.
Indeed, if $M$ is the hereditarily separable ENH-manifold of Example \ref{ex:hersep},
then its first floor is Hausdorff and
we may choose by induction $U_\alpha$ such that $U_{\alpha+1}-U_\alpha$ contains 
a countable dense subspace of $M-U_\alpha$.
It seems more difficult to exhibit an example such that $U_{\alpha+1}-U_\alpha$ is H-maximal 
in $M-U_\alpha$, as in the proof of Lemma \ref{lemma:SbSDa}.
We believe that a variant of Example \ref{ex:hersep}
with $\omega_1$-many floors could yield such a construction (under {\bf CH}), but 
did not really try to check whether it works, and the case is still open.


\section{Quasi-countably compact subsets of NH-manifolds}\label{sec:quasi-ctbly-cpct}

In this section, we prove two theorems and exhibit two examples involving quasi-countable compactness.

\subsection{There are no quasi-countably compact ENH-manifolds}\label{subsec:noENHqcc}

In a previous version of these notes, we wrote (more or less)
``our gut feeling is that there is a very simple argument lurking somewhere which implies that there are no
quasi-countably compact ENH-manifolds''.
It appears that our gut feeling was right.\footnote{This could be a reason of personal pride, but the said pride
is kind of tarnished by the fact that it took us so long to find the argument in question.}

\begin{thm}
   \label{thm:noENHqcc}
   Let $X$ be first countable and locally Hausdorff 
   such that $NH_X(X)$ is dense in $X$.
   Then $X$ is not quasi-countably compact.
\end{thm}
Notice that if $X$ is everywhere non-Hausdorff, then $NH_X(X)=X$.
In passing, this theorem shows that Question \ref{ques:cdu} (a) or (c) has a positive answer
if one changes ``uncountable'' to ``infinite''.
The proof makes repeated uses of Lemmas \ref{lemma:twolimits} and \ref{lemma:NHseq}.
\begin{proof}
   Suppose that $X$ is quasi-countably compact.
   \begin{claim}
      \label{claim:qcc1}
      $X$ is everywhere non-Hausdorff.
   \end{claim}
   \begin{proof}[Proof of the claim]
      By definition, $NH_X(X)=\{x\in X\,:\,NH(x)\not=\varnothing\}$.
      Let $x\in X$ and $x_n\in X$ be a sequence converging to $x$ such that $NH(x_n)\not=\varnothing$.
      Take $y_n\bumpeq x_n$.
      By quasi-countable compactness, there is a cluster point $y$ of the $y_n$,
      then $y\bumpeq x$ by Lemma \ref{lemma:NHseq}, hence $NH(x)\not=\varnothing$.
   \end{proof}   
   The proof now basically consists in taking subsequences of limits of subsequences of limits (etc).
   Let us first present a semi-formal sketch of the argument.
   We start with two points $x\bumpeq y^{(0)}$ and take a sequence $x_n$ converging to both of them
   (using Lemma \ref{lemma:twolimits}). Choose then $y_n^{(1)}\bumpeq x_n$ for each $n$.
   Since $X$ is quasi-countably compact, there is a subsequence of the $y_n^{(1)}$  
   converging to some $y^{(1)}$. By Lemma \ref{lemma:NHseq}, $x\bumpeq y^{(1)} \bumpeq y^{(0)}$.
   For each $n$, choose a sequence $x_{n,m}$ converging to both $x_n$ and $y^{(1)}_n$,
   and then choose $y_{n,m}^{(2)}\bumpeq x_{n,m}$.
   Some subsequence of $y_{n,m}^{(2)}$ converges to some $y_{n}^{(2)}$, and
   a subsequence of the $y_{n}^{(2)}$ to some $y^{(2)}$.
   Applying Lemma \ref{lemma:NHseq} twice, we first see that 
   $x_n\bumpeq y_{n}^{(2)}\bumpeq y_{n}^{(1)}$, and then that $x\bumpeq y^{(2)} \bumpeq y^{(1)}$
   and $y^{(2)} \bumpeq y^{(0)}$. Proceeding by induction, we end up defining
   $y^{(n)}$ for each $n\in\omega$ such that $y^{(n)}\bumpeq y^{(m)}$ for each $n\not= m$.
   Hence, $Y=\{y^{(n)}\,:\,n\in\omega\}$ is an infinite closed discrete subspace by 
   Lemma \ref{lemma:disc}, a contradiction.
   \\
   We now give a (much more tiresome) formal proof for those not convinced by the above lines. 
   To make things clearer, we first take some time to 
   explain the notations and conventions.
   The indices of our sequences belong to trees (in the set theoritical sense)
   $T^{(n)}$ ($n\in\omega$) which are subtrees of $\omega^{<n+1}$, that is,
   the sequences of integers of length $<n+1$ ordered by end-extension.
   If $\sigma\in\omega^{<n}$ and $k\in\omega$, then $\sigma k\in\omega^{<n+1}$ is the end-extension
   of $\sigma$ by $k$. Denote by $h(\sigma)$ the length of $\sigma$, and by 
   $\sigma\upharpoonright m$ the sequence containing the first $m$ entries of $\sigma$
   (with $\sigma\upharpoonright m=\sigma$ if $m\ge h(\sigma)$). If $A$ is a subtree of $\omega^{<\omega}$
   then $A\upharpoonright m = \{\sigma\upharpoonright m\,:\,\sigma\in A\}$.
   The trees $T^{(n)}$ have the following properties:
   \vskip .2cm
   (T1) Each $\sigma\in T^{(n)}$ with $h(\sigma)<n$ has infinitely many immediate successors,
       i.e. $\{k\in\omega\,:\,\sigma k\in T^{(n)}\}$ is infinite.
       
   (T2) If $m<n$, then $T^{(n)}\upharpoonright m \subset T^{(m)}$.
   \vskip .2cm
   \noindent
   We shall define $T^{(n+1)}$ by first taking all the $\sigma k$ ($k\in\omega$) 
   for all $\sigma$ of length $n$ in $T^{(n)}$,
   and then go down and ``prune'' some members of $T^{(n)}$ of length $m$, $m$ going from $n$ to $0$. That is:
   we delete some $\sigma$ and all their successors from $T^{(n+1)}$. We thus assume implicitely that
   whenever $\sigma$ is deleted from $T^{(n+1)}$, all its successors are deleted as well.
   \\
   At the $n$-th stage, we
   will define points $x_\sigma,y^{(n)}_\sigma$ for $\sigma\in T^{(n)}$.
   In particular, $n\ge h(\sigma)$.
   When these are defined, we let 
   $\overrightarrow{\mathbf{x}}_{n,\sigma}$ and $\overrightarrow{\mathbf{y}}^{(m)}_{n,\sigma}$
   (with $h(\sigma)\le m\le n$)
   respectively denote the sequences 
   $\langle x_{\sigma k}\,:\,k\in\omega\,,\sigma k \in T^{(n)}\rangle$ and
   $\langle y_{\sigma k}^{(m)}\,:\,k\in\omega\,,\sigma k \in T^{(n)}\rangle$ with $k$ increasing.
   By (T1) and (T2) above, 
   if $h(\sigma)\le n$ and $\sigma\in T^{(n)}\upharpoonright h(\sigma)$, then 
   $\overrightarrow{\mathbf{x}}_{n,\sigma}$ and
   $\overrightarrow{\mathbf{y}}_{n,\sigma}^{(m)}$ 
   are respectively (infinite) subsequences of 
   $\overrightarrow{\mathbf{x}}_{h(\sigma),\sigma}$ and 
   $\overrightarrow{\mathbf{y}}_{h(\sigma),\sigma}^{(m)}$ 
   and hence the former have the same limits as the latter if these do converge.
   \\
   We may now start the construction
   with two points $x_\varnothing \bumpeq y_\varnothing^{(0)}$
   and set $T^{(0)}=\varnothing$.
   Now, 
   fix $n$ and
   assume that $T^{(m)}$ satisfies (T1) and (T2) above for each $m\le n$,
   and that $x_\sigma,y_\sigma^{(m)}$ are defined for each $\sigma\in T^{(n)}$ and $h(\sigma)\le m\le n$, such that
   the following holds:
   \vskip .2cm
    (XY1) For all $i<j\le n$ with $h(\sigma)\le i$, $x_\sigma\bumpeq y_\sigma^{(i)}\bumpeq y_\sigma^{(j)}$,
    
    (XY2) $\overrightarrow{\mathbf{x}}_{m,\sigma}$ converges to both $x_\sigma$ and $y_\sigma^{(h(\sigma))}$
          when $h(\sigma)<m\le n$.
          
    (XY3) $\overrightarrow{\mathbf{y}}_{j,\sigma\upharpoonright i}^{(m)}$ converges to 
          $y_{\sigma\upharpoonright i}^{(m)}$ when $0\le i \le m \le j\le n$.
   \vskip .2cm
   \noindent
   We now define $T^{(n+1)}$ from the top, pruning away some points when we go down.
   First, set $T^{(n+1)}$ as $\{\sigma k\,:\,\sigma\in T^{(n)}, h(\sigma) = n, k\in\omega\}$.
   For each $\sigma\in T^{(n)}$ with $h(\sigma) = n$ and $k\in\omega$
   choose $x_{\sigma k}$ such that 
   $\overrightarrow{\mathbf{x}}_{n+1,\sigma}$ converges to both $x_\sigma$ and $y_\sigma^{(n)}$
   (which is possible by Lemma \ref{lemma:twolimits} and (XY1)).
   We now define $y_{\sigma\upharpoonright j}^{(n+1)}$ for each $\sigma\in T^{(n)}$
   by reverse induction on $j$.
   Start by choosing $y_{\sigma k}^{(n+1)}\bumpeq x_{\sigma k}$. 
   Fix $j<n+1$ and
   assume that the $y_{\sigma\upharpoonright j+1}^{(n+1)}$ are defined.
   Take a converging subsequence of 
   $\overrightarrow{\mathbf{y}}_{n,\sigma\upharpoonright j}^{(n+1)}$
   and call its limit $y_{\sigma\upharpoonright j}^{(n+1)}$, then,
   in $T^{(n+1)}$, delete all the immediate successors
   of $\sigma\upharpoonright j$ that are not in this subsequence. This shrinks $T^{(n+1)}$ 
   while (T1) and (T2) remain true. 
   The $n+1$-th stage finishes
   when $j$ reaches $0$,
   defining $T^{(n+1)}$ and $y_\eta^{(n+1)}$ for each $\eta\in T^{(n+1)}$.
   Then (XY2) and (XY3) hold by construction, we now show that (XY1) holds as well.
   \begin{claim}
      \label{claim:qcc2}
      Let $\sigma\in T^{(n+1)}$ with $h(\sigma)= n$. 
      Then $x_\sigma\bumpeq y_\sigma^{(n+1)}\bumpeq y_\sigma^{(n)}$.
   \end{claim}
   \begin{proof}
      By (XY2), 
      $\overrightarrow{\mathbf{x}}_{n+1,\sigma}$ converges to both $x_\sigma$ and $y_\sigma^{(n)}$
      and $\overrightarrow{\mathbf{y}}_{n+1,\sigma}^{(n+1)}$ converges to $y_\sigma^{(n)}$.
      Since the members of these sequences are $x_{\sigma k}\bumpeq y_{\sigma k}^{(n+1)}$,
      the conclusion holds by Lemma \ref{lemma:NHseq}.
   \end{proof}   
   \begin{claim}
      \label{claim:qcc3}
      Let $\sigma\in T^{(n+1)}$ with $h(\sigma)= n$ and let $j\le n$.
      Then $x_{\sigma\upharpoonright j}\bumpeq
      y_{\sigma\upharpoonright j}^{(n+1)}\bumpeq y_{\sigma\upharpoonright j}^{(k)}$
      for each $k$ with $j\le k\le n$.
   \end{claim}
   \begin{proof}[Proof of the claim] 
      For all $\sigma$ with $h(\sigma)= n$ at once
      and reverse induction on $j$.
      Claim \ref{claim:qcc2} provides the case $j=n$ (hence $k=n$).
      Assume that it holds for all $n\ge i>j$.
      By (XY3) $y_{\sigma\upharpoonright j}^{(n+1)}$ is the limit of 
      $\overrightarrow{\mathbf{y}}_{n+1,\sigma\upharpoonright j}^{(n+1)}$
      whose members are of the form $y_{(\sigma\upharpoonright j)k}^{(n+1)}$.
      Since $(\sigma\upharpoonright j)k = \eta\upharpoonright j+1$ for some $\eta\in T^{(n+1)}$
      with $h(\eta)=n$,
      by induction we have that 
      $y_{(\sigma\upharpoonright j)k}^{(n+1)}\bumpeq y_{(\sigma\upharpoonright j)k}^{(k)}$
      for each $j+1\le k\le n$.
      By (XY3) and Lemma \ref{lemma:NHseq}, 
      $y_{\sigma\upharpoonright j}^{(n+1)}\bumpeq y_{\sigma\upharpoonright j}^{(k)}$
      for each $j+1\le k\le n$. 
      The only case left is to check that 
      $y_{\sigma\upharpoonright j}^{(n+1)}\bumpeq y_{\sigma\upharpoonright j}^{(j)}$.
      But by (XY2) $y_{\sigma\upharpoonright j}^{(j)}$ is the limit of 
      $\overrightarrow{\mathbf{x}}_{j,\sigma\upharpoonright j}$
      which has $\overrightarrow{\mathbf{x}}_{n+1,\sigma\upharpoonright j}$
      as a subsequence. By induction,
      $x_{\eta\upharpoonright j+1}\bumpeq y_{\eta\upharpoonright j+1}^{(n+1)}$,
      for all $\eta$ with $h(\eta) = n$,
      hence the limits of $\overrightarrow{\mathbf{x}}_{n+1,\sigma\upharpoonright j}$
      and of $\overrightarrow{\mathbf{y}}_{n+1,\sigma\upharpoonright j}^{(n+1)}$,
      that is 
      $y_{\sigma\upharpoonright j}^{(j)}$ and
      $y_{\sigma\upharpoonright j}^{(n+1)}$, cannot be separated.
   \end{proof}
   After this painstaking indices management, we may finally finish the proof
   by noting that (XY1) implies that
   $Y = \{y_\varnothing^{(n)}\,:\,n\in\omega\}$ satisfies $z_1\bumpeq z_2$ for each $z_1,z_2\in Y$,
   $z_1\not= z_2$. Lemma \ref{lemma:disc} implies that $Y$ is infinite closed discrete,
   a contradiction.
\end{proof}

\begin{cor}
   \label{cor:noENHqcc}
   Let $M$ be a manifold and $Y\subset X$ be quasi-countably compact.
   Then $NH(Y)$ is closed and nowhere dense.
\end{cor}
\begin{proof}
   Suppose that $NH(Y)$ is dense in some Euclidean $U\subset M$. 
   By Lemma \ref{lemma:NHK}, $NH(Y)$ is closed and by density $NH(Y)\supset U$.
   Let $V\subset U$ be such that 
   $\wb{V}\cap U$ is compact. 
   Since $Y\cap \wb{V}$ is closed in $Y$, it is quasi-countably compact.
   Hence, $Z = (\wb{V}\cap U)\cup (Y\cap \wb{V})$ is a quasi-countably compact first countable
   locally Hausdorff space, but $\{x\in Z\,:\,NH(x)\not=\varnothing\}\supset V$, a contradiction.
\end{proof}

We may also phrase the same result as follows.
\begin{cor}
   Let $M$ be a manifold such that $\{x\in U\,:\,NH(x)\not=\varnothing\}$ is dense in some open subset $U$.
   Then $M$ is not quasi-countably compact. 
\end{cor}

These result directly imply the
following variant of Corollary \ref{cor:ClH} (see Section \ref{sec:cover} for the relevant definitions).
\begin{thm}
   \label{thm:Scountcpct}
   Let $S$ be a scan of some $\kappa$-Baire subset $A$ 
   of an Euclidean subset $U$ contained in an ENH-manifold $M$.
   Then $S$ cannot be covered by $<\kappa$-many subsets $K_\alpha$
   such that each $K_\alpha$ is quasi-countably compact.
   In particular, if $A=U$, this holds for $\kappa = \text{cov}(\mathcal{B})$.
\end{thm}
\begin{proof}
   Let $K_\alpha$ ($\alpha<\lambda<\kappa$) be 
   as in the statement with $\cup_{\alpha<\lambda} K_\alpha\supset S$.
   By Corollary \ref{cor:noENHqcc}, 
   $B = A-\cup_\alpha(NH(K_\alpha))\not=\emptyset$, any $x\in B$ can be separated from any
   point of $S$ and hence $NH(x) = \varnothing$.
\end{proof}

This theorem also holds with $\wb{U}$ instead of $S$ in the statement.
The remaining of the section contains two examples showing that (quasi-)countably compact subspaces of manifold
do not behave as well as (quasi-)compact ones.


\subsection{A contrarian countably compact subset}\label{subsec:contrarian}

In this subsection we show that Lemma \ref{lemma:compactinchart} does not 
hold if $K$ is only countably compact.

\begin{example}
   \label{ex:omega1notinchart}
  There is a NH-$1$-manifold $M$ with a copy $B$ of $\omega_1$ satisfying $NH(B)\cap B=\varnothing$
  such that no open Hausdorff $U$ contains $B$, actually $\mathcal{OH}(B,M)=\aleph_1$.
\end{example}
We shall use Fodor's lemma and basic properties of stationary subsets in the proof.
See e.g. \cite[pp. 126--127]{JustWeeseII} if in need of details.
(Note: there was a flaw in the space presented in a previous version of these notes.)
\begin{proof}[Details]
  We start with a remark. 
  Let $E_n = \{\frac{1}{m}\,:\,m>n\}$.
  Topologize the space $P = \R\times\{0\}\cup E_0\times\{1\}$
  by declaring the sets $(U-E_n)\times\{0\}\cup (U\cap E_n)\times\{1\}$ 
  and the sets $(U\cap(0,+\infty))\times\{0\}$
  open whenever $U\subset\R$ is open.
  Then $P$ is a NH-$1$-manifold, and the points in $E_0\times\{1\}$ converge to $\langle 0,0\rangle$.
  We inductively build a space using this simple idea repeatedly.
  Let $\LL_+$ be the long ray. 
  As in Example \ref{ex:nocloseddis}, we shall often consider $\omega_1$ as 
  a subspace of $\LL_+$ by identifying $\alpha$ with 
  $\langle\alpha,0\rangle$, and write directly $\alpha\in\LL_+$,
  and this holds as well for subsets of $\omega_1$ such as the $L_\alpha$ defined below.
  In the following description, the intervals $(\alpha,\beta)$ are understood to be intervals in 
  $\LL_+$. Write $\Lambda$ for the set of limit ordinals in $\omega_1$
  and $\Lambda_2$ for the limits of limits. 
  We fix a ladder system $\langle L_\alpha\,:\,\alpha\in\Lambda\rangle$ on $\omega_1$, 
  that is: $L_\alpha$ 
  is a strictly increasing $\omega$-sequence $\nu_{\alpha,n}$
  with supremum $\alpha$. 
  We may (and do) assume
  that $\nu_{\alpha,n}>0$ is successor for each $n>0$ and $\alpha\in\Lambda$, 
  and 
  that $\nu_{\alpha,n} = \mu_{\alpha,n} + 1$, where $\mu_{\alpha,n}\in\Lambda$,
  when $\alpha\in\Lambda_2$. It is also convenient to assume that $\mu_{\alpha,0}=0$
  for each $\alpha\in\Lambda_2$.
  We define $M$ to be the set 
  $$
     \LL_+\times\{0\}\cup\cup_{\alpha\in\Lambda}L_\alpha\times\{\alpha\}
  $$
  topologized as follows. Write $\pi$ for the projection on the first factor.
  {\em Caution}: the topology on $\LL_+\times\{0\}$
  is not that of $\LL_+$,
  as $L_\alpha\times\{0\}$ does not have $\langle\alpha,0\rangle$
  as a limit point (it is actually closed discrete). 
  The topology on $(\LL_+ -\Lambda)\times\{0\}$ is that of a disjoint union of intervals in $\LL_+$.
  A base of neighborhoods of $\langle\nu_{\alpha,n},\alpha\rangle$
  is given by the sets
  $(V-F)\times\{0\}\cup (V\cap F)\times\{\alpha\}$ for $V$ open containing $\nu_{\alpha,n}$
  and $F\subset L_\alpha$ finite.
  This defines the open sets containing all the points at height $>0$ 
  and those whose projection is not in $\Lambda$.
  We shall define $I_\alpha\ni\langle \alpha,0\rangle$ by induction on 
  $\alpha\in\Lambda$, with the property that $I_\alpha$ is homeomorphic to $\R$.
  Set $I_0 = (0,1)\times\{0\}$.
  If $\alpha = \beta + \omega\in\Lambda-\Lambda_2$, set 
  $$ I_\alpha = I_\beta\cup \left((\beta,\alpha+1)-L_\alpha\right)\times\{0\}
                \cup \left( L_\alpha\cap (\beta,\alpha)\right)\times\{\alpha\}.$$
  In words: take the interval $(\beta,\alpha+1)$, put it at height $0$, lift the
  points of $L_\alpha$ at height $\alpha$ and glue it to the right of $I_\beta$.
  Suppose now that $\alpha\in\Lambda_2$.
  Set 
  $$ J_{\alpha,n} = \left(I_{\mu_{\alpha,n+1}}\cap (\mu_{\alpha,n},\mu_{\alpha,n+1}]\right)\times\{0,1\}.$$
  Define now $I_\alpha$ by taking 
  $$ \left[\left(\bigcup_{n\in\omega} J_{\alpha,n}\right) -\pi^{-1}(\pi(L_\alpha)) \right]\cup 
     \left( L_\alpha\times\{\alpha\} \right)
     \cup \left( [\alpha,\alpha+1)\times\{0\}\right).
  $$ 
  In words: take the union of the $J_{\alpha,n}$, lift all the points with first coordinate in $L_\alpha$
  at height $\alpha$, and add the half open interval $[\alpha,\alpha+1)\times\{0\}$.
  Now, declare open the sets $I_\alpha\cap (x,y)\times\{0,1\}$,
  where $x<\alpha<y$, $x,y\in\LL_+$. 
  \begin{claim}
     $\pi\upharpoonright I_\alpha:I_\alpha\to (0,\alpha+1)$ is an homeomorphism and $M$ is a NH-$1$-manifold.
  \end{claim}
  \begin{proof}[Proof of the claim]
     If $\alpha = \beta+\omega$, then $h$ restricted to $H_\alpha = I_\alpha\cap(\beta,\alpha+1)\times\{0,1\}$ 
     is clearly an homeomorphism.
     Then the result follows since
     $I_\alpha = H_\alpha\cup I_\beta$ and $H_\alpha\cap I_\beta = (\beta,\beta+1)\times\{0\}$.\\
     Let now $\alpha\in\Lambda_2$.
     By induction, $h\upharpoonright J_{\alpha,n}$ is a homeomorphism onto $(\mu_{\alpha,n},\mu_{\alpha,n+1}]$,
     and hence $h$ restricted to $\bigcup_{n\in\omega} J_{\alpha,n}$ 
     is an homeomorphism onto $(0,\alpha)$.
     But when we lift 
     the points of $L_\alpha$, only finitely many of them lie in each $J_{\alpha,n}$, and
     they are not the right endpoint since $\mu_{\alpha,n}\not=\nu_{\alpha,m}$. Hence we end up with the same
     construction as in the space $P$ above, and the result follows.
  \end{proof}
  \begin{claim}
     $B=\Lambda\times\{0\}$ is a copy of $\omega_1$ and $NH(B) = \varnothing$.
  \end{claim}
  \begin{proof}[Proof of the claim]
     If $x,y\in\LL_+$ with $x<y$,
     then $B\cap \pi^{-1}((x,y)) = \bigl(\Lambda\cap(x,y)\bigr)\times\{0\}$.
     The induced topology on $\Lambda\times\{0\}$ is thus the same as that 
     of $\Lambda$ as a subset of  $\omega_1$, which is a copy of $\omega_1$.
     Points that cannot be separated in $M$ are those with same first coordinate,
     but $\pi^{-1}(\pi(B))=B$, hence $NH(B) = \varnothing$.
  \end{proof}
  \begin{claim}
     Any open set containing $\Lambda\times\{0\}$ is non-Hausdorff.
  \end{claim}
  \begin{proof}[Proof of the claim]
     Let $U\supset\Lambda\times\{0\}$ be open.
     For each $\alpha\in\Lambda$ there is some $\beta<\alpha$ such that 
     $I_\alpha\cap (\beta,\alpha]\times\{0,1\}\subset U$.
     By Fodor's lemma, there is a stationary $S_0$ and some $\beta$ such that 
     $I_\alpha\cap (\beta,\alpha]\times\{0,1\}\subset U$ for each $\alpha\in S_0$.
     Since a countable union of non-stationary sets is non-stationary,
     there is some $n\in\omega$ and a stationary $S_1\subset S_0$
     such that $\nu_{\alpha,n}>\beta$ for each $\alpha\in S_1$.
     Again by Fodor, there is a stationary $S_2\subset S_1$ and $\gamma\ge\beta$ such that
     $\nu_{\alpha,n}=\gamma$ for each $\alpha\in S_2$.
     Let $\alpha$ be a limit point of $S_2$ (which exists by stationarity)
     and let $\delta\in S_2$ be such that $\gamma < \delta < \alpha$.
     Then by definition of $I_\delta$ and $I_\alpha$,
     both $\langle \gamma,\delta\rangle$ and $\langle \gamma,\alpha\rangle$
     are in $U$.
     But these points cannot be separated in $M$ since their neighborhoods
     intersect at the $0$-th floor.
  \end{proof}
  This shows that $M$ and $B$ have the required properties.
  \end{proof} 

Notice in passing that the $1$-manifold just described is not (homeomorphic to) a towel-rack manifold,
since the neighborhoods of some points intersect more that $2$ ``floors''.
We end this subsection with an easy lemma about copies of $\omega_1$ is manifolds which
did not find any use in these notes but which we record for future reference.

\begin{lemma}
   Let $e:\omega_1\to M$ be an homeomorphism on its image, where $M$ is a manifold.
   Then there is $\alpha\in\omega_1$ such that 
   for each $\beta>\alpha$ there is an Euclidean $U_\beta\subset M$ such
   that $e([\alpha,\beta])\subset U_\beta$ (where the interval is taken in $\omega_1$).
\end{lemma}
In particular, $e(\beta)$ and $e(\gamma)$ can be separated in $M$ for each $\gamma,\beta>\alpha$.
\begin{proof}
   For each limit $\beta$ choose an Euclidean $V_\beta\ni e(\beta)$.
   Since $e$ is a homeomorphism on its image, there is some $\alpha(\beta)<\beta$ such that
   $V_\beta\supset e([\alpha(\beta),\beta])$.
   By Fodor's lemma, there is some $\alpha$ and a stationary $S\subset\omega_1$ such that
   $\alpha(\beta) = \beta$ for each $\beta\in S$. Set $U_\beta = V_\gamma$ for the smallest
   $\gamma\in S$ such that $\beta\le\gamma$.
\end{proof}


\subsection{No copy of $\omega_1$ under $\clubsuit$, following Nyikos}\label{subsec:noomega_1}

Whether a countably compact (hence, Hausdorff) non-compact manifolds contains a copy of $\omega_1$ is 
independant of {\bf ZFC}: the answer is yes under {\bf PFA} (Balogh \cite{Balogh:1989})
and no under $\diamondsuit$ (Rudin, see e.g. \cite[Example 3.14]{Nyikos:1984}).
Actually, the late Peter Nyikos 
even gave an example of a countably compact non-compact $2$-manifold without a copy of $\omega_1$
under the weaker axiom $\clubsuit_C$ (a weakening of $\clubsuit$ compatible with {\bf MA + $\neg$CH}) 
in an unpublished draft 
(see \cite[Section 5]{Nyikos:Antidiamond}).
(Notice that we do not know whether this construction really works as there are gaps in the presentation
that we are not able to overcome.)
In another unpublished draft, he gives the ideas 
(see \cite[Theorem 8.1]{Nyikoscoherent})
for building the so-called {\em sprat} that we decribe
below under $\clubsuit$. This construction is simple and immediately yields
a quasi-countably compact non-quasi-compact NH-$1$-manifold containing no copy of $\omega_1$.
In this short subsection, we fill the details in Nyikos' sprat construction and 
describe the said manifold.
Of course, it is a bit of an overkill to use $\clubsuit$ to obtain a non-Hausdorff example while Nyikos' example 
in his draft is Hausdorff and needs only $\clubsuit_C$,
but the fact that our example is $1$-dimensional 
(and Nyikos'one unclear status) made it worth describing in our opinion.

{\em Sprat} is a nickname Nyikos uses for what is
usually called a {\em $2$-to-$1$ closed preimage of $\omega_1$}. That is,
a Hausdorff space $X$ equipped with a closed map $\pi:X\to\omega_1$ such that
$|\pi^{-1}(\{\alpha\})|=2$ for each $\alpha\in\omega_1$.
We may assume that the underlying set of $X$ is $\omega_1\times\{0,1\}$.
Since $\pi$ is closed, the space is countably compact. It implies that
for any limit $\alpha$, the union of any two neighborhoods of $\langle \alpha,0\rangle$ and $\langle \alpha,1\rangle$
must cover all of $\pi^{-1}((\gamma,\alpha])$ for some $\gamma<\alpha$
(the intervals are always taken in $\omega_1$
in this description).
We say that $X$ is {\em symmetrical} iff for each $\alpha$ there is a clopen $B^0_\alpha\ni \langle \alpha,0\rangle$
such that $|B^0_\alpha\cap\pi^{-1}(\{\beta\})|=1$ for each $\beta\le\alpha$. 
This implies that $B^1_\alpha = \pi^{-1}([0,\alpha])-B^0_\alpha$ is a clopen set containing 
$\langle \alpha,1\rangle$ and such that 
$|B^1_\alpha\cap\pi^{-1}(\{\beta\})|=1$ for each $\beta\le\alpha$ as well.
A base of neighborhoods of $\langle \alpha,i\rangle$ is then given by 
$B^i_\alpha \cap \pi^{-1}\left((\gamma,\alpha]\right)$
for $\gamma<\alpha$. In particular, $\langle\alpha+1,i\rangle$ is isolated in $X$.  
\\
Let $\mathscr{L} = \{L_\alpha\,:\,\alpha<\omega_1\}$ be a ladder system on $\omega_1$.
As above, we let $\{\nu_{\alpha,n}\,:\,n\in\omega\}$ be an increasing enumeration of $L_\alpha$.
We assume that $\nu_{\alpha,0}=0$ for each $\alpha$.
We shall now define a symmetrical sprat $X^\mathscr{L}$ 
by specifying the sets $B^i_\alpha$ for each $\alpha<\omega_1$.
Notice that this defines a sprat topology iff whenever $\beta<\alpha$ are limit ordinals
with $\langle\beta,j\rangle\in B^i_\alpha$, then there is $\gamma<\beta$ with
\begin{equation}
   \label{eq:symsprat}
   \tag{$\bullet$}
   B^i_\alpha\cap\pi^{-1}((\gamma,\beta]) = B^j_\beta\cap\pi^{-1}((\gamma,\beta]).
\end{equation}
If $\alpha<\omega$ we let $B^i_\alpha$ be $[0,\alpha]\times\{i\}$.
We now proceed by induction.
If $\alpha=\beta+1$, we let $B^i_\alpha = \{\langle \alpha,i\rangle\}\cup B^i_{\beta}$.
For $\alpha$ limit, we 
write $R^\alpha_n = \pi^{-1}\left((\nu_{\alpha,n},\nu_{\alpha,n+1}]\right)$ and set
\begin{equation}
\label{eq:2-1}
\tag{$\bullet\bullet$}
B^i_\alpha = \cup_{n\in\omega} \left(
                (B^i_{\nu_{\alpha,2n+1}}\cap R_{2n}^\alpha)
                \cup
                (B^{1-i}_{\nu_{\alpha,2n+2}}\cap R_{2n+1}^\alpha)\right)
                .
\end{equation}
\begin{claim}
   For any $i,j\in\{0,1\}$ and 
   $\beta<\alpha$, both limit, with $\langle\beta,j\rangle\in B^i_\alpha$, 
   there is $\gamma<\beta$ such that (\ref{eq:symsprat}) holds.
\end{claim}
\begin{proof}[Proof of the claim]
   We prove it by induction on $\alpha$. 
   The result is clear if $\alpha = \omega$. Let  
   $n$ be maximal such that $\nu_{\alpha,n}<\beta$.
   Then (\ref{eq:2-1}) yields that 
   $$ B^j_{\nu_{\alpha,n+1}}\cap\pi^{-1}\left((\nu_{\alpha,n},\beta]\right)
      = B^i_{\alpha}\cap\pi^{-1}\left((\nu_{\alpha,n},\beta]\right).
   $$
   The result follows by induction since $\nu_{\alpha,n+1}<\alpha$.
\end{proof}

\begin{claim}
   \label{claim:2-1}
   For $i\in\{0,1\}$ and $\alpha\in\Lambda$, the closure of 
   $L_\alpha\times\{i\}$ contains $\pi^{-1}(\{\alpha\})$.
\end{claim}
\begin{proof}[Proof of the claim]
   Indeed, by definition $B^j_\beta$ contains $\langle\beta,j\rangle$, hence 
   by (\ref{eq:2-1}) any neighborhood of $\langle\alpha,i\rangle$ contains both
   $\langle\nu_{\alpha,2n+1},i\rangle$ and $\langle\nu_{\alpha,2n+2},1-i\rangle$
   for each big enough $n$.
\end{proof}

We let {\em club} be a shorthand for closed and unbounded.
\begin{lemma}
    \label{lemma:notomega_1inX}
   Let $X^\mathscr{L}$ be the symmetrical sprat defined above, where
   $\mathscr{L}=\{L_\alpha\,:\,\alpha\in\Lambda\}$ is a $\clubsuit$-sequence.
   Then any club subset $C$ of 
   $X^\mathscr{L}$ (i.e. $\pi(C)$ is club)
   contains $\pi^{-1}(\{\alpha\})$ for a club subset of $\alpha$.
\end{lemma}
\begin{proof}
   The fact that the set of such $\alpha$ is closed is immediate, let us see that is is unbounded. 
   Let $C_i=C\cap (\omega_1\times\{i\})$, then one of $C_0,C_1$ (say $C_0$)
   must be uncountable (stationary, actually).
   Then there are stationary many $\alpha$ such that $L_\alpha\subset C_0$, and $C\supset\pi^{-1}(\{\alpha\})$
   for these $\alpha$
   by Claim \ref{claim:2-1}. 
\end{proof}

\begin{cor}
    \label{cor:notomega_1inX}
   If $\mathscr{L}=\{L_\alpha\,:\,\alpha\in\Lambda\}$ is a $\clubsuit$-sequence, 
   $X^\mathscr{L}$ does not contain any copy of $\omega_1$.
\end{cor}
\begin{proof}
   It is easy to see that if $Y$ is a sprat, any $1$-to-$1$ map $\omega_1\to Y$
   is strictly increasing on a club subset. Hence any copy
   of $\omega_1$ contains one which intersects 
   each $\pi^{-1}(\alpha)$ in at most one point, contradicting the previous lemma.
\end{proof}

We may now define our NH-$1$-manifold.
   Take a symmetrical sprat $X^\mathscr{L}= \omega_1\times\{0,1\}$ with $B^j_{\alpha,n}$ defined as above.
   Let $M^\mathscr{L} = \bigl((\LL_{\ge 0}-\omega_1)\times\{-1\}\bigr)\cup X^\mathscr{L}$.
   (As in Example \ref{ex:omega1notinchart}, we treat $\omega_1$ as a subset of $\LL_{\ge 0}$.)
   The topology is as follows. On the $(-1)$-th floor, we take the usual topology as a discrete union of open intervals.
   If $u<v$ are points of $\LL_{\ge 0}$, write $(u,v)_\bullet$ for the subset 
   $\{\langle x,-1\rangle\,:\,x\in\LL_{\ge 0}-\omega_1,\,u<x<v\}$.
   The neighborhoods of $\langle\alpha,i\rangle\in X^\mathscr{L}$
   are $(u,v)_\bullet\cup \left(B^i_\alpha\cap \pi^{-1}((\gamma,\alpha])\right)$ where 
   $\gamma<u<\alpha<v<\alpha+1$, $\gamma\in\omega_1$ and $u,v\in\LL_{\ge 0}-\omega_1$.
   In words: we insert an open interval between the points in $\pi^{-1}(\{\alpha\})$ and
   $\pi^{-1}(\{\alpha+1\})$, as we do in $\LL_{\ge 0}$, but so
   that $\langle\alpha,0\rangle\bumpeq\langle\alpha,1\rangle$ and 
   $\langle\alpha+1,0\rangle\bumpeq\langle\alpha+1,1\rangle$.
   By construction, the subspaces $(0,\alpha)_\bullet\cup B^i_\alpha$ are copies of
   the compact interval $[0,\alpha]\subset\LL_{\ge 0}$ and 
   thus $M^\mathscr{L}$ is quasi-countably compact.

\begin{example}
   If $\clubsuit$ holds, there is a quasi-countably compact non-compact NH-$1$-manifold $M$
   which does not contain any copy of $\omega_1$. In fact, any quasi-countably compact
   non-compact subspace is non-Hausdorff, and $\mathcal{OH}(M^\mathscr{L})=\omega_1$.
\end{example}
\begin{proof}[Details]
   Build $X^\mathscr{L}$ and $M=M^\mathscr{L}$ as above.
   Any quasi-countably compact non-quasi-compact subset of $M^\mathscr{L}$ 
   intersects $X^\mathscr{L}$ on a club subset,
   and is thus non-Hausdorff by Lemma \ref{lemma:notomega_1inX}. 
   The claim about $\mathcal{OH}(M^\mathscr{L})$ follows easily.
   \end{proof}

An observant reader might complain that $M^\mathscr{L}$ 
is not exactly a manifold because the points $\langle 0,0\rangle$
and $\langle 0,1\rangle$ have neighborhoods homeomorphic to $[0,1)$ and not to $(0,1)$.
In order not to destroy quasi-countable compactness, we might for instance take two copies of $M^\mathscr{L}$ 
and identify these
points, or take a copy of the interval $[0,1]$ and glue $0$ to $\langle 0,0\rangle$
and $1$ to $\langle 0,1\rangle$.

\begin{ques}
  Is there in {\bf ZFC} a quasi-countably compact non-quasi-compact NH-manifold which does not contain $\omega_1$~?
\end{ques} 

As said above, if the manifolds are Hausdorff the answer is ``no''  under the {\bf PFA}. 
We do not know whether the proof(s) can be adapted to non-Hausdorff manifolds.


\section{Examples of NH-$1$-manifolds.}
\label{sec:1dim}

We already saw that 
NH-manifolds can be quite complicated even in dimension $1$. This section's goal is to 
exhibit other specimens that we found interesting. 

\subsection{Homogeneous and inhomogenous manifolds with ``shortcuts''}

\vskip .3cm
\noindent
While $NH(x)$ is closed nowhere dense,
$NH^\infty(x)$ can be almost anything. 
Recall that $[E]^2$ denotes the $2$-elements subsets of $E$, which we identify with
the pairs $\langle x,y\rangle$ with $x<y$ if $E$ is totally ordered.

\begin{example}\label{ex:E}
   Let $E\subset \R$. 
   There is a NH-$1$-manifold $W_E = \R\times\{0\}\sqcup [E]^2\times\{1\}$
   with $\mathcal{E}(W_E)=\max(|E|,\omega)$
   such that, for each $x,y\in E\times\{0\}$, we have $x\Bumpeq y$ in $W_E$.
\end{example}
More precisely, the $\Bumpeq$-classes of $W_E$ are $E\times\{0\} \cup [E]^2\times\{1\}$
   and singletons.
\begin{proof}[Details]
   We topologize $W_E=\R\times\{0\}\sqcup [E]^2\times\{1\}$ as 
   follows. 
   Write $c_{x,y}$ for the point $\{ x,y\}\times\{1\}\in [E]^2\times\{1\}$.
   Neighborhoods of $x\in\R\times\{0\}$ are the usual neighborhoods, while a neighborhood of 
   $c_{x,y}$ is $\left((x-\epsilon,x)\cup (y,y+\epsilon)\right)\times\{0\}\cup \{c_{x,y}\}$ for $\epsilon >0$.
   We see $c_{x,y}$ as a shortcut between $x$ and $y$. 
   By definition, $x\bumpeq c_{x,y}\bumpeq y$, hence $x\Bumpeq y$, whenever $x,y\in E$. 
   Any other member of $W_E$ can be separated from all the other points.
\end{proof}

\begin{example}
    Any continuous map $f:W_\Q\to H$ is constant
    when $H$ is an Hausdorff space.
\end{example}

\begin{proof}[Details]
    By Lemma \ref{lemma:denseconstant}, since $\Q\times\{0\}$ is dense in $W_\Q$.
\end{proof}

Pushing the same idea further yields the following homogeneous example.
The symbol \fatS is to be understood as a ``$\mathbb{S}^1$ with specks of dust''.

\begin{example}
   \label{ex:FatS1}
   A HNH-$1$-manifold \fatS such that $x\Bumpeq y$ for each $x,y\in$\fatS
   and any continuous $f:\text{\fatS}\to H$ is constant when $H$ is Hausdorff. 
   Moreover, each point of \fatS has a $2$-rim-simple neighborhood.
\end{example}
\begin{proof}[Details]
   As a set, \fatS is the union of the circle $\mathbb{S}^1$ and of four points 
   $c_{x^+,y^+},c_{x^+,y^-},c_{x^-,y^+},c_{x^-,y^-}$ for each $\{x,y\}\in [\mathbb{S}^1]^2$.
   Since order does not matter, $c_{x^\pm,y^\pm} = c_{y^\pm,x^\pm}$.
   We see $\mathbb{S}^1$ as the real numbers modulo $1$, and expressions such as $x\pm\epsilon$ have
   to be understood modulo $1$. 
   The neighborhoods of points in $\mathbb{S}^1$ are the usual neighborhoods in $\mathbb{S}^1$.
   Neighborhoods of the other points are defined as:
   \begin{align*}
      c_{x^+,y^+}:& \{c_{x^+,y^+}\}\cup (x,x+\epsilon)\cup (y,y+\epsilon),\\
      c_{x^+,y^-}:& \{c_{x^+,y^-}\}\cup (x,x+\epsilon)\cup (y-\epsilon,y),\\
      c_{x^-,y^+}:& \{c_{x^-,y^+}\}\cup (x-\epsilon,x)\cup (y,y+\epsilon),\\
      c_{x^-,y^-}:& \{c_{x^-,y^-}\}\cup (x-\epsilon,x)\cup (y-\epsilon,y),
   \end{align*}
   where $\epsilon > 0$ and each interval is a subset of $\mathbb{S}^1$. See Figure \ref{fig:1}.
   \begin{figure}[h]
     \begin{center}
      \epsfig{figure=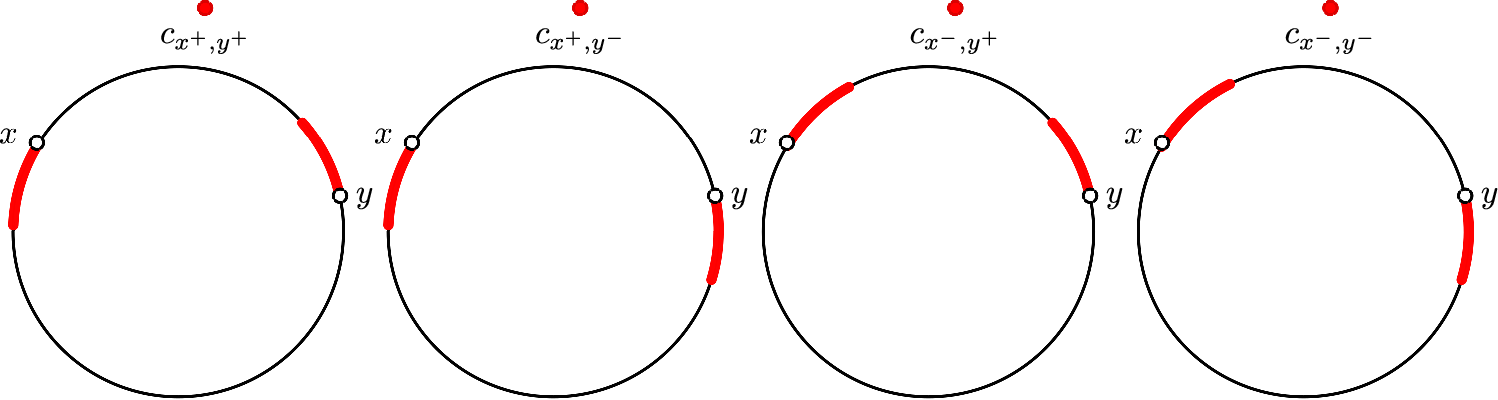, width=.9\textwidth}
      \caption{Neighborhoods of $c_{x^+,y^+},c_{x^+,y^-},c_{x^-,y^+},c_{x^-,y^-}$ in Example \ref{ex:FatS1}.}
      \label{fig:1}
      \end{center}
   \end{figure}
   It is clear that whenever $x,y\in\mathbb{S}^1$, then $x\bumpeq c_{x^\pm,y^\pm}\bumpeq y$, hence
   $a\Bumpeq b$ for each $a,b\in$\fatS. 
   The claim about $f$ being constant then follows from Lemma \ref{lemma:denseconstant}.
   We now show that \fatS is homogeneous. First, any rotation $x\mapsto x+\theta$, 
   $c_{x^\pm,y^\pm}\mapsto c_{(x+\theta)^\pm,(y+\theta)^\pm}$ is a homeomorphism.
   We will show that for each $x\not=0$ in $\mathbb{S}^1$, the following hold.
   \begin{itemize}
   \itemsep -.1cm
       \item[(1)] There is a homeomorphism interchanging $0$ with $c_{x^-,0^-}$ and $x$ with $c_{x^+,0^+}$.      
       \item[(2)] There is a homeomorphim
             interchanging $c_{x^+,0^-}$ with $c_{x^+,a^+}$ for some $a\in\mathbb{S}^1$.
       \item[(3)]There is a homeomorphim interchanging $c_{x^-,0^+}$ with $c_{a^-,b^-}$ for some $a,b\in\mathbb{S}^1$.    
   \end{itemize}
   This implies that \fatS is homogeneous: (1) and rotations shows that any point of $\mathbb{S}^1$
   can be sent to $c_{x^+,y^+}$ or $c_{x^-,y^-}$ for any $x,y\in\mathbb{S}^1$;
   (2)--(3) show how to reach $c_{x^+,y^-}$ and $c_{x^-,y^+}$.\\
   So, let us now show (1); that is,
   given $0\not=x\in\mathbb{S}^1$, describe a homeomorphism $h$
   which interchanges $c_{x^+,0^+}$ with $x$ and $c_{x^-,0^-}$ with $0$. 
   The images of these points being decided, let us see what happens to the other points.
   Traveling counterclockwise from $0$ in $\mathbb{S}^1$, consider the open interval $(0,x)$. 
   Then $h$ flips this interval by symmetry around its middle point $\frac{x}{2}$,
   and is the identity on the rest of $\mathbb{S}^1$ (except of course $x$ and $0$),
   as in Figure \ref{fig:2}.
   \begin{figure}[!t]
     \begin{center}
      \epsfig{figure=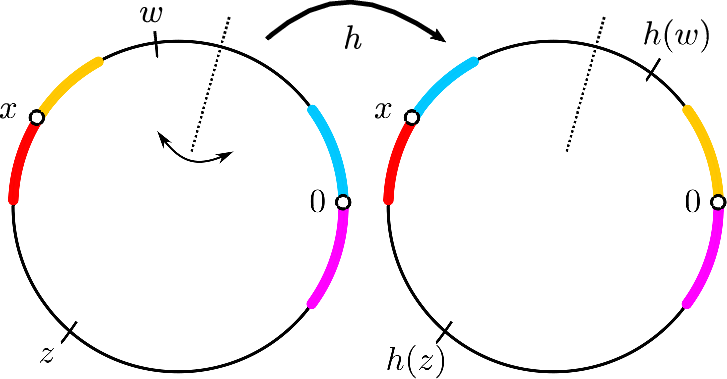, width=.45\textwidth}
      \caption{The homeomorphism $h$}
      \label{fig:2}
      \end{center}
   \end{figure}
   This turns a neighborhood of $x$ (in red and yellow) into one of $c_{x^+,0^+}$  
   and a neighborhood of $0$ (in cyan and magenta) into one of $c_{x^-,0^-}$, and vice-versa.
   We set $h(c_{x^+,0^-}) = c_{x^+,0^-}$ and $h(c_{x^-,0^+}) = c_{x^-,0^+}$,
   because neighborhoods of these points are preserved (that of $c_{x^-,0^+}$ is flipped, however).
   Now, given $z,w\in\mathbb{S}^1$, the images of $c_{z^\pm,w^\pm}$ can be infered from
   how the intersection of their neighborhoods with $\mathbb{S}^1$ changes under $h$.
   If both $z,w$ are on the flipped interval, then we move and 
   change the signs: $c_{z^+,w^-}$ is sent to
   $c_{h(z)^-,h(w)^+}$, for instance. 
   If both are on the unflipped interval, then $h$ is the identity on them.
   If $z$ is on the flipped interval and $w$ is not, we only move and change the sign of $z$:
   $h(c_{z^+,w^-}) = h(c_{h(z)^-,w^-})$, for instance.
   The cases of $c_{z^\pm,x^\pm}$ and $c_{z^\pm,0^\pm}$ are similar: 
   by looking at their neighborhoods, one sees that 
   if $z$ is in the unflipped interval, then we must set:
   \begin{align*}
      h(c_{z^+,x^-}) & = c_{z^+,0^+} & h(c_{z^-,x^-}) & = c_{z^-,0^+} \\
      h(c_{z^+,x^+}) & = c_{z^+,x^+} & h(c_{z^-,x^+}) & = c_{z^-,x^+} \\
      h(c_{z^+,0^+}) & = c_{z^+,x^-} & h(c_{z^-,0^+}) & = c_{z^-,x^-} \\
      h(c_{z^+,0^-}) & = c_{z^+,0^-} & h(c_{z^-,0^-}) & = c_{z^-,0^-}.
   \end{align*}
   Finally, if $w$ is in the flipped interval, 
   \begin{align*}
      h(c_{w^+,x^-}) & = c_{h(w)^-,0^+} & h(c_{w^-,x^-}) & = c_{h(w)^+,0^+} \\
      h(c_{w^+,x^+}) & = c_{h(w)^-,x^+} & h(c_{w^-,x^+}) & = c_{h(w)^+,x^+} \\
      h(c_{w^+,0^+}) & = c_{h(w)^-,x^-} & h(c_{w^-,0^+}) & = c_{h(w)^+,x^-} \\
      h(c_{w^+,0^-}) & = c_{h(w)^-,0^-} & h(c_{w^-,0^-}) & = c_{h(w)^+,0^-}. 
   \end{align*}
   This defines $h$ and gives (1).\\ 
   Notice that given $z$ in the unflipped interval and $w$ in the flipped interval, $h$ interchanges 
   $c_{z^+,x^-}$ with $c_{z^+,0^+}$ and $c_{w^-,x^+}$ with $c_{h(w)^+,x^+}$.
   By pre- or post-composing with a rotation, this actually yields (2) and (3): 
   change the roles of $x$ and $z$, and
   given $z$, choose $x$ either ``before'' or ``after'' $z$ defining $h$ as above.
\end{proof}

\begin{rems}\ \\
   (1) An open Hausdorff subset of \fatS can contain only countably many $c_{x^\pm,y^\pm}$, hence
       $\mathcal{OH}(\text{\fatS}) = \mathfrak{c}$. \\
   (2) $NH^2(x) = \text{\fatS}$ for each $x\in\text{\fatS}$.
   Indeed, for any $x,y,z\in\mathbb{S}^1$,
   we have:
   \begin{align*} c_{x^+,y^+}&\bumpeq x \bumpeq c_{x^\pm,z^\pm} &
                  c_{x^+,y^\pm}&\bumpeq c_{x^+,z^\pm} &
                  c_{x^-,y^\pm}&\bumpeq c_{x^-,z^\pm} \\
                  c_{x^\pm,y^+}&\bumpeq c_{z^\pm,y^+} &
                  c_{x^\pm,y^-}&\bumpeq c_{z^\pm,y^-}
   \end{align*}
   Hence, for each $u,v\in\text{\fatS}$, there is a $\bumpeq$-sequence of length at most $3$ 
   between $u$ and $v$.\\
   (3) Let 
        $$ A = \bigcup_{x\in\mathbb{S}^1-\{0\}}\{c_{0^+,x^+},c_{0^+,x^-},c_{0^-,x^+},c_{0^-,x^-}\}.$$ 
       Then $A$ is discrete and hence
       Hausdorff in the subspace topology but $NH(A)=$\fatS, hence in particular 
       $NH(A)$ is closed but not nowhere dense. This shows that Lemma \ref{lemma:sep} (d)
       does not hold if $A$ is only assumed to be Hausdorff in the subspace topology.\\
   \label{rem:FatS1}
\end{rems}

\begin{ques}
   \label{ques:BumpeqOHc}
   Is there a HNH-manifold $M$ with $\mathcal{OH}(M)<\mathfrak{c}$ such that $x\Bumpeq y$ for each $x,y\in M$~?
\end{ques}
Of course, this cannot happen if $\text{cov}(\mathcal{B}) = \mathfrak{c}$ by Theorem \ref{thm:ClH},
hence examples (if they exist at all) can only be consistent.


\subsection{NH-$1$-Manifolds with ``complicated'' $NH(x)$.}\label{subsection:NH1nd}

Let us now look for NH-$1$-manifolds with ``complicated'' $NH(x)$ (for some $x$).
First, let us see that $NH(x)$ can be homeomorphic to spaces which are locally 
similar to the converging sequence space $\omega+1$.
There is a general construction which is actually rather simple to obtain this kind of spaces
(see Example \ref{ex:NHclnd}).
But since we do like to use images, let us first start with a 
more picturesque (under our opinion) specimen.

\begin{example}\label{ex:NHordinal}
   A NH-$1$-manifold $M$ with $\mathcal{E}(M) = 2$ 
   such that for some $x\in M$, $NH(x)$ is homeomorphic to a countable ordinal $\alpha$
   (with the order topology). 
\end{example}

\begin{figure}[h]
     \begin{center}
        \epsfig{figure=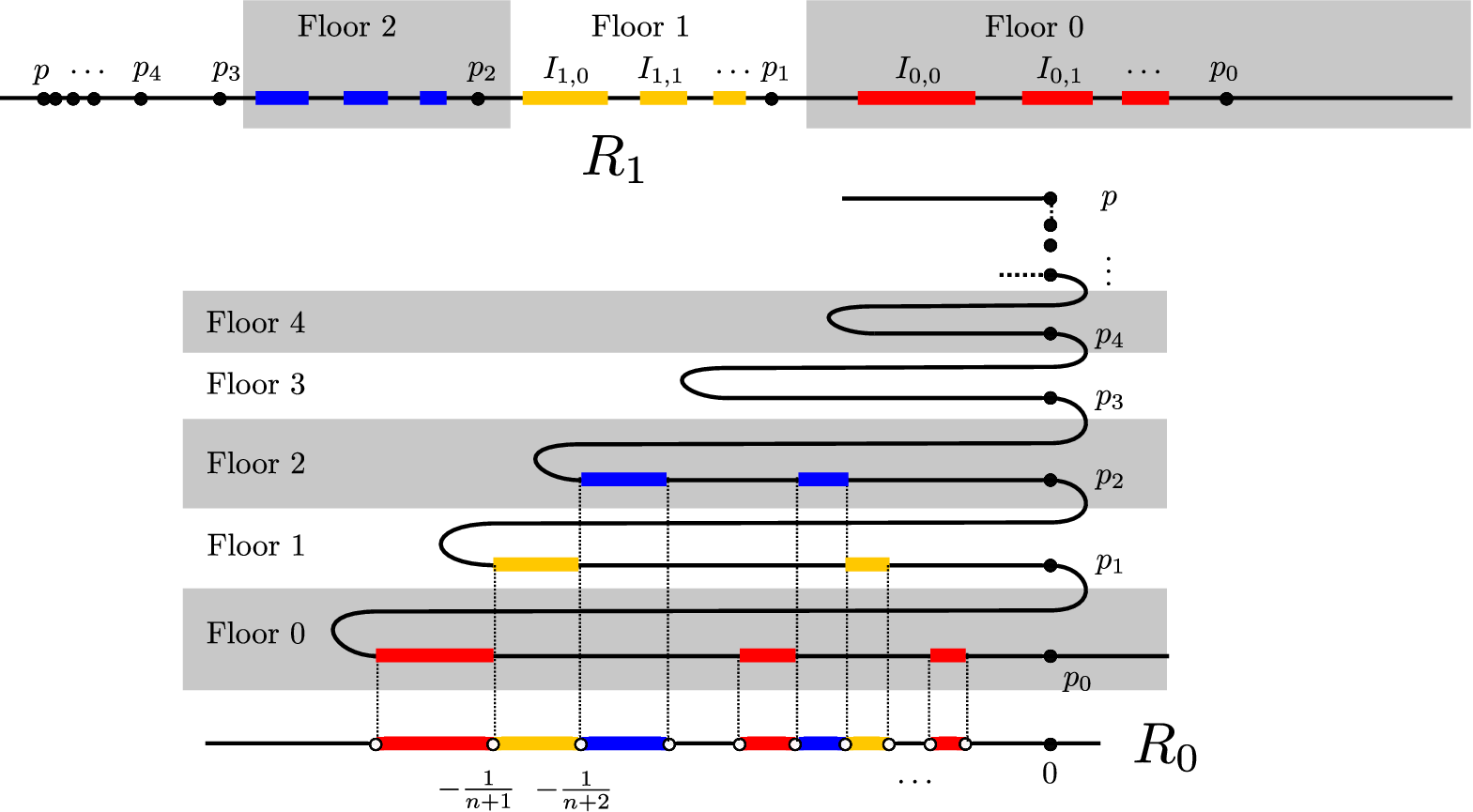, width=.9\textwidth}
         \caption{A NH-$1$-manifold such that $NH(0)$ is a converging sequence.}
         \label{fig:3}
      \end{center}
\end{figure}

\begin{proof}[Details]
   We first explain how to obtain $\omega+1$.
   Let $R_0$, $R_1$ be two copies of $\R$.
   In $R_1$, take a decreasing sequence $p_n$ ($n\in\omega$) converging to some $p$,
   and between $p_{n+1}$ and $p_n$, take an increasing sequence of disjoint
   open intervals $I_{n,m}$ ($m\in\omega$)
   which converge to $p_n$, as in the top part of Figure \ref{fig:3}.
   Let $\sigma:\omega\to\omega\times\omega$  be a bijection.
   Let $\varphi_n: I_{\sigma(n)}\to (-\frac{1}{n+1} , -\frac{1}{n+2})$ be any homeomorphism. 
   We let $M$ be $R_0\sqcup R_1$ where we identify $I_{\sigma(n)}\subset R_1$ with 
   $(-\frac{1}{n+1} , -\frac{1}{n+2})\subset R_0$ through $\varphi_n$. 
   A way to look at these identifications in pictured in Figure \ref{fig:3}.
   (The word ``floor'' does not have a real signification, but we find it
   easier to think about this construction this way.)
   We write $0$ for the origin in $R_0$.
   Then $NH(0) = \{p_n\,:\,n\in\omega\}\cup\{p\}\subset R_1$.\\
   To obtain $NH(0)$ homeomorphic to a countable ordinal $\alpha$ is almost the same.
   Start by embedding $\alpha$ in $R_1$ (sending $\beta$ to some $p_\beta$).
   We may choose the embedding to be decreasing.
   Take sequences of open intervals $I_{\alpha,n}$ between $p_{\alpha+1}$ and $p_\alpha$ converging to the latter.
   Then, choose a bijection $\sigma:\omega\to\alpha\times\omega$ and define $M$ as above.
\end{proof}

This construction can be generalized as follows.

\begin{example}
   \label{ex:NHclnd}
   For any closed nowhere dense set $E\subset\R$, there is
   a NH-$1$-manifold $M$ with $\mathcal{E}(M)= 2$ and a point $x$ for which $NH(x)$ is a copy of $E$.
\end{example}
\begin{proof}[Details]
   The complement of $E$ in $\R$ is dense and the union of countably many disjoint open intervals $(a_n,b_n)$ 
   ($n\in\omega$). Let $I_{n,m,0}$, $I_{n,m,1}$ ($m\in\omega$) be the intervals 
   $\left( a_n + \frac{b_n-a_n}{2^{m+2}} , a_n + \frac{b_n-a_n}{2^{m+1}}\right)$
   and $\left( b_n - \frac{b_n-a_n}{2^{m+1}} , b_n - \frac{b_n-a_n}{2^{m+2}}\right)$.
   As before, let $\sigma:\omega\to\omega\times\omega\times\{0,1\}$  be a bijection and
   $\varphi_n: I_{\sigma(n)}\to (-\frac{1}{n+1} , -\frac{1}{n+2})$ be an homeomorphism;
   let $M$ be the union of two copies $R_0,R_1$ of $\R$ with identifications through the $\varphi_n$'s.
   Since any neighborhood of a point $p$ of $E$ contains either
   infinitely many $I_{n,m,0}$ or infinitely many $I_{n,m,1}$ (or both),
   its copy in $R_1$ cannot be separated from the copy of $0$ in $R_0$,
   while points in $R_1-E$ can.
\end{proof}
Notice that whether some (or each) $\varphi_n$ preserves the local orientation is irrelevant for the argument.
The same idea can be used to obtain even more complicated $NH(x)$.
Recall that a family of sets is {\em almost disjoint} iff the intersection of any two members is finite.
Call a space $\Psi$-like if it is the union of a countable set (which we may identify with $\omega$)
and an almost disjoint family
$\mathcal{N}$ of infinite subsets of $\omega$, where the points in $\omega$ are isolated and a
neighborhood of $p\in\mathcal{N}$ is $\{p\}$ union $p-F$, 
where $p$ is seen as a subset of $\omega$ and $F\subset\omega$ is finite.
Any $\Psi$-like space is locally compact and separable.
If $\mathcal{N}$ is maximal (for the inclusion) almost disjoint, that is, if any 
infinite subset of $\omega$ not in 
$\mathcal{N}$ has infinite intersection with each member of $\mathcal{N}$, then the space is 
called a $\Psi$-space. Every $\Psi$-space is pseudocompact, as well known. Their original construction
is due to Isbell and Mr\'owka, for more on these spaces, see e.g. \cite{Mrowka:1977}.

\begin{example}
   \label{ex:NH-psi2}
   For any $\Psi$-like space $Y = \omega\cup\mathcal{N}$, there is
   a NH-$1$-manifold $M$ with $\mathcal{E}(M)= |\mathcal{N}|$ and a point 
   $x$ for which $NH(x)$ is a copy of $Y$.
\end{example}
\begin{proof}[Details]
   Our space will be a quotient of the disjoint union of copies $R_p$ of $\R$ for each $p\in\mathcal{N}$
   together with two other disjoint copies $R_A,R_B$ of $\R$.
   For any point $x\in\R$, write $x_p$, $x_A$, $x_B$ for its respective copies in $R_p,R_A,R_B$,
   and the same for subsets of $\R$.
   If $p\in\mathcal{N}$, then $p$ is an infinite subset of $\omega$, and 
   for each $n\in p$, we identify the interval
   $(\frac{1}{n+2},\frac{1}{n+1})_p\subset R_p$ with the interval $(2n-1,2n+1)_A\subset A$, for instance
   with an affine homeomorphism. Then $0_p$ is a limit of the 
   sequence $\{2n_A\,:\,n\in p\}$. Moreover, if $p,q$ are distinct members of $\mathcal{N}$,
   then their intersection is finite and hence $0_p$ and $0_q$ may be separated in the space.\\
   We now wish to make $2n_A\bumpeq 0_B$ for each $n\in\omega$, which
   automatically yields $0_p\bumpeq 0_B$ for each $p\in\mathcal{N}$ as well 
   by Lemma \ref{lemma:NHseq} (take $y_i=0_B$ for each $i$).
   To achieve this, fix as before a bijection $\sigma:\omega\to\omega\times\omega$.
   Set $I^A_{n,m}$ to be the interval $(2n + \frac{1}{m+2}, 2n + \frac{1}{m+1})_A$.
   Identify $I^A_{\sigma(n)}$ with $(\frac{1}{n+2},\frac{1}{n+1})_B$ for each $n$ with an affine homeomorphism.
   (Of course, this identifies subintervals of $(\frac{1}{n+2},\frac{1}{n+1})_p$ with 
   subintervals of $(\frac{1}{n+2},\frac{1}{n+1})_B$ as well.)
   Each neighborhood of $0_B$ intersects 
   $(\frac{1}{n+2},\frac{1}{n+1})_B$ for all but finitely many $n$,
   hence $0_B\bumpeq 2n_A$ for each $n$.
   If $x\not= 0$ and $y\not=2n$, then both $x_p$ and $y_A$ can be separated from $0_B$. 
   This shows that $NH(0_B) = \{2n_A\,:\,n\in\omega\}\cup \{0_p\,:\,p\in\mathcal{N}\}$.
   The induced topology on $NH(0_B)$ is exactly that of $Y$ if we interpret
   $\{2n_A\,:\,n\in\omega\}$ as $\omega$ and $\{0_p\,:\,p\in\mathcal{N}\}$ as $\mathcal{N}$.
\end{proof}

A variation of the same idea yields a $NH(x)$ homeomorphic to a special tree, as we shall now see.
Recall that a (set theoretic) {\em tree} is a set 
$T$ partially ordered by $\le$ such that the subset of predecessors of each
member is well ordered. For $\alpha$ an ordinal, 
we say that $x\in\text{Lev}_\alpha(T)$ iff the predecessors of $x$ in $T$
have order type $\alpha$, and the height of $T$ is the smallest $\alpha$ such that $\text{Lev}_\alpha(T) = \varnothing$.
An {\em antichain} is a subset of $T$ containing pairwise $\le$-incomparable elements.
A tree is {\em special} iff it is the union of countably many antichains,
and {\em rooted} if it contains a unique minimal element, called the root. 
It is well known that special trees of height $\omega_1$ do exist,
some of them having only countable levels.
Given a rooted tree $T$ with root $r$, 
we let $H_T$ be $T\times[0,1) - \{\langle r,0\rangle\}$ with lexicographic order topology.
We abuse the notation and write $\le$ as well for the lexicographic order on $H_T$.

\begin{lemma}
   Let $T$ be a rooted tree of height $\le\omega_1$.
   Then $H_T$ is a $1$-manifold, $T$ embeds in $H_T$ and if $\le$ is not a total order, then $H_T$ is non-Hausdorff.
\end{lemma}
\begin{proof}
   To show that $H_T$ is a $1$-manifold, one argues exactly as in the proof that $\LL_+$
   is a Hausdorff $1$-manifold (see e.g. \cite[Lemma 1.10]{GauldBook}).
   Notice that $T$ embeds in $H_T$ as $T\times\{0\}$ (with the exception of the
   root $r$, that we may embed as $\langle r,\frac{1}{2}\rangle$).
   If $\le$ is not a total order, let $\alpha$ be minimal such that there are incomparable $x,y\in\text{Lev}_\alpha(T)$.
   So there is an homeomorphism $\varphi:\{\langle x,t\rangle\in H_T\,:\,x\text{ at level }<\alpha\}\to(0,1)$,
   and given two neighborhoods of $x$ and $y$, there is $\epsilon>0$ such that both contain 
   $\varphi^{-1}((1-\epsilon,1))$. 
\end{proof}

\begin{example}
  \label{ex:NHspecialtree}
  Given a rooted special tree $T$, there is 
  a NH-$1$-manifold $M$ with $\mathcal{E}(M)=|T|$ such that there is $x\in M$ with
  $NH(x)$ homeomorphic to $T$.
\end{example}
\begin{proof}[Details]
  Let $T$ be a special tree with root $r$, then $T$ has height $\le\omega_1$, as easily seen.
  Define $\wt{T}$ as $T$ with an added
  point $\wt{r}$ below the root, so that $T$ embeds in $H_{\wt{T}}$
  as $T\times\{0\}$.
  Let $A_n$ be disjoint antichains in $\wt{T}$ such that $\cup_{n\in\omega}A_n = \wt{T}$.
  We define $M$ as $\{0\}\times \R\cup H_{\wt{T}}$, quotiented as follows.
  Let $I_n\subset (0,1)$ be the interval $(\frac{1}{2^{n+1}},\frac{1}{2^{n}})$ and
  fix a bijection $\sigma:\omega\times\omega\to\omega$. 
  If $x\in T$, let $J_x\subset H_{\wt{T}}$ be the interval $\{\langle y,t\rangle\in H_T\,:\,y\le x,\,0\le t<1\}$.
  If $x\in A_n$, we identify with an affine homeomorphism the
  interval $\{x\}\times I_m \subset J_x$ with $\{0\}\times I_{\sigma(n,m)}$.
  Since each $A_n$ is an antichain, when $y<z\le x$ we have that $\{y\}\times I_m$ and 
  $\{z\}\times I_k$ are identified with disjoint subintervals of $\{0\}\times\R$, 
  hence $J_x$ is still homeomorphic to an interval
  in the quotient space. Moreover, no points of $T\times\{0\}\subset H_{\wt{T}}$ are 
  identified and the induced topology
  is still that of $T$. If $x\in T$, every neighborhood of $\langle x,0\rangle$ contains
  $\{x\}\times I_m$ for infinitely many $m$, hence $\langle 0,0\rangle\bumpeq \langle x,0\rangle$.
  Any other point of $M$ can be separated from $\langle 0,0\rangle$, hence 
  $NH(\langle 0,0\rangle) = T\times\{0\}\simeq T$.
\end{proof}

Notice that all of $H_T$ is identified with points of $\{0\}\times\R$ except those
$\langle x,t\rangle$ where $t$ is either $0$ or of the form $\frac{1}{2^{n+1}}$.
In $H_{\wt{T}}$, the interval $\{\wt{r}\}\times(0,1)$ 
is not identified, but we may set 
$\langle \wt{r},t\rangle \sim \langle 0,t-1\rangle$
so that in the resulting quotient space, $\{0\}\times\R$ is H-maximal and CH-maximal.
Our last $1$-dimensional example has much simpler $NH(x)$ but also more homeomorphisms. 

\begin{example}
  \label{ex:NHnotdiscrete}
  A NH-$1$-manifold $M$ with $\mathcal{E}(M)=2$ such that the following hold.\\
  $\bullet$ There are $x,y\in M$ such that $NH(x)$ and $NH(y)$ are homeomorphic to
  $\omega+1$, $y$ is the limit point of $NH(x)$ and $x$ that of $NH(y)$. \\
  $\bullet$ There is a homeomorphism which interchanges $x$ and $y$.
\end{example}
\begin{figure}[h]
     \begin{center}
        \epsfig{figure=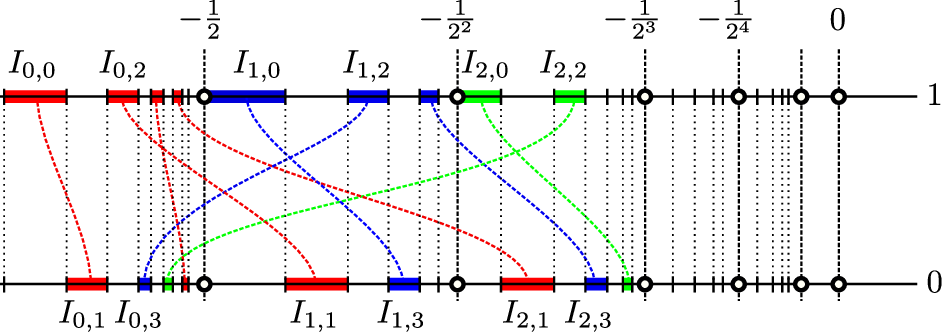, width=.8\textwidth}
         \caption{The identifications in Example \ref{ex:NHnotdiscrete}. Intervals linked by a dashed line
         are identified. The identifications of non-colored intervals
         (not shown to avoid overcharging the picture) are given by the horizontal 
         mirror image, for instance $I_{0,0}\times\{0\}$ is identified 
         with $I_{0,1}\times\{1\}$, etc.}
         \label{fig:NHnotdiscrete}
      \end{center}
\end{figure}
\begin{proof}[Details]
  Define the intervals 
  $I_{n,m}=\left(-\frac{1}{2^{n+1}}(1+\frac{1}{2^{m}}) , -\frac{1}{2^{n+1}}(1+\frac{1}{2^{m+1}})\right)\subset\R$.
  Notice that the left boundary of $I_{n,0}$ is $- \frac{1}{2^n}$,
  that the $I_{n,m}$ accumulate to $-\frac{1}{2^{n+1}}$ when $m\to\infty$ and 
  to $0$ when $n\to\infty$.\\
  $M$ will be $(-1,1)\times\{0,1\}$ with $I_{n,m}\times\{0\}$ identified to some 
  $I_{n',m'}\times\{1\}$, in such a way that in the resulting space, the $I_{n,m}\times\{0\}$ accumulate to
  both $\langle -\frac{1}{2^n},0\rangle$ and $\langle 0,1\rangle$ when $m\to\infty$; and reciprocally, the
  $I_{n,m}\times\{1\}$ accumulate to
  both $\langle -\frac{1}{2^n},1\rangle$ and $\langle 0,0\rangle$ when $m\to\infty$.
  This ensures that $NH(\langle 0,0\rangle)= \{\langle -\frac{1}{2^n},1\rangle\,:\,n\in\omega\}\cup \{\langle 0,1\rangle\}$,
  which is of course homeomorphic to $\omega+1$.
  The identification is done symmetrically: if $I_{n,m}\times\{0\}$ is identified to 
  $I_{n',m'}\times\{1\}$, then $I_{n,m}\times\{1\}$ is identified to 
  $I_{n',m'}\times\{0\}$. This implies that the map $h:M\to M$ which flips the levels $(-1,1)\times\{0\}$
  and $(-1,1)\times\{1\}$ is well defined and a homeomorphism. Note that $h$ interchanges
  $\langle 0,0\rangle$ and $\langle 0,1\rangle$, 
  and that $\langle 0,1-i\rangle$ 
  is the limit point of $NH(\langle 0,i\rangle)$ ($i=0,1$).
  Figure \ref{fig:NHnotdiscrete} tries to describe graphically the situation.
  \\
  To make precise the identification rules, let $A_n$ ($n\in\omega$) be a partition of $\omega$ into infinite 
  subsets, and let $a_{n,m}$ be an increasing enumeration of $A_n$.
  Fix a bijection $\tau:\omega\to\omega\times\omega$, with $\tau(n) = \langle \tau_0(n),\tau_1(n)\rangle$.
  We now define the bijection $\sigma:\omega\times\omega\to\omega\times\omega$ as
  $$ \sigma(n,m) = \langle \sigma_0(n,m) , \sigma_1(n,m)\rangle = \langle \tau_0(m) , a_{n,\tau_1(m)}\rangle.$$
  If $n$ is fixed, then $\sigma(n,m)$ is cofinal in each $\{k\}\times\omega$ for each $k$ when $m\to\infty$.
  Then, identify $I_{n,2m}\times\{0\}$ with $I_{\sigma_0(n,m) , 2\sigma_1(n,m) + 1}\times\{1\}$,
  and reciprocally identify 
  $I_{n,2m}\times\{1\}$ with $I_{\sigma_0(n,m) , 2\sigma_1(n,m) + 1}\times\{0\}$. In words:
  we identify (through $\sigma$) 
  the $I_{n,m}$ on one level with even $m$ with $I_{n',m'}$ one the other level with odd $m'$. 
  \\
  We can make the same identifications to the right of the $0$ points on each level by simply changing the
  signs in all the intervals $I_{n,m}$, so that the two levels have only countably many points not in common. 
\end{proof}
The manifold just described is quite far from being homogeneous,
and $\langle 0,0\rangle$ does not belong to any sorted neighborhood nor has any f-rim-simple neighborhoods. 

\vskip .3cm\noindent
Notice that in every example in this subsection, we can make further identifications of intervals
so that in the resulting $1$-manifold, the $0$-th floor 
(or the $B$ floor in \ref{ex:NH-psi2}) is actually CH-maximal and H-maximal. This is already the case 
in Example \ref{ex:NHspecialtree}.


\section{Examples of non-homogeneous NH-$2$-manifolds.}
\label{sec:2dim}

We saw that NH-$1$-manifolds with $\mathcal{OH}(M) = 2$ may be quite complicated.
This worsens in dimension $2$.
We try to stick to examples for which we can draw pictures
whenever possible (and give the more abstract variations after the picturable ones). 
A simple way to obtain a non-discrete $NH(x)$ is the following example, which we give in any dimension.
\begin{example}
  \label{ex:NHS1}
  An $(n+1)$-manifold $M$ with $\mathcal{E}(M) = 2$ and a point $p$ such that $NH(p)$ is homeomorphic to 
  the $n$-sphere $\mathbb{S}^n$.
\end{example}
\begin{proof}[Details]
   The idea is to use the fact that if one identifies all the points in a
   closed ball in $\R^{n+1}$, the quotient space is still homeomorphic to $\R^{n+1}$.
   Write $0$ for the origin in $\R^{n+1}$.
   Take $M=\R^{n+1}\cup\{0^*\}$, and define the neighborhoods of $0^*$ as $\{0^*\}\cup (B(0,r)-\wb{B(0,1)})$
   for any $r>1$, where
   $B(x,r)$ is the ball of radius $r$ centered at $x$. (Points in $\R^{n+1}$ have their usual neighborhoods.)
   Then $NH(0^*)$ is the boundary of $B(0,1)$, that is the copy of $\mathbb{S}^n$ of radius $1$ centered at the origin.
   \\
   Another way of defining $M$ is with a quotient:
   $M = \R^{n+1}\times\{0,1\}/\sim$, where $\langle x,0\rangle\sim\langle y,1\rangle$ iff $x=y$ and $|x|>1$,
   and $\langle x,1\rangle\sim\langle y,1\rangle$ iff $x,y\in\wb{B(0,1)}$. 
\end{proof}

There is a simple generalization of Example \ref{ex:NHS1} that 
yields more complicated $NH(x)$ in dimension $2$. 
We use a 100 years old theorem of R.L. Moore, whose original source is \cite[Theorem 25]{Moore:1925}.
As usual, $\partial A$ denotes the boundary $\wb{A}-\text{int}(A)$.

\begin{example}
   \label{ex:MooreThm}
   Let $C\subset\R^2$ be compact connected with connected complement. 
   Then there is a NH-$2$-manifold $M$ with $\mathcal{E}(M)=2$
   and $x\in M$ such that $NH(x)$ is homeomorphic to $\partial C$.
\end{example}
\begin{proof}[Details]
   Moore Theorem implies that $\R^2/C$ (i.e. the quotient space where points in $C$ are identified) 
   is homeomorphic to $\R^2$. 
   In particular, there is a homeomorphism $\phi:\R^2-\{0\}\to\R^2-C$. 
   We may then let $M$ be $\R^2\cup\{0^*\}$, with neighborhoods of $\{0^*\}$
   given by $\{0^*\}\cup\phi( B(0,r)-\{0\})$ for $r>0$. Then $0^*$ cannot cannot be
   separated from the points in the boundary of $C$. 
\end{proof}

Notice that this example (and Example \ref{ex:NHS1} in dimension $2$, which is a particular case) 
does have a CH-maximal sorted neighborhood:
set $G = \{0^*\}\cup \R^2- C$, then 
the map $g:M\to G$ defined as $g\upharpoonright (\R^2-C) = \phi$, $g\upharpoonright C \equiv 0^*$, $g(0^*)=0^*$ 
is a sorting map. However, if $C$ is not trivial (i.e. a point),
$\partial C\cap U$ is uncountable for any open set $U$, hence
no neighborhood of any point of $\partial C$ is weakly sorted
by Lemma \ref{lemma:sortBumpeq}.
$G$ is H-maximal iff $C$ has no interior.
If we repeat the construction and add new points $x^*$ for each $x\in\R^2$ simply by translating the neighborhoods,
then there can be no sorted neighborhood containing $x^*$ and $y^*$ if $NH(x^*)\cap NH(y^*)$ is non-empty 
(by Proposition \ref{prop:NHdisjoint}). If $C$ is not ``too simple'', there will be such intersections.

In Example \ref{ex:MooreThm}, $NH(x)$ is a singleton for any $x\in\R^2$,
but $NH(0^*)$ is complicated.
The constructions in the next two subsections can be thought as generalizations of this example,
but reversing the roles: we start with the plane and add new points with complicated 
neighborhoods, so that exactly one point $x$ in the plane ends up with a complicated $NH(x)$, while the 
added points $y$ all have $NH(y)=\{x\}$. This then enables to have a sorted neighborhood 
containing all the (initial) 
plane.\footnote{We admit not being sure that it is a desirable feature, but think that the examples
in question are interesting in their own, so we chose to keep them.}


\subsection{Pr\"ufer-like constructions}\label{subsubsection:Prufer}

We start with adaptations of known constructions of Hausdorff surfaces.\footnote{The constructions 
in this subsection (except Example \ref{ex:Nyikosize}) can be seen as
particular cases of those described in the next one using conformal theory, but we prefer 
to describe them on their own.}

\begin{example}
   \label{ex:BadPrufer}
   A NH-$2$-manifold \Prufer with $\mathcal{E}(\text{\Prufer}) = 2$ and a point $p$ such that 
   $NH(p)$ is homeomorphic to $\R$.
   Moreover, \Prufer has a CH-maximal sorted neighborhood and a base of rim-simple open Hausdorff sets.
\end{example}
\begin{proof}[Details]
We first remark that if in the following construction we remove the point
$p\in R_0$ defined below, what we do is exactly the process of
Pr\"uferizing a surface at some point
(see e.g. \cite[Example 1.25]{GauldBook}).
In a sense, we follow the same process but do it badly, so that
instead of obtaining an honest Hausdorf manifold we end up with a non-Hausdorff one.
 
\begin{figure}[h]
     \begin{center}
        \epsfig{figure=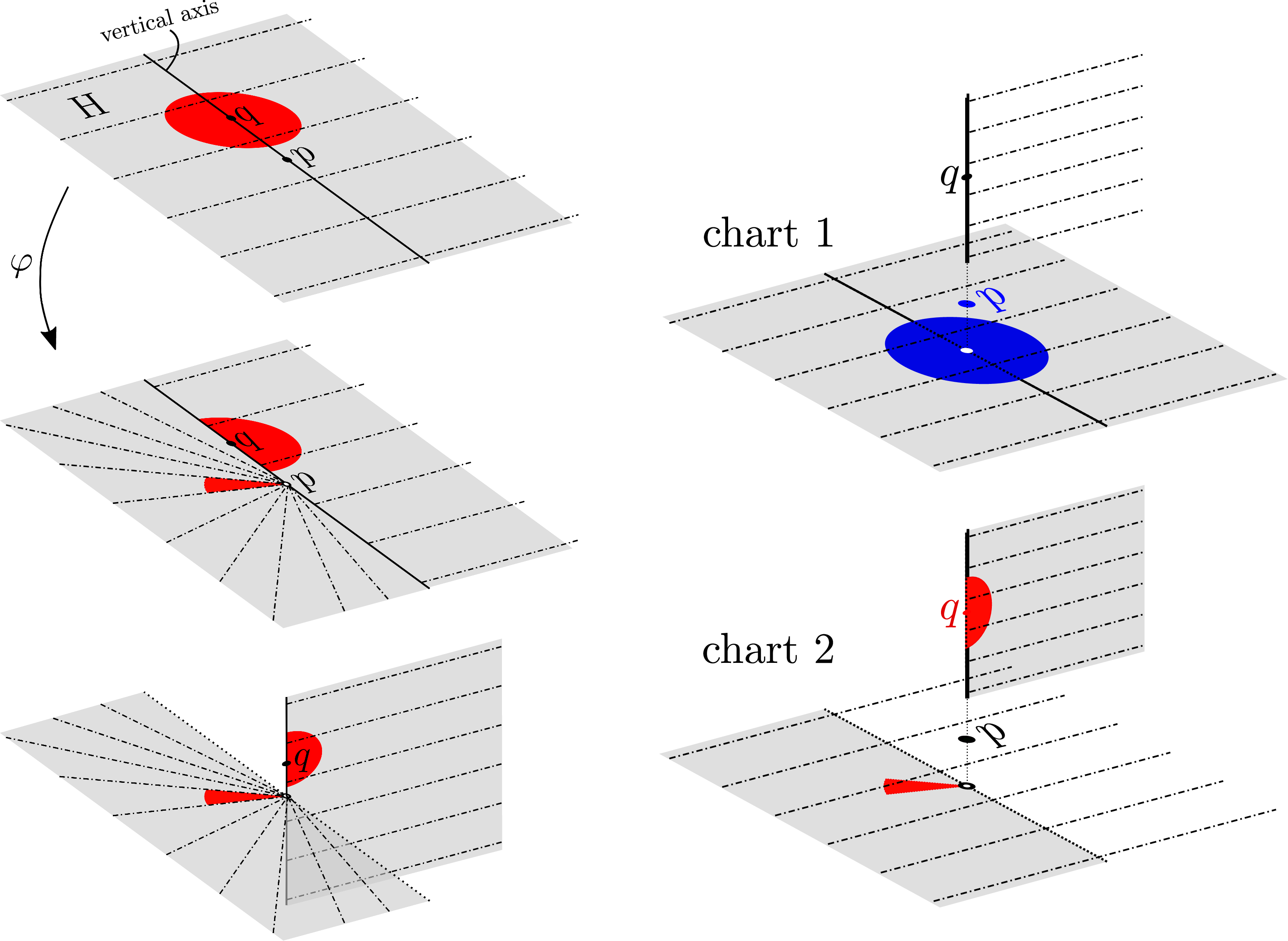, width=.9\textwidth}
         \caption{A NH-$2$-manifold such that $NH(p)$ is a copy of $\R$.}
         \label{fig:5}
      \end{center}
\end{figure}

\noindent
First, consider a point $p$ in the vertical axis in the plane.
Set $H = (-\infty,0)\times\R$ and 
let $\varphi:H\to H$ be the homeomorphism $(x,y)\mapsto (x, yx + c)$ (for some $c\in\R$).
This sends a horizontal line at height $y$ in $H$ into a line of slope $y$ pointing at $p=\langle 0,c\rangle$.
This homeomorphism transforms a neighborhood of $q = \langle 0 , y\rangle$
as depicted in the upper left part of Figure \ref{fig:5}. The bottom left picture is another way to look at the same
plane which is fit for our later purposes,
and which explains the symbol \Prufer for our manifold.\\
Now, take two copies $P_0,P_1$ of the plane 
with their two copies $H_0,H_1$ of $H$ and identify $(x,y)\in H_1$ with $\varphi(x,y)$ in $H_0$.
This yields a NH-$2$-manifold, where the point $p\in P_0$ cannot be separated from any point
in the vertical axis in $P_1$.
The righthandside of Figure \ref{fig:5} tries to picture the 
entire manifold and its two charts, with the neighborhoods
of $p\in P_0$ and $q$ in the vertical axis of $P_1$ 
are shown (in their respective charts) in blue and red, respectively. Since there is 
a kind of ``branching'' at $p$ we lifted it from the  bottom plane. 
By construction, only the points in the vertical axis of $P_1$ cannot be separated
from $p\in P_0$.\\
The claim about having a rim-simple basis should be clear for every points
except maybe $p$ and those in $NH(p)$.
But looking at Figure \ref{fig:5} without blinking during $10$ seconds should convince
the reader that it also holds for these points, as every intersection of neighborhoods
depicted there yields connected sets. $P_0$ is a sorted neighborhood in \Prufer: 
send all of $P_1-H_1$ to the point $p$.
\end{proof}

\begin{rem}
\label{rem:Prufer}
\
\\
(a)
Given two points in $NH(p)$, there is a homeomorphism sending one to the other: just use a 
vertical translation $\langle x,y\rangle\mapsto\langle x,y+u\rangle$ in $P_1$, 
set $\langle x,y\rangle\mapsto \langle x, y +xu\rangle$ in $H_0$
and fix the rest of $P_0$. This sends a neighborhood of $p$ to another neighborhood of $p$.\\
(b) We could also identify the copies of $(0,+\infty)\times\R$ pointwise
(or again through $\varphi$), so that 
all of $P_0,P_1$ except their vertical axis are identified.
\end{rem}

It is possible to push the idea further.
To simplify, let us assume that $p$ is the origin in $P_0$.
Suppose now that one deletes $(0,+\infty)\times\R$ in $P_1$ in Example \ref{ex:BadPrufer}.
Hence, a member of the vertical axis in $P_1$ has now neighborhoods homeomorphic to $(-\infty,0]\times\R$,
as for (Hausdorff) surfaces with boundary. Now, identify $\langle 0,x\rangle \in P_1$ with 
$\langle 0,-x\rangle \in P_1$.\footnote{In an honest Hausdorff surface, this is known as {\em Mooreizing}
this boundary component, see e.g. {\cite[Ex. 1.27]{GauldBook}}.}
This has the effect of folding the vertical axis of $P_1$ around the origin, 
and its points now have
true Euclidean neighborhoods.

\begin{figure}[h]
     \begin{center}
        \epsfig{figure=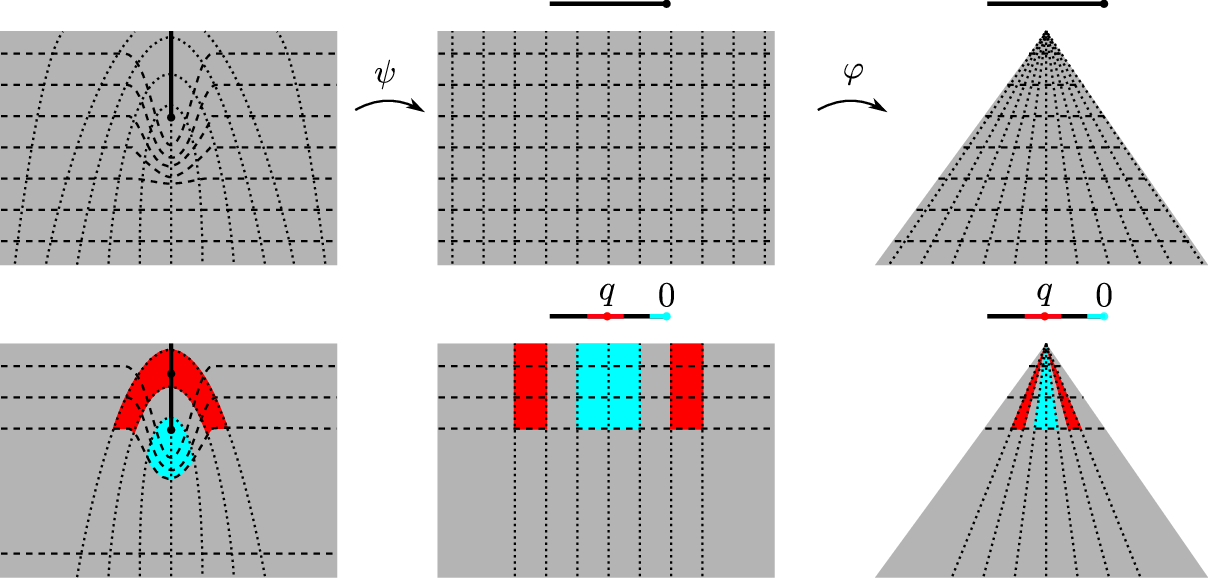, width=.9\textwidth}
         \caption{How to Moorize the plane at a point.}
         \label{fig:Pruferize_Square}
      \end{center}
\end{figure}

Another way of understanding this trick is to look at Figure \ref{fig:Pruferize_Square}, where
the image is turned at $90^\circ$. It describes how to insert a half open interval at a given
point. The maps $\psi$ and $\varphi$ are homeomorphisms, $\varphi$
being as in Example \ref{ex:BadPrufer}.
In cyan is shown a neighborhood of the origin in the folded vertical axis, while
a neighborhood of another point $q$ is shown in red. 

\begin{figure}[h]
     \begin{center}
        \epsfig{figure=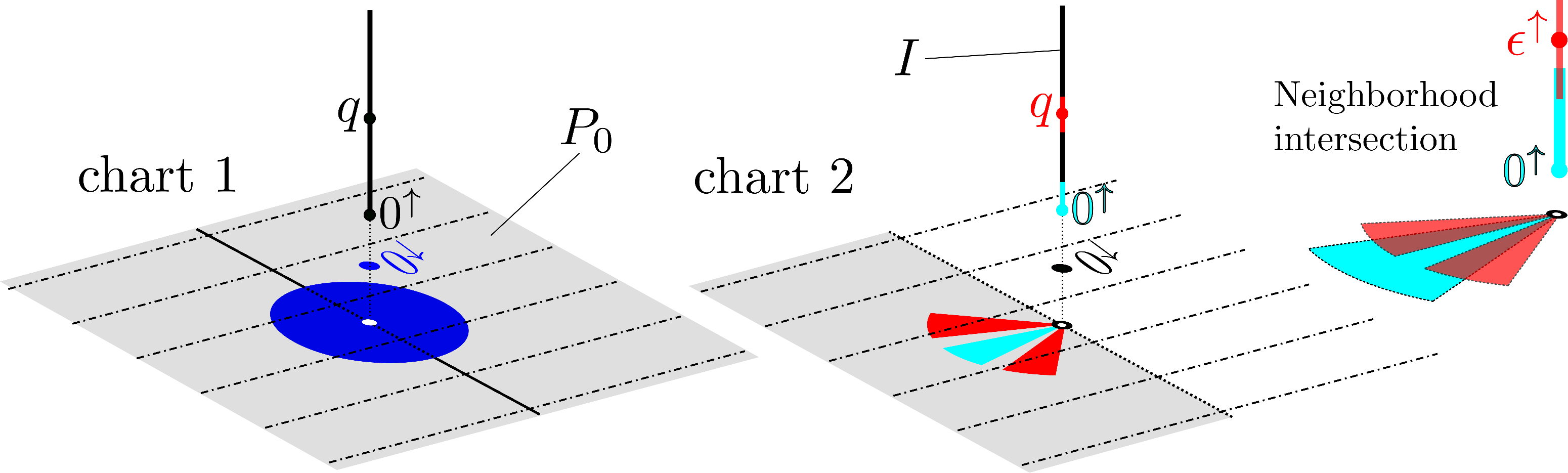, width=.9\textwidth}
         \caption{Example \ref{ex:BadMoore}.}
         \label{fig:BadMoore}
      \end{center}
\end{figure}

The resulting manifold thus consists of the plane $P_0$ (thought as downstairs)
together with a half open interval $I$ (upstairs), which may be thought as $[0,+\infty)$. 
The $0$ point of the vertical axis downstairs is written $0^\downarrow$, and if
$a\in[0,+\infty)$, we write $a^\uparrow$ for the corresponding point of $I$.
Neighborhoods of
$0^\downarrow,0^\uparrow$ and any other point $q\in I$ are shown on Figure \ref{fig:BadMoore}
with respective colors blue, cyan and red.
Other points downstairs' neighborhoods are standard balls in $P_0$. The manifold just described 
yields then the following.

\begin{example}
   \label{ex:BadMoore}
   There is a 
   NH-$2$-manifold \Moore and points $u,v\in\text{\Moore}$ with the following properties.
   \\
   $\bullet$ \Moore has a H-maximal sorted neighborhood which contains $u$;\\
   $\bullet$ $NH(u)$ is a copy of $[0,+\infty)$;\\
   $\bullet$ there is a basis $\mathcal{B}$ of Euclidean neighborhoods
   at $v$ which are rim-simple;\\
   $\bullet$ for each member $U$ of $\mathcal{B}$, there is a point $q\in U$
   and an Euclidean neighborhood $V$ of $q$ contained in $U$ such that any open $W\subset V$
   containing $q$ is not rim-simple at some point of its boundary.  
\end{example}
Hence, in particular, there is no basis for the topology of \Moore which consists of open sets
that are rim-simple, and the fact that an open set is rim-simple does not prevent
the existence of subsets that are ``complicated'' at some point of their boundary. 
\begin{proof}[Details]
   $P_0$ is again a sorted neighborhood in \Moore for the same reason as before.
   Take $u$ to be $0^\downarrow$ and $v$ to be $0^\uparrow$. 
   A base for the neighborhoods of $0^\uparrow$ is given by triangles pointing at the $0$ point in the vertical 
   axis in $H_0$, with side's slope $\pm\epsilon$, 
   symmetric with respect to the horizontal axis, together with an interval $[0^\uparrow,\epsilon^\uparrow)$ in $I$.
   Let $U$ be such a neighborhood.
   Since any member of $I$ is contained in $H_0\cup I$, which is Euclidean, $U$
   is simple at the only point in its boundary intersected with $I$, that is: $\epsilon^\uparrow$.
   (Note: one may first think that it is not the case, as the
   their intersection with $P_0$ gives a disconnected bowtie-like figure as in Figure \ref{fig:BadMoore} (right).
   But these two parts are connected together in $I$, as seen in Figure \ref{fig:Pruferize_Square},
   going right to left.)\\
   Moreover, $U\cap P_0$ being convex, $U$ is simple at any point in $(\wb{U}-U)\cap P_0$.
   Hence, $U$ is rim-simple. 
   Take now $a$ small enough so that $a^\uparrow\in U$, and let $V$ be a basic neighborhood of $a^\uparrow$
   not containing $0^\uparrow$. By reducing $V$ further, we may assume that its intersection with $U$
   is the union of two triangles (in the form of a partially folded bowtie). Since $a^\uparrow\bumpeq 0^\downarrow$
   by construction, $0^\downarrow$ is in the boundary of $V$.
   But no neighborhood of $0^\downarrow$ has a connected intersection with $V$ (see again Figure 
   \ref{fig:BadMoore}).
   (This also shows in passing that no Euclidean neighborhood of $0^\downarrow$ is rim-simple.)
\end{proof}

We can push the idea further. The long line $\LL$ is obtained by taking two copies
of $\LL_{\ge 0}$ and gluing them at their $0$ point.

\begin{example} 
   \label{ex:Nyikosize}
   NH-$2$-manifolds $M$ with $\mathcal{OH}(M) = 2$ and a point $p$ such that $NH(p)$ is homeomorphic to the
   either $\LL_+$ or $\LL_{\ge 0}$ or $\LL$. 
   Moreover, $M$ has a H-maximal sorted neighborhood containing $p$.
\end{example}
This time we do not delve into the rim-simplicity of neighborhoods.
\begin{proof}[Outline of the construction]
There is a procedure similar to that of Pr\"uferizing 
(and/or Moorizing) at a point $p$ in the plane which is called
{\em Nyikosizing} at $p$,
described in detail in \cite[Example 1.29]{GauldBook}.
This has the effect of inserting a copy $L$ of either $\LL_+$, $\LL_{\ge 0}$
or $\LL$ in (the boundary of) 
the plane, and the neighborhoods of points in $L$ are ``wedges'' or ``bowties''
(similar to the neighborhoods of the point $q$ in Figure \ref{fig:5}) between curves
in $H$ that ``point'' to $p$.
Now, take again two copies $P_0,P_1$ of the plane, and Nyikosize $P_1$ at some point $p$ in the vertical axis.
Then identify pointwise $H_0\subset P_0$ with $H_1\subset P_1$. (This time, we do not need any homeomorphism
between $H_0$ and $H_1$, because the work of ``putting $L$ at $p$'' is done by the Nyikosization process.)
Then as above, any neighborhood of $p\in P_0$ will intersect any curve pointing at $p$ in $H$,
which means that no point of $L\subset P_1$ can be separated from $p$.
The sorting map is given by sending all of $NH(p)$ to $p$, and the identity elsewhere.
\end{proof}

Notice that the constructions in the previous three examples 
can be performed at any $p$ in the plane, without interfering with each other: we just define the neighborhoods of the
points in the new added piece above $p$ by specifying how they intersect the plane and this new added piece.
This yields the following.
\begin{example}
   \label{ex:BadAll}
   Given a (non-continuous) map $\sigma:\R^2\to\{\{*\},\R,\R_{\ge 0},\LL_+,\LL_{\ge 0},\LL\}$ 
   (where $\{*\}$ denotes a singleton),
   there is an ENH-$2$-manifold $M$ which consists of $\R^2\cup\cup_{p\in\R^2}NH(p)$,
   where $NH(p)$ is a copy of $\sigma(p)$ for each $p$, and $NH(x)=\{p\}$ for each $x\in NH(p)$.
   Moreover, the map that sends $NH(p)$ to $p$ (and the identity in $\R^2$) is a sorting map.
\end{example}
Notice that for $\R$, we need to mirror of the construction in \Prufer around the vertical axis as in point (b) of 
Remark \ref{rem:Prufer}.

\subsection{Constructions using conformal theory}\label{subsubsection:conformal}

Start with a bounded domain (i.e. open connected and simply connected subset of the complex plane) $U$,
where the boundary $\partial U = \wb{U}-U$ is not ``too complicated''.
With the help of Riemann and Caratheodory's conformal theory, it is possible
to build a $2$-manifold with a sorted neighborhood
such that $NH(x)$ is a copy of a part of $\partial U$ for some point $x$.
First, Riemann mapping theorem tells us that $U$ is conformally equivalent to the interior $D^\circ$
of the closed unit disk $D$, that is, there is a conformal homeomorphism
$\phi:U\to D^\circ$. If the boundary $\partial U = \wb{U}-U$ is moreover locally connected,
Caratheodory's theorem yields that
the inverse of $\phi$ extends continuously to a map $\psi:D\to\wb{U}$ (see e.g. \cite[Theorem 16.6]{MilnorDynamics:1990}).

\begin{figure}
     \begin{center}
        \epsfig{figure=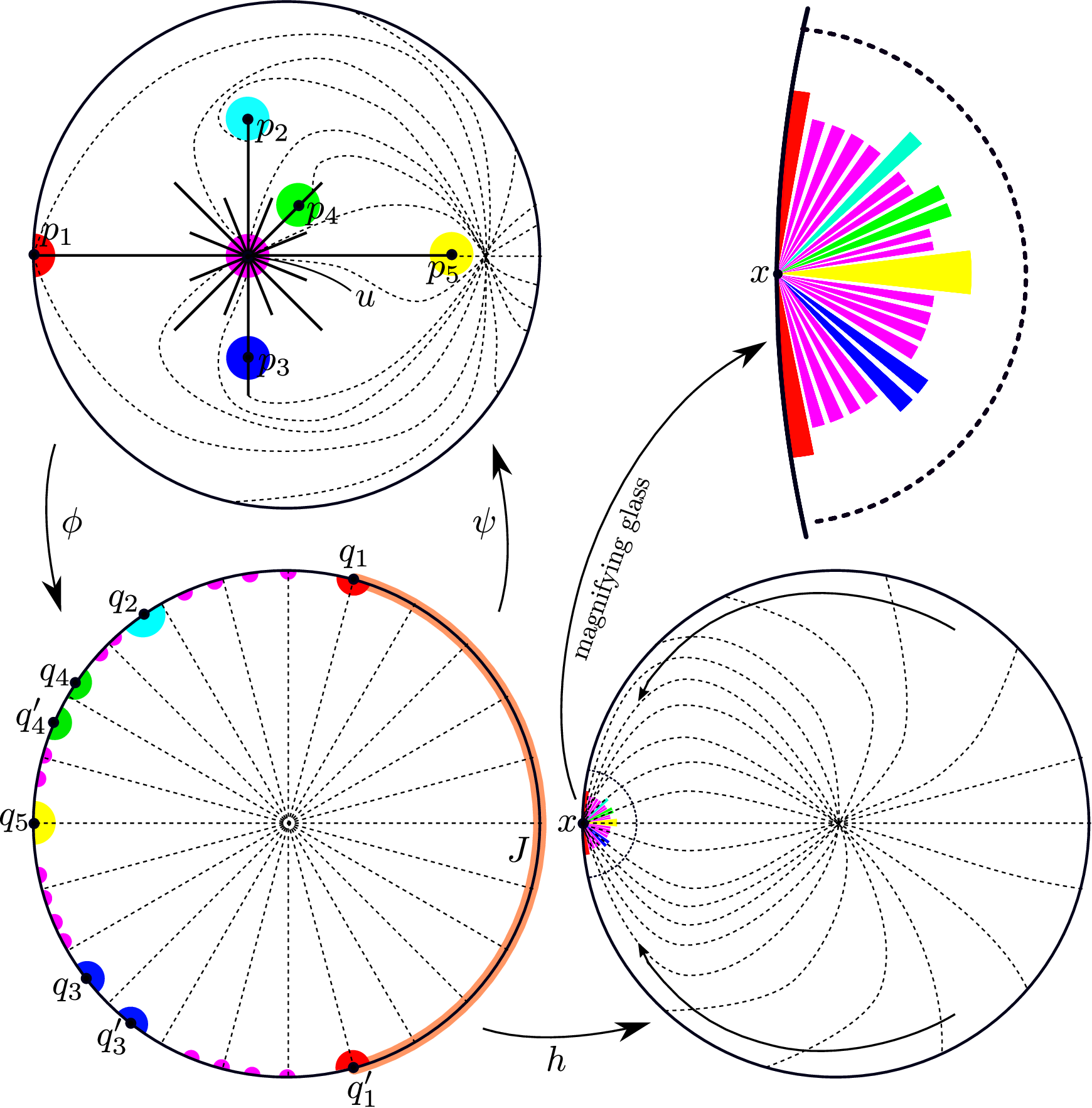, width=.75\textwidth}
         \caption{Pinching a conformal map: a domain with locally connected boundary.}
         \label{fig:caratheodory}
      \end{center}
\end{figure}

\begin{figure}
     \begin{center}
        \epsfig{figure=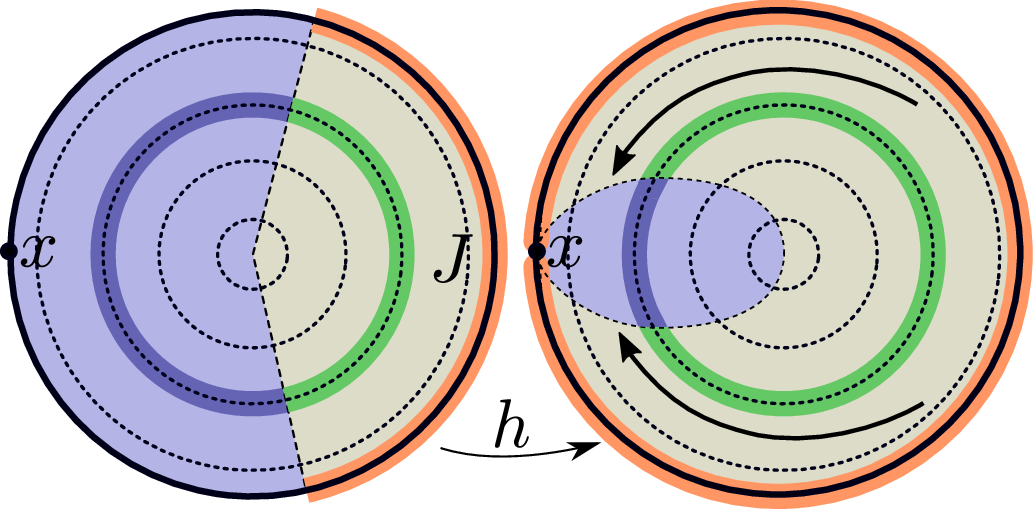, width=.45\textwidth}
         \caption{The pinching map $h:D\to D$.}
         \label{fig:h}
      \end{center}
\end{figure}

The constructions we are about to describe are possible for any domain
whose boundary is ``sufficiently nice'', but rather than
trying to formally state what we mean by that, 
we shall concentrate on two particular examples.
The first is shown in Figure \ref{fig:caratheodory},
to which the reader is refered througout the following description.
We start with a domain $U$ defined by the interior of a disk with a kind of ``infinite crossing star'' 
deleted. In words: 
Start with a horizontal segment of length $1$, take another vertical one of length $\frac{1}{2}$
and make them intersect at their center. Then, bissect each created angle with segments of half
length again crossing at their center, and so on. 
This set is locally connected (due to the fact that the added segments
have length going to $0$), call it \starbig.
Attach then the leftmost point of \starbig { }
to the boundary of the disk,
see the top left part of Figure \ref{fig:caratheodory}. By construction, $U\subset D^\circ$ and $\wb{U}=D$.
Caratheodory's theorem tells us that the inverse $\psi$ of the Riemann map $\phi:U\to D^\circ$ extends continuously to 
all of $D$, defining $\psi:D\to \wb{U}=D$. 
We have pictured neighborhoods of some points of \starbig, and how they are mapped under $\phi$.
Points $p_i$ have one or two preimages $q_i,q_i'$ under $\psi$, while $u$,
the point ``in the middle of the star'', has infinitely many (actually,
its preimages form a Cantor set). The number of connected components of a neighborhood of $u$ minus \starbig { }is 
finite but unbounded by construction. 
\\
The preimages of \starbig{ }under $\psi$ do not cover the entire boundary of $D$, since
by e.g. \cite[Theorem 15.5]{MilnorDynamics:1990}, $\phi$ extends homeomorphically to
the boundary circle minus $p_1$, whose image is an open interval $J$.
We may thus define a continuous map $h$ of the boundary $\partial D$ to itself
which stretches $J$ so that it covers the entire boundary $\partial D$ except one point
which we call $x$, and extend it to all of $D$ such that $h\upharpoonright D^\circ$ is an homeomorphism of $D^\circ$,
and $h\upharpoonright J$ an homeomorphism on its image as well.
Figure \ref{fig:h} tries to describe this map, which is akin to ``pinching'' $D$ around $x$.
(See e.g. \cite[Lemma 6.2]{GauldBook} for a related construction. 
Note in passing that this $h$ is non-conformal.)
As in Examples \ref{ex:BadPrufer}--\ref{ex:BadMoore}, this has the effect  
of transforming the collection of neighborhoods of the $p_i$ into a bouquet of wedges, as
seen in the magnifying glass on top right of Figure \ref{fig:caratheodory}.

\begin{example}
   \label{ex:*}
   A NH-$2$-manifold $M$ with $\mathcal{E}(M) = 2$, such that the following hold.
   \\
   $\bullet$ There is a point $b\in M$ with $NH(b)$ homeomorphic
   to \starbig. \\
   $\bullet$ Each point of $M$ has a basis of f-rim-simple Hausdorff open neighborhoods.\\
   $\bullet$ For each $n\in\omega$, $b$ has no Hausdorff neighborhood which is $n$-rim-simple.\\
   $\bullet$ There is a sorted neighborhood containing $b$.
\end{example}
\begin{proof}[Details]
Recall that $U\subset D^\circ$ and 
$\wb{U}=D$.
Set $M= \R^2\times\{0\}\cup (\R^2-U)\times\{1\}$.
Neighborhoods on the $0$-th floor and in $(\R^2-D)\times\{1\}$ are the usual ones, while those of points 
$\langle v,1\rangle$, with $v$
in
the boundary of $U$, are given by 
$$(W - U)\times\{1\}\cup \Bigl(h\circ\phi( W\cap U )\Bigr)\times\{0\},$$
where $W$ is a neighborhood of $v$ in the plane.
(This is the same as identifying $\langle u,1\rangle$ with $\langle h\circ\phi(u),0\rangle$
for each $u\in U$.)
Since $h\circ\phi:U\to D^\circ$ is a homeomorphism, this defines a $2$-manifold.
Now, any neighborhood of $b = \langle x,0\rangle$ intersects any neighborhood of any point of 
\starbig$\times\{1\}$. In particular, 
if $n$ is the number of connected components of
its intersection with a (nice enough) neighborhood of $\langle u,1\rangle$
($u$ the middle of the star \starbig), then $n$ is finite
but there is no finite bound on $n$.
Since $u$ is the only (really) problematic point, $M$ is locally f-rim-simple but not 
locally $n$-rim-simple for any $n\in\omega$.
\\
The sorted neighborhood in $M$ is the $0$-th floor, 
and actually the sorting map may be defined on all of 
$M$.
Let us define this map $f:M\to\R^2$. 
It is the projection on the first coordinate on the $0$-th floor.
Set $f(\langle v,1\rangle) = x$ whenever
$v\in\text{\starbig}$. For those $v$ in the boundary of $U$ (as a subset of $\R^2$) 
not in \starbig, hence in the unit circle except $p_1$, 
we know that $\phi$ extends homeomorphically to them, and $h\circ\phi$
is an homeomorphism from the unit circle minus $p_1$ to the unit circle minus $x$,
we may then set $h\circ\phi(p_1)= x$
and write $h\circ\phi$ in angular coordinates as $\mu(\theta)$.
Denote the polar coordinates in $\R^2$ as $\langle r,\theta\rangle_\text{pol}$
For $r\ge 1$, then 
$f( \langle \langle r,\theta\rangle_\text{pol} ,1\rangle ) = \langle 1,\mu(\theta)\rangle_\text{pol}$.
This defines a continuous function, and we are over.
\end{proof}

\begin{figure}
     \begin{center}
        \epsfig{figure=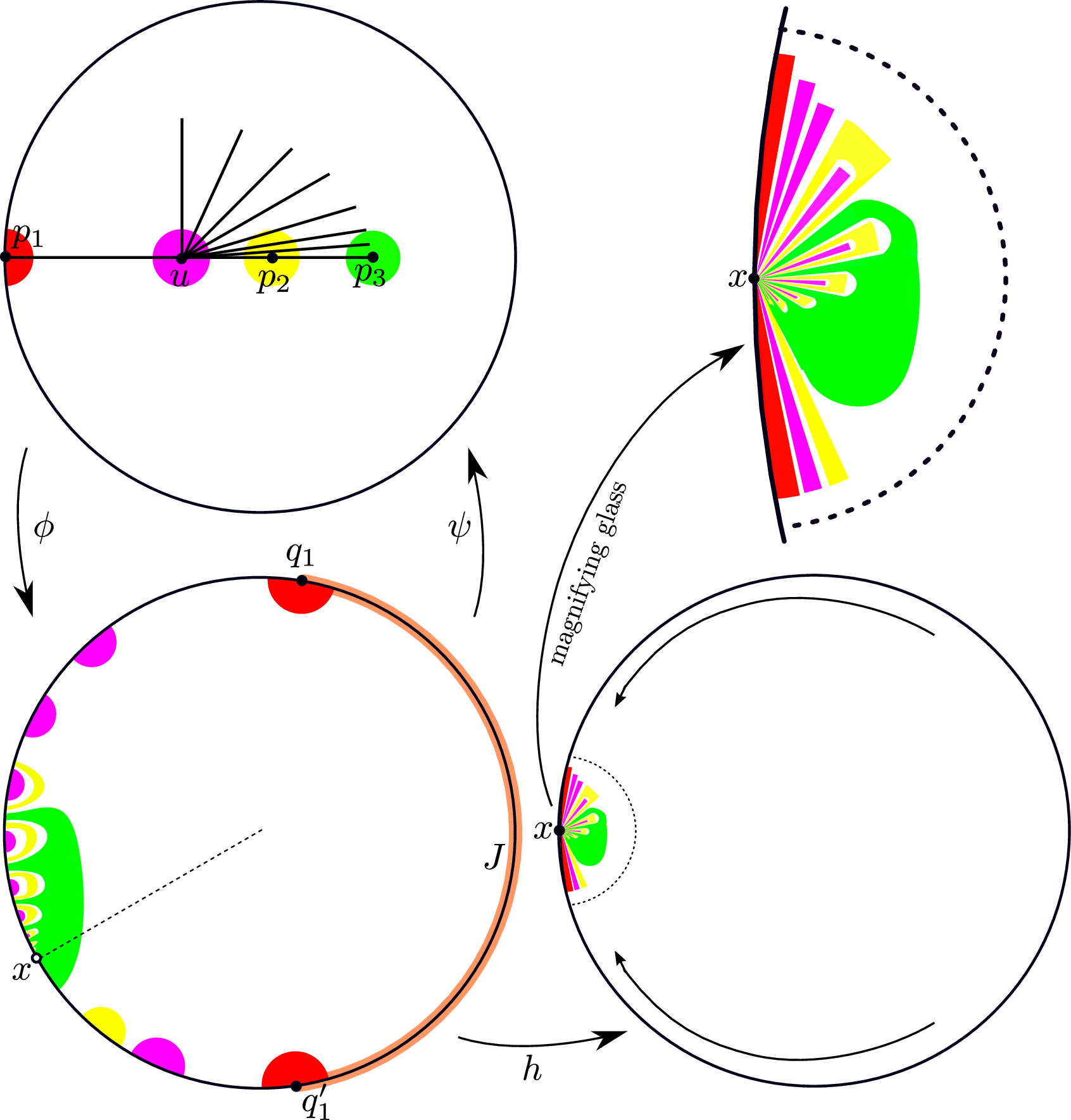, width=.75\textwidth}
         \caption{Pinching a conformal map: a domain with a non-locally connected boundary.}
         \label{fig:caratheodory2}
      \end{center}
\end{figure}

Let us now take another domain. It is shown in Figure \ref{fig:caratheodory2} (top left), to which
we now refer for what follows.
(The rays are not shown in this figure.)
It consists of the disk minus a countable union of line segments of fixed length, accumulating to one of them,
and attached to the boundary circle.
Call this set of accumulating segments \fanfig. Consider again the Riemann map
$\phi:U\to D^\circ$ and its inverse $\psi$.
By Caratheodory's theory (e.g. \cite[Theorem 15.5]{MilnorDynamics:1990}),
we see that $\psi$ extends homeomorphically to the boundary $\partial D$ minus one point, which we call 
again $x$.
(The image under $\psi$ of the ray landing at $x$ does not land on \fanfig.)
The neighborhoods of any point of \fanfig{ }(except $p_3$, that is, the extremal point)
in the segment to which the others accumulate are sent by $\phi$
to countably many disjoint open sets, infinitely many being ``engulfed'' by any neighborhood 
of the extremal point $p_3$,
as seen on the bottom left. As before, $\phi$ extends to the points of the boundary circle minus $p_1$
and its image is an open interval $J\not\ni x$;
we might thus again define a map $h:D\to D$, which is an homeomorphism on $D^\circ$ and
such that $h(J)=\partial D-\{x\}$. The effect on the neighborhoods of the various points of \fanfig{ }
are again shown on the right part of the figure.
Since $h\circ\phi:U\to D^\circ$ is again an homeomorphism, we may define another NH-$2$-manifold
exactly as in the previous example, with two floors, the difference being that this time, any 
(small enough) neighborhood of $\langle x,0\rangle$ intersected
with either $\langle u,1\rangle$ or $\langle p_2,1\rangle$ has infinitely many components.
With the same definition for the sorting map, we obtain the following.

\begin{example}
   \label{ex:fan}
   A NH-$2$-manifold $M$ with $\mathcal{E}(M) = 2$, such that the following hold.
   \\
   $\bullet$ There is a point $b\in M$ with $NH(b)$ homeomorphic
   to \fanfig. \\
   $\bullet$ No neighborhood of $b$ is f-rim-simple.\\
   $\bullet$ There is a sorted neighborhood containing $b$.
\end{example}

It seems clear\footnote{Actually, this entire last paragraph is perhaps only wishful thinking.} 
that these constructions may be done with even more complicated domains
as long as they have the property of having at least one part of their boundary
that is smooth, so that $\phi$ does extend to its points to an homeomorphism onto some interval $J\subset\partial D$.
By construction, the rest of of $\partial U$ is ``sent to $NH(b)$'', and we may end up with
rather complicated $NH(b)$ containing for instance the boundary of the Mandelbrot set or other delicacies.  
Also,
by deleting one point in the boundary circle, one may work in the plane instead of the disk.
Moreover, we may mirror the construction and ``only'' have two copies 
of \starbig and \fanfig{ } glued at their respective $p_1$ points and the circle (or line)
on the first floor, instead 
of all of $\R^2-U$. Then, 
as in Example \ref{ex:BadAll},
we may repeat the same process by doing it at each $x$ on the circle (or line),
and using concentric circles (or parallel lines),
so that each $x$ in the bottom plane has the same $NH(x)$.
None of these variations would yield an homogeneous manifold, though, because
it is not possible to send a point of $NH(b)$ to $b$ with an homeomorphism.
We let the very motivated reader fill the details (and deal with the fatal mistakes that we overlooked).


\section{HNH-manifolds and minimal flows}
\label{sec:flows}

In this section, we show how to use flows on a Hausdorff $n$-manifold $N$ to obtain a NH-$(n+1)$-manifold $M$
with $NH(x)$ a copy of a closed subspace of $N$. 
If the flow has more properties, we obtain that $NH(x)$ is a copy of the entire $N$ 
(or a discrete countable collection of copies of $N$) and/or that $M$ is homogeneous.
We start with the simplest example to illustrate the technique.

\begin{nota}
   Given $t\in S\subset\R$, $S^{*t}$ denotes $S-\{t\}$ and $\R^*$ abbreviates $\R^{*0}$. 
\end{nota}

\begin{example}
   \label{ex:HNHS1}
   A HNH-$2$-manifold $M$ such that $NH(z)$ is a discrete countable family of copies of $\mathbb{S}^1$ for each $z\in M$. 
\end{example}
\begin{proof}[Details]
   We view $\mathbb{S}^1$ as the real numbers modulo $1$.
   We let $M$ be $\mathbb{S}^1\times\R\times\Z$, topologized as follows.
   If $E\subset\mathbb{S}^1$, $n\in\Z$ and $s\in\R$, 
   define $U(E,t,\epsilon)\subset\mathbb{S}^1\times\R$ as $E\times(t-\epsilon,t+\epsilon)^{*t}$ and
   $\psi_{t,n}:\mathbb{S}^1\times\R^{*t}\to \mathbb{S}^1\times\R^{*t}$
   as 
   $$\psi_{t,n}(x,s) = \left\langle x + \frac{n}{s-t}, s \right\rangle.$$ 
   (Of course, the addition is done modulo $1$.)
   Then $\psi_{t,n}$ is an autohomeomorphism of $\mathbb{S}^1\times\R^{*t}$,
   $\psi_{t,n}\circ\psi_{t,m}=\psi_{t,n+m}$ and $\psi_{t,0}$ is the identity.
   Notice also that 
   \begin{equation}
      \label{eq:psi}
      \tag{$\psi$}
      \psi_{t_0,n}\circ\psi_{t_1,m} = \psi_{t_1,m}\circ\psi_{t_0,m}
      \text{ on }\mathbb{S}^1\times(\R-\{t_0,t_1\}).
   \end{equation}
   We define the system of neighborhoods of $\langle x,t,n\rangle$ as
   \begin{equation}
       \label{eq:W}
       \tag{W}
       W(V,t,\epsilon,n) = \Bigl(V\times\{\langle t,n\rangle\}\Bigr)
                           \cup 
                           \Bigl(\psi_{t,n}(U(V,t,\epsilon))\times\{0\}\Bigr),
   \end{equation}
   where $V$ ranges over the open subsets of $\mathbb{S}^1$ containing $x$.
   In words: 
   we start with $V\times(t-\epsilon,t+\epsilon)$,
   lift $V\times\{t\}$ at the $n$-th floor and put $V\times(t-\epsilon,t+\epsilon)^{*t}$
   on the $0$-th floor after applying $\psi_{t,n}$.
   Since $\psi_{t,0}$ is the identity, if $n=0$ we just obtain $V\times(t-\epsilon,t+\epsilon)\times\{0\}$.
   Hence, $\mathbb{S}^1\times\R\times\{0\}$ is a copy of $\mathbb{S}^1\times\R$.
   \\
   Now, given $t\in\R$, $n\in\Z$, define
   \begin{equation}
      \label{eq:h}
      \tag{h}
      h_{t,n}( \langle x,s,m \rangle ) =\left\{
                 \begin{array}{lc} 
                    \langle \psi_{t,n}(x,s),m\rangle & \text{ if }  s\not= t \\
                    \langle x,t,n+m \rangle        & \text{ if }  s = t
                 \end{array}
                 \right.     
   \end{equation}
   Then $h_{t,n}$ is an homeomorphism. To see that, notice first that 
   $h_{t,-n}\circ h_{t,n}=h_{t,n}\circ h_{t,-n}=id$, hence we just need to show that $h_{t,n}$ is continuous.
   Given any open $V\subset\mathbb{S}^1$,
   $\epsilon>0$ and $m\in\Z$, by (\ref{eq:W}) and (\ref{eq:h}),
   $h_{t,n}$ 
   turns $W(V,t,\epsilon,m)$ into $W(V,t,\epsilon,n+m)$,
   hence neighborhoods of points with second coordinate $t$ are well taken care of.
   Now, take $\langle x,s,m\rangle\in M$ with $s\not= t$ and a small neighborhood $V$ of $x$ in $\mathbb{S}^1$.
   We wish to understand what does $h_{t,n}$ do to $W(V,s,\epsilon,m)$.
   If $\epsilon<|t-s|$, its action on $V\times(s-\epsilon,s+\epsilon)\times\{0\}$ is just that
   of the homeomorphism $\psi_{t,n}$, and by (\ref{eq:psi}), 
   $h_{t,n}(\psi_{s,m}(U(V,t,\epsilon))\times\{0\})$ is equal to
   $\psi_{s,m}(\psi_{t,n}\left(U(V,t,\epsilon))\right)\times\{0\}$.
   On the $m$-th floor part of $W(V,s,\epsilon,m)$, $h_{t,n}$ just applies $\psi_{t,n}$ to $V\times\{s\}$.
   Hence, the image of $W(V,s,\epsilon,m)$ is an homeomorphic image of 
   $V\times(s-\epsilon,s+\epsilon)$ put at the $0$-th floor after applying $\psi_{s,m}$, except
   it central slide $V\times\{s\}$ which is lifted at the $m$-th floor.
   This is an open set in $M$
   (see Claim \ref{claim:Hopen} below if in need of more details). 
   It follows that the inverse $h_{t,-n}$ of $h_{t,n}$ is continuous, 
   and arguing with $h_{t,-n}$ shows that this also holds for $h_{t,n}$.
   By definition, $h_{t,n}\circ h_{t,m} = h_{t,n+m}$, hence each point of $M$ can be put on the $0$-th
   floor by an homeomorphism, which shows that $M$ is indeed a NH-manifold.\\
   Notice that for any $n\not=0$ and any $x\in\mathbb{S}^1$, the subset of $\mathbb{S}^1\times\R$
   \begin{equation}
      \label{eq:A}
      \tag{A}
       A(n,x) = \wb{\left\{\left\langle x + \frac{n}{-t},t\right\rangle\,:\,t\in\R^{*}\right\}}
                \cap \mathbb{S}^1\times\{0\}
   \end{equation}
   is equal to $\mathbb{S}^1\times\{0\}$, because it ``turns around the circle'' faster and faster as $t$ approaches $0$.
   This implies that 
   $z=\langle x,t,n\rangle$ cannot be separated from $\langle y,t,m\rangle$ whenever $n\not= m$,
   but may be separated from all the other points.
   It follows that $NH(z)$ is a countable discrete collection of copies of $\mathbb{S}^1$.\\
   To show that $M$ is homogeneous, it is enough to see that for each $t,s\in\R$ and $x,y\in\mathbb{S}^1$
   there is a homeomorphism of $M$ sending $\langle x,t,0\rangle$ to $\langle y,s,0\rangle$,
   the homeomorphisms $h_{t,n}$ enabling to travel beween distinct floors.
   For that, note first that $\langle x,t,n \rangle\mapsto\langle x+\theta,t,n\rangle$ defines an homeomorphism.
   The same holds for the map $\langle x,t,n \rangle\mapsto\langle x,t+r,n\rangle$,
   since it transforms $W(V,t,\epsilon,n)$ into $W(V,t+r,\epsilon,n)$. 
\end{proof}

\vskip .3cm
\noindent
Analyzing this example, we realized that
the features below are sufficient to make the construction work. 
\begin{itemize}
    \itemsep -.1cm
    \item[(i)]
    There is a (very simple) continuous group action of $\R$ on $\mathbb{S}^1$ given by $x\mapsto x+t$,
    such that the forward and backward orbit (i.e. taking all $t>0$ or all $t<0$) 
    of any $x$ passes through any neighborohhod of any $y\in\mathbb{S}^1$
    infinitely many times. 
    \item[(ii)]
    The group action is homogeneous (trivially in this case) 
    in the sense that for each $x,y\in\mathbb{S}^1$, there is a homeomorphism
    of $\mathbb{S}^1$ which sends the orbit of $x$ to the orbit of $y$.
\end{itemize}
The combination of (i) and the fact that $\displaystyle\frac{n}{t}$ goes to infinity when $t$ approaches $0$ imply that 
the subset defined by (\ref{eq:A}) is indeed all of $\mathbb{S}^1\times\{0\}$, which in turn yields
that $NH(z)$ is a discrete family of copies of $\mathbb{S}^1$.
The fact that the $h_{t,n}$ are indeed homeomorphisms follows from (i), in the sense that 
the commutations in (\ref{eq:psi}) hold
thanks to the fact that the group action of $\R$ on $\mathbb{S}^1$ is commutative as well.
Then, in the proof of the homogeneity of $M$, we use (ii) to show that there is an homeomorphism sending
$\langle x,t,n\rangle$ to $\langle y,t,n\rangle$. \\
The next subsection establishes a general setting for building NH-manifolds similar to Example \ref{ex:HNHS1}.
These manifolds may or may not share the same properties, depending on whether
variations of (i) and (ii) hold or not.

\subsection{Recipes for building NH-manifolds out of a flow}\label{subsec:FlowsGen}

Recall that a {\em flow} on a space $X$ is a continuous group action of $\R$ on $X$.
Flows are usually written as a map $\Phi:X\times\R\to X$, writing $\Phi_t$ for $\Phi(\cdot,t)$, such 
that $\Phi_t\circ\Phi_s = \Phi_{t+s}$ (hence, $\Phi_t$ is a homeomorphism for each $t\in\R$
and $\Phi_0$ is the identity).
The {\em orbit} $O(x)$ of a point $x\in X$ is $\{\Phi_t(x)\,:\,t\in\R\}$, if one takes only positive or negative $t$
we obtain the {\em forward} orbit $O^+(x)$ and {\em backward} orbit $O^-(x)$. 

\begin{defi}
   Let $\Phi_t$ be a flow on a Hausdorff manifold $N$.
   We say that $\Phi_t$ is transitive iff 
   there is some $x\in N$ with a dense orbit,
   and that it is homogeneous iff for each $x,y\in N$ there is
   a homeomorphism $h$ of $N$ which preserves the orbits of $\Phi_t$
   and such that $h(x)=y$.
   Finally, the flow is minimal iff the orbit of each $x\in X$ in dense in $X$.
\end{defi}

The simplest example of a non-trivial 
minimal flow on a Hausdorff surface is the irrational flow on the $2$-dimensional
torus $\mathbb{T}^2$:
let $\mathbb{T}^2$ be modeled as $[0,1]^2$ with horizontal and vertical boundaries identified, 
choose any irrational $\theta\in\R$, then $\Phi_t$ is defined by following the line of slope $\theta$
at unit speed. (It is well known that the torus is the only compact surface which admits a minimal flow).
The definition of an irrational flow obviously generalizes to the $n$-dimensional torus $\mathbb{T}^n$.

The irrational flow on $\mathbb{T}^n$ is homogeneous, as easily seen. 
On the other hand, a transitive flow may have fixed points
(see below), in which case it cannot be homogeneous (a non-fixed point cannot be sent to a fixed point
by an homeomorphism which preserves the orbits).
Notice that if a flow is transitive, the set of points with a dense orbit is dense in $N$: just take the
points in the said dense orbit.
Notice also that a flow that is both transitive and homogeneous is minimal.

\begin{defi}
   A recipe for a NH-manifold $\mathscr{R}$ (recipe for short) is a tuple $\langle G,N,\Phi_t\rangle$,
   where:
 
      -- $G$ is a finite set $\{g_0,\dots,g_{k-1}\}$ of continuous functions $\R^*\to\R$;
      
      -- $N$ is a Hausdorff $n$-manifold;
      
      -- $\Phi_t$ is a flow on $N$.
\end{defi}

Let $\mathscr{R}=\langle G,N,\Phi_t\rangle$ be such a recipe. 
For any $\sigma\in\Z^k$, we let $g_\sigma:\R^*\to\R$ be $\sum_{i=0}^{k-1}\sigma_i\cdot g_i$, where $\sigma_i$ is of 
course the $i$-th coordinate of $\sigma$. 
By definition, $g_{\sigma+\tau}=g_\sigma + g_\tau$ (where $\sigma+\tau$ is the coordinatewise addition in $\Z^k$).\\
We now define the NH-$(n+1)$-manifold $M^\mathscr{R}$ as follows. Its underlying set is $N\times\R\times\Z^k$.
We denote the origin $\langle 0,\dots,0\rangle$ in $\Z^k$ by $\langle 0\rangle$.
Given $E\subset N$, $\sigma\in\Z^k$, $s\in\R$ and $\epsilon>0$, we let 
$U(E,s,\epsilon)$ be $E\times (s-\epsilon,s+\epsilon)^{*s}$ and define 
   $\psi_{s,\sigma}:N\times\R^{*s}\to N\times\R^{*s}$
   as 
   $$ 
      \psi_{s,\sigma}(x,t) = \langle \Phi_{g_\sigma(t-s)}(x), t \rangle.
   $$ 
   Notice that since $\Phi_t$ is a flow, we again have that
   \begin{equation}
     \label{eq:Rpsi}
     \tag{$\mathscr{R}\psi$}
     \psi_{s,\sigma}\circ\psi_{t,\tau} = \psi_{t,\tau}\circ\psi_{s,\sigma}
   \end{equation}
   wherever these functions are defined
   (i.e. when the second coordinate is not $s$ or $t$).
   Moreover, 
   $\psi_{t,\langle 0\rangle}$ is the identity,
   $\psi_{t,\sigma}\circ\psi_{t,\tau}=\psi_{t,\sigma+\tau}$ and  
   $\psi_{s,\sigma}$ is an autohomeomorphism of $N\times\R^{*s}$.
   We define the system of neighborhoods of $\langle x,t,\sigma\rangle$ as in (\ref{eq:W}), that is:
   \begin{equation}
       \label{eq:RW}
       \tag{$\mathscr{R}$W}
       W(V,t,\epsilon,\sigma) = \Bigl(V\times\{\langle t,\sigma\rangle\}\Bigr) \cup 
                                \Bigl( \psi_{t,\sigma}(U(V,t,\epsilon))
       \times\{\langle 0\rangle\}\Bigr),
   \end{equation}
   where $\epsilon>0$ and $V$ ranges over the open subsets of $N$ containing $x$.
   By definition, $N\times\R\times\{\langle 0\rangle\}$ is a homeomorphic copy of $N\times\R$.
   \begin{claim}
       \label{claim:Hopen}
       Let $\vartheta:N\times(s-\epsilon,s+\epsilon)\to N\times(s-\epsilon,s+\epsilon)$
       be an homeomorphism such that $\vartheta(N\times\{t\})\subset N\times\{t\}$ for each $t\in(s-\epsilon,s+\epsilon)$.
       Then 
       $$ H = \Bigl(\psi_{s,\sigma}\circ\vartheta(U(V,s,\epsilon))\times \{\langle 0\rangle\}\Bigr)
              \cup 
              \Bigl(\vartheta( V\times\{s\})\times\{\sigma\}\Bigr)$$
       is open for each open $V\subset N$ and $\sigma\in\Z^k$
   \end{claim}
   \begin{proof}[Proof of the claim]
       If $\sigma=\langle 0\rangle$, this is immediate, hence we assume that $\sigma\not=\langle 0\rangle$.
       We only need to check that a point $\langle y,s,\sigma\rangle\in\vartheta( V\times\{s\})\times\{\sigma\}$ has an
       open neighborhood inside $H$,
       because it is immediate for points with third coordinate equal to $\langle 0\rangle$. 
       Let $V_0\subset N$ be open, $V_0\ni y$, 
       and $\delta>0$ be small enough so that 
       $\vartheta^{-1}(V_0\times(s-\delta,s+\delta))\subset V$.   
       Then $W(V_0,s,\delta,\sigma)\subset H$. 
   \end{proof}   
   Notice that for any $E\subset N$, $t\in\R$ and $\sigma\in\Z^k$,
   $E\times\{\langle t,\sigma\rangle\}$ is a copy of $E$ in $M^\mathscr{R}$.
   We define the maps $h_{t,\sigma}$ for $t\in\R$, $\sigma\in\Z^k$, exactly as before:
   \begin{equation}
      \label{eq:Rh}
      \tag{$\mathscr{R}$h}
      h_{t,\sigma}( \langle x,s,\tau \rangle ) =\left\{
                 \begin{array}{lc} 
                    \langle\psi_{t,\sigma}(x,s),\tau\rangle & \text{ if }  s\not= t \\
                    \langle x,t,\sigma+\tau \rangle        & \text{ if } s=t 
                 \end{array}
                 \right.     
   \end{equation}   

\begin{claim}
   \label{claim:M^Rmanifold}
   For each $t\in\R$, $\sigma\in\Z^k$,
   $h_{t,\sigma}$ is a homeomorphism of $M^\mathscr{R}$ and 
   $h_{t,\sigma}^{-1}=h_{t,-\sigma}$. 
   Moreover, the map $\langle x,t,\sigma\rangle\mapsto\langle x,t+r,\sigma\rangle$ is a homeomorphism
   for each $r\in\R$.
\end{claim}
\begin{proof}[Proof of the claim]
   The map $\langle x,t,\sigma\rangle\mapsto\langle x,t+r,\sigma\rangle$
   sends $W(V,t,\epsilon,\sigma)$ to $W(V,t+r,\epsilon,\sigma)$ and is thus clearly an homeomorphism.
   The fact that $h_{t,\sigma}\circ h_{t,-\sigma}=h_{t,-\sigma}\circ h_{t,\sigma}=id$
   follows from the definition and properties of $\psi_{t,\sigma}$.
   As in Example \ref{ex:NHS1}, we see that $h_{t,\sigma}$ 
   turns $W(V,t,\epsilon,\tau)$ into $W(V,t,\epsilon,\tau+\sigma)$, hence acts homeomorphically on
   $N\times\{t\}\times\Z^k$.
   If $t\not=s$, take $\epsilon<|t-s|$.
   By (\ref{eq:Rh}), (\ref{eq:RW}) and (\ref{eq:Rpsi}), we see that
   \begin{align*}
      h_{t,\sigma}(W(V,s,\epsilon,\tau)) &= \psi_{t,\sigma}\circ\psi_{s,\tau}(U(V,s,\epsilon))\times\{\langle 0\rangle\}
                                            \cup \psi_{t,\sigma}(V\times\{s\})\times\{\tau\} \\
                                         &= \psi_{s,\tau}\circ\psi_{t,\sigma}(U(V,s,\epsilon))\times\{\langle 0\rangle\}
                                            \cup \psi_{t,\sigma}(V\times\{s\})\times\{\tau\}
   \end{align*}
   which is open by Claim \ref{claim:Hopen}.
   If follows that $h_{t,-\sigma}$ is an homeomorphism, hence so is $h_{t,\sigma}$.
\end{proof}
Notice that the previous claim implies that $M^\mathscr{R}$ is a NH-manifold.
The equivalent of (\ref{eq:A}) in this setting is a bit more cumbersome to write down.
Given $V$ open in $N$ and $\sigma\in\Z^k$, set $B(V,\sigma)$ to be the closure in $N\times\R$
of $\{\langle \Phi_{g_\sigma(t)}(y),t\rangle\,:\,y\in V,\,t\in\R^*\}$,
and $A(V,\sigma)$ be the projection on the first coordinate of $B(V,\sigma)\cap N\times\{0\}$.
Finally, set $A(x,\sigma)$ be the intersection of $A(V,\sigma)$ for all open $V\ni x$.
\begin{claim}
   $NH(\langle x,t,\sigma\rangle)$ is the discrete union of the $A(x,\sigma-\tau)\times\{\langle t,\tau\rangle\}$
   for each $\tau\not=\sigma$.
\end{claim}
\begin{proof}[Proof of the claim]
   Fix $\sigma\not=\tau$ in $\Z^k$ and $\epsilon>0$.
   Let $x,y\in N$ and choose open $V_x\ni x$, $V_y\ni y$.
   Consider the intersection $W(V_x,t,\epsilon,\sigma)\cap W(V_y,t,\epsilon,\tau)$.
   By (\ref{eq:RW}), these two sets may intersect only on the $\langle 0\rangle$-th floor, where the intersection
   is equal to
   $\psi_{t,\sigma}(U(V_x,t,\epsilon))\cap \psi_{t,\tau}(U(V_y,t,\epsilon))$.
   Applying the homeomorphism $\psi_{t,-\tau}$ to this intersection yields
   $\psi_{t,\sigma-\tau}(U(V_x,t,\epsilon))\cap U(V_y,t,\epsilon)$.
   By definition, this intersection is non empty for all $\epsilon>0$ iff $y\in A(x,\sigma-\tau)$,
   which shows that 
   $$ NH(\langle x,t,\sigma\rangle)\cap N\times\{\langle t,\tau\rangle\}
      =A(x,\sigma-\tau)\times\{\langle t,\tau\rangle\}.$$
   If $x\not=y$, choose disjoint open $V_x\ni x$ and $V_y\ni y$, then for each $\epsilon>0$ we have that
   $W(x,t,\epsilon,\sigma)\cap W(y,t,\epsilon,\sigma)=\varnothing$, hence
   $\langle x,t,\sigma\rangle$ can be separated from $\langle y,t,\sigma\rangle$.
   If $t\not=s$, choose $\epsilon<|t-s|$, then 
   $W(N,t,\epsilon,\sigma)\cap W(N,s,\epsilon,\tau)=\varnothing$ for any $\tau$.
   Thus, $\langle x,t,\sigma\rangle$ can be separated from any point with second coordinate $s\not= t$.
\end{proof}

\begin{claim}
   If $\Phi_t$ is transitive,
   $\displaystyle\lim_{x\to 0^+}g_\sigma(x) = -\lim_{x\to 0^-}g_\sigma(x)$ and 
   $\displaystyle\lim_{x\to 0^+}g_\sigma(x)=\pm\infty$, 
   then $A(x,\sigma) = N$ for each $x\in N$.
\end{claim}
\begin{proof}[Proof of the claim]
   We may assume that $\displaystyle\lim_{x\to 0^+}g_\sigma(x) = +\infty$.
   Given any open $V\ni x$, there is some $y\in V$ such that either $O^+(y)$ or $O^-(y)$ is dense in $N$
   (because there is a point of $V$ in the dense orbit).
   By definition, $A(V,\sigma)$ is the accumulation points of the graph of the function $t\mapsto \Phi_{g_\sigma(t)}(y)$ 
   for all $y\in V$, when $t$ goes to $0$.
   Since $g_\sigma$ has opposite sign on each side of the origin (if near enough), both forward and
   backward orbits of $y$ accumulate. It follows that $A(V,\sigma) = N$.
\end{proof}

\begin{claim}
   If $\Phi_t$ is an homogeneous flow, then $M^\mathscr{R}$ is homogeneous.
\end{claim}
\begin{proof}[Proof of the claim]
   Claim \ref{claim:M^Rmanifold} shows that there are homeomorphisms 
   $\langle x,t,\sigma\rangle\to\langle x,t,\tau\rangle$ and $\langle x,t,\sigma\rangle\to\langle x,s,\sigma\rangle$
   for each $x\in N$, $s,t\in\R$ and $\sigma,\tau\in\Z^k$.
   It is thus enough to exhibit homeomorphisms
   sending $\langle x,t,\sigma\rangle$ to $\langle y,t,\sigma\rangle$ for each $x,y\in N$.
   But this follows from the flow's homogeneity: given an homeomorphism $h:N\to N$
   preserving the orbits, 
   then $\langle x,t,\sigma\rangle\mapsto \langle h(x),t,\sigma\rangle$ 
   is a homeomorphism of $M^\mathscr{R}$.
\end{proof}

At the risk of repeating ourselves once more, we
summarize all of this in a Lemma.
Only (c) is not proved above, but is an immediate consequence of (a) since the $\langle 0\rangle$-th floor is 
H- and CH-maximal.
\begin{lemma}[{\bf Properties of $M^\mathscr{R}$}]
   \label{lemma:propMR}
   \ \\
   Let $\mathscr{R}=\langle G,N,\Phi_t\rangle$ be a recipe, with $N$ a Hausdorff $n$-manifold.
   Then $M^\mathscr{R}$ is an NH-$(n+1)$-manifold with underlying set $N\times\R\times\Z^k$,
   such that:
   \begin{itemize}
   \itemsep -.1cm
   \item[(a)] for all $x\in N$, $t,s\in\R$ and $\sigma,\tau\in\Z^k$, there is a homeomorphism
      sending $\langle x,t,\sigma\rangle$ to $\langle x,s,\tau\rangle$;
   \item[(b)] $NH(\langle x,t,\sigma\rangle)\cap 
               \Bigl(N\times\R\times\{\tau\}\Bigr)$ is equal to 
       $A(x,\sigma-\tau)\times\{\langle t,\tau\rangle\}$
       when $\tau\not=\sigma$ (and is empty when $\tau=\sigma$); moreover, the family of these components is
       discrete;
   \item[(c)] each $x\in M^\mathscr{R}$ has a H- and CH-maximal open neighborhood homeomorphic to $N\times\R$.
   \end{itemize}
   \noindent
   Moreover, the following hold.
   \begin{itemize}
   \itemsep -.1cm
   \item[(d)]
   If $\Phi_t$ is transitive,
   $\displaystyle\lim_{x\to 0^+}g_\sigma(x) = -\lim_{x\to 0^-}g_\sigma(x)$ and 
   $\displaystyle\lim_{x\to 0^+}g_\sigma(x)=\pm\infty$,
   then $A(x,\sigma)=N$ for each $x$;
   \item[(e)]
   if $\Phi_t$ is homogeneous, then $M^\mathscr{R}$ is homogeneous.
   \end{itemize}
\end{lemma}


\begin{rem}
   The choice of $\Z^k$ as a third coordinate is somewhat arbitrary, as we could choose 
   any additive subgroups $P_0,\dots,P_k$ of $\R$, take the product group $P=\Pi_{i=0}^k P_k\subset\R^k$
   and define $M^\mathscr{R}$ with underlying set $N\times\R\times P$.
   The only changes would be in point (b) above. 
\end{rem}

\subsection{Examples}
\label{subsec:FlowsEx}

In order to shorten a bit some statements, we say that a manifold $M$ is {\em $X$-like} iff
any $x\in M$ has a H- and CH-maximal neighborhood homeomorphic to $X$.
All the claimed properties of the examples below follow from Lemma \ref{lemma:propMR}.
The next one was the first one we discovered.

\begin{example}
   \label{ex:sines}
   An $\R^2$-like HNH-manifold $M^\mathscr{R}$ 
   such that $NH(z)$ is a countable discrete family of copies of $[0,1]$
   for each $z\in M^\mathscr{R}$, so $NH(z)$ is non-homogeneous.
\end{example}
\begin{proof}[Details]
   Take $G=\{g_0\}$, where $g_0(x) = \sin(1/x)$, $N = \R$ and $\Phi_t(x) = x+t$.
   This flow is obviously minimal and homogeneous, hence $M^\mathscr{R}$ is homogeneous.
   Since $G$ has only one element, we write $n$ instead of $\langle n\rangle$.
   If $n\ge 0$,
   then $A(x,\pm n) = [x-n,x+n]$  for each $x$.
   It follows that
   $NH(\langle x,t,n\rangle)$ 
   is a discrete countable family of compact intervals.
\end{proof}

By changing $g_0$ we obtain the following.
\begin{example}
   \label{ex:sines2}
   An $\R^2$-like HNH-manifold $M^\mathscr{R}$  
   such that $NH(z)$ is a countable discrete family of copies of $\R$
   for each $z\in M^\mathscr{R}$.
\end{example}
\begin{proof}[Details]
   We take again $N=\R$, $\Phi_t(x) = x+t$, but set
   $G=\{g_1\}$, where $\displaystyle g_1(x) = \frac{\sin(1/x)}{x}$.
   Thus, $A(x,n)=\R$ for each $n\not= 0$.
\end{proof}

Combining Examples \ref{ex:sines}--\ref{ex:sines2} yields even worse. 
\begin{example}
   \label{ex:sines3}
   An $\R^2$-like HNH-manifold $M^\mathscr{R}$ 
   such that $NH(z)$ is a countable discrete family of copies of $\R$, $[0,1]$ and singletons
   for each $z\in M^\mathscr{R}$. In particular, connected components of $NH(z)$ are not all homeomorphic.
\end{example}
\begin{proof}[Details]
   Same $N$ and $\Phi_t$ as in the previous examples, 
   but take $G=\{g_0,g_1,g_2\}$, where $g_0,g_1$ are as above and $g_2$ is constant on $0$.
   This time, $A(x,\langle m,n,k\rangle) =\R$ when $n\not= 0$, $A(x,\langle \pm m,0,k\rangle) =[x-m,x+m]$
   and $A(x,\langle 0,0,k\rangle) = \{x\}$. 
\end{proof}

This example shows the usefulness of considering a family of functions $G$ rather that only one.
Let us turn to higher dimensions. Recall that $\mathbb{T}^n$ denotes the $n$-dimensional
torus.

\begin{example}
   There is a
   $(\mathbb{T}^n\times\R)$-like
   HNH-manifold such that $NH(z)$ is a discrete collection of countably many copies of  
   $\mathbb{T}^n$ for each point $z$. 
\end{example}
\begin{proof}[Details]
   The irrational flow $\Phi_t$
   on $\mathbb{T}^n$ is homogeneous and minimal, we just need to set $g(t) = 1/t$ and $G=\{g\}$.
\end{proof}

\begin{example}
   There is a
   $(\mathbb{T}^n\times\R)$-like
   HNH-manifold such that $NH(z)$ is a countable discrete collection of subspaces, each homeomorphic to 
   either $\mathbb{T}^n$, or $[0,1]$, or a point.
\end{example}
Hence, not only the connected components of $NH(x)$ are not homeomorphic, they also have distinct dimensions
(in the naive sense).
\begin{proof}[Details]
   Take a minimal homogeneous flow on $\mathbb{T}^n$, set $g_0(t) = 1/t$, $g_1(t) = \sin(1/t)$, $g_2(t) = 0$ (constant) 
   and et $G=\{g_0,g_1,g_2\}$. Then $A(x,\langle m_0,m_1,m_2\rangle)$ is equal to $\mathbb{T}^n$ whenever $m_0\not=0$,
   to $\{\Phi_{t}(x)\,:\,|t|\le|m_1|\}$ 
   (which is a copy of a compact interval since $\Phi_t$ is a flow such $\Phi_t(x)\not=\Phi_s(x)$ when $s\not= t$)
   when $m_0=0,m_1\not=0$, and to the singleton $\{x\}$ if
   $m_0=m_1=0$.  
\end{proof}

We may also obtain non-Hausdorff and non-homogeneous $NH(z)$ in an HNH-manifold.
\begin{example}
   A HNH-manifold such that $NH(z)$ is non-homogeneous and non-Hausdorff for each point $z$.
\end{example}
\begin{proof}[Details]
   Take any $M^\mathscr{R}$ of the previous examples
   and let $N$ be $\mathbf{G}(M^\mathscr{R},2)$, then $N$ is a HNH-manifold by Lemma \ref{lemma:THNH}.
   Then one sees easily that $NH(z)$ is non-Hausdorff.
\end{proof}

Notice that in every example above, $\Phi_t$ is homogeneous (and minimal).
If one gives up the flow's homogeneity but keeps its transitivity, 
then the only property lost for $M^\mathscr{R}$ is also homogeneity.
In that case, keeping all the floors of $M^\mathscr{R}$ seems a bit useless (since these floors are there 
precisely to ensure homogenity). 
If $0\le\ell\le k-1$, we let $\langle 1\rangle_\ell\in\Z^k$ denote $\langle 0,\dots,0,1,0,\dots,0\rangle$
where $1$ is at the $\ell$'s index. 
We let now  
$M^\mathscr{R}_\ell$ be the ``zeroth and first floors in the $\ell$-th coordinate'' of $M^\mathscr{R}$,
that is, the submanifold $M^\mathscr{R}\cap N\times\R\times \{\langle 0\rangle,\langle 1\rangle_\ell\}$.
Then $M^\mathscr{R}_\ell$ is not homogeneous, but is still $(N\times\R)$-like,
and (b) of Lemma \ref{lemma:propMR} still holds.
Since there are only two ``floors'', $NH(x)$ is contained only in one.
\\
We may thus look for manifolds which admit transitive flows. 
To say that our knowledge of this mathematical field
is lacunary is an understatement, so we content ourselves with well known results
and let the more educated readers find other examples.
We will use the following result of V. Jim\'enez L\'opez and G. Soler L\'opez
\cite[Theorems A--B]{JimenezLopezSolerLopez:2004}.
\begin{thm}\ \\
   \label{thm:transitiveflows}
   (a)
   All metrizable surfaces which are not
   homeomorphic to the sphere, the projective plane or an open subset of the Klein bottle
   admit a transitive flow.\\
   (b) Any metrizable manifold of dimension $\ge 3$ admits a transitive flow.
\end{thm}

These yield the following.
\begin{example}
   Let $N$ be either a surface as in (a) of Theorem \ref{thm:transitiveflows}
   or a Hausdorff manifold of dimension $\ge 3$. 
   Then there is an $(N\times\R)$-like manifold 
   such that $NH(z)$ is a copy of $N$ for each $z\in M$.
\end{example}
\begin{proof}[Details]
   Take the transitive flow given above, and set again $g(t) = 1/t$ and $G=\{g\}$.
   Then $M^\mathscr{R}_0$ has the required properties.
\end{proof}

Let us end this section with a couple of remarks and questions.
In our constructions, flow's homogeneity is definitely needed to ensure $M^\mathscr{R}$'s homogeneity, which
is moreover plugged in $M^\mathscr{R}$ through the use of the distinct floors $\langle\sigma\rangle$.
Moreover, we need at least the flow's transitivity to obtain that $NH(x)$ is a copy 
(or a countable discrete collection of copies) of the manifold
we started with. But
the mere existence of non-trivial flows on non-metrizable Hausdorff manifolds
(for instance, those with a copy of $\omega_1$ inside)
is not automatic, see e.g. \cite{GGDynamics}. 
We also do not possess any non-trivial criteria for {\em a priori} deciding that 
a (copy of a) given Hausdorff manifold $N$ cannot appear as $NH(x)$ for some $x$ in an HNH-manifold 
of strictly bigger dimension.
It is thus probable that new ideas 
(or a maybe just a clearer mind) are necessary in order to answer (in the positive
or the negative) any of the following questions.

\begin{ques}
   Is there an HNH-$2$-manifold $M$ such that, for any $z\in M$, $NH(z)$ is homeomorphic to
   a single copy of $\mathbb{S}^1$~? a single copy of $\R$~?
\end{ques} 

\begin{ques}
   Is there an HNH-$3$-manifold $M$ such that, for any $z\in M$, $NH(z)$ is homeomorphic to:\\ 
   (a) a single copy of $\mathbb{T}^2$~?\\
   (b) a single copy or a (countable) discrete family of copies of $\Sigma_{g,0}$
       (the compact surface of genus $g$), when $g\ge 2$~?\\
   (c) a discrete family of copies of $\mathbb{S}^2$~? 
\end{ques} 

\begin{ques}
   Is there an HNH-manifold $M$ 
   such that for each $z\in M$, $NH(z)$ contains a copy of $\omega_1$~? a copy of $\LL_+$~?  
\end{ques}

Note in passing that none of our examples have a weakly sorted neighborhood 
(e.g. by Proposition \ref{prop:homogensort}) 
or a base of rim-simple neighborhoods (by construction).


\section{Some more questions}\label{sec:ques}

Besides those already encontered in this note,
we list here questions that seem natural to us but that we did not investigate
(or for which we made no progress at all).

\vskip .3cm
\noindent
First, any attempt at a general classification of HNH-manifolds seems frivolous, 
but maybe something can be done in dimension $1$.

\begin{ques}
   Can we somehow classify the HNH-$1$-manifolds for which each point is contained in:\\
   (a) A sorted neighborhood~? An H-maximal sorted neighborhood~?\\
   (b) An open subset which is both CH-maximal and H-maximal~?\\
   (c) A sorted neighborhood which is both CH-maximal and H-maximal~?
\end{ques}

Recall that Examples \ref{ex:hersep} and \ref{ex:nocloseddis} satisfy (c) (but are not homogeneous).
Hence, any classification in (c) depends on the answers to Question \ref{q:homkun}.
Note also that Example \ref{ex:FatS1} satisfies (b).\\
We saw in Examples \ref{ex:NHclnd}, \ref{ex:NH-psi2} and \ref{ex:NHspecialtree} that 
there are $1$-manifolds such $NH(x)$ contains the Cantor space, $\Psi$-spaces
or a copy of special tree for some point $x$.

\begin{ques}
   Are there NH-$1$-manifolds
   with some $x$ such that $NH(x)$ contains a copy of a Suslin tree~? Of $\omega_1$~? 
\end{ques}

If $M$ is an ENH-manifold
then $|\wb{U}-U|\ge\text{cov}(\mathcal{B})$ for any open $U\subset M$ by Corollary \ref{cor:ClH}.
Since there are models of set theory with $\text{cov}(\mathcal{B})$ strictly smaller than the continuum, the
following 
(somewhat related to Question \ref{ques:BumpeqOHc}) is not ruled out a priori.
\begin{ques} 
   Are there consistent examples of ENH- or HNH-manifolds with some open subset $U$ such that
   $|\wb{U}-U|<\mathfrak{c}$~?  
\end{ques}

A completely different type of question, more akin to the interest of physicists:

\begin{ques}
   Do our examples with ``complicated'' $NH(x)$, or the homogeneous ones, admit a differentiable structure~? 
\end{ques}

And we finish with the question asked countless times at the end of countless talks and PhD defenses:

\begin{ques}
   What about higher dimensions~?
\end{ques}

That is: can we say more than what is shown in Subsection \ref{subsec:FlowsEx}~?


\end{document}